%% file: spectral.tex
\documentclass[11pt,final]{iuthesis}

\usepackage{amssymb}
\usepackage{accents}
\usepackage{epic}
\usepackage{verbatim}
\usepackage{paralist} 
\usepackage{url}
\usepackage[numbers]{natbib}

\title{Noncomputable Spectral Sets}
\author{Jason Teutsch}
\department{Mathematics}
\advisor{Lance Fortnow}
\secondreader{David McCarty}
\thirdreader{Mike Dunn}
\fourthreader{Janos Simon}
\fifthreader{Kevin Zumbrun}
\department{Mathematics}
\submitdate{January 19, 2007}
\copyrightyear{2007}
\makeindex

\newcommand{\INF}{\mathrm{INF}}
\newcommand{\TOT}{\mathrm{TOT}}

\newcommand{\FUN}{\mathrm{FUN}}
\newcommand{\FIN}{\mathrm{FIN}}
\newcommand{\COF}{\mathrm{COF}}
\newcommand{\mCOMP}{\mathrm{mCOMP}}
\newcommand{\LOW}{\mathrm{LOW}}
\newcommand{\HIGH}{\mathrm{HIGH}}

\newcommand{\K}{\mathit{K}}
\newcommand{\fMIN}{\mathrm{f\text{-}MIN}}
\newcommand{\MIN}{\mathrm{MIN}}
\newcommand{\RAND}{\mathrm{RAND}}
\newcommand{\R}{\mathrm{R}}
\newcommand{\fR}{\mathrm{fR}}
\newcommand{\sR}{\mathrm{sR}}
\newcommand{\mn}{\mathrm{min}}
\newcommand{\FUNfR}{\FUN\text{-}\fR}

\newcommand{\T}{\mathrm{T}}
\newcommand{\m}{\mathrm{m}}
\newcommand{\btt}{\mathrm{btt}}
\newcommand{\Tt}{\mathrm{tt}}
\newcommand{\bT}{\mathrm{bT}}
\newcommand{\tm}{\mathrm{time}}
\newcommand{\p}{\mathrm{p}}

\newcommand{\converge}{\mathop\downarrow}
\newcommand{\diverge}{\mathop\uparrow}

\newcommand{\ijoin}{\mathop{\oplus}}
\newcommand{\join}{\mathrel{\oplus}}
\newcommand{\phe}{\varphi}
\newcommand{\zero}{\emptyset}
\newcommand{\0}{\mathbf{0}}
\DeclareMathOperator{\dom}{dom}
\DeclareMathOperator{\rng}{range}
\newcommand{\thick}{\mathrm{Thick\text{-}}}
\newcommand{\aev}{\mathrm{a.e.\text{-}}}
\newcommand{\omegaINF}{\omega\mathrm{\text{-}}\INF}
\newcommand{\REA}{\mathrm{REA}}

\newcommand{\restr}{\mathrel{\upharpoonright\nolinebreak[4]\hspace{-0.65ex}\upharpoonright}}
\newcommand{\restrs}{\mathrel{\upharpoonright\nolinebreak[4]\hspace{-0.50ex}\upharpoonright}}

\newcommand{\Tincomp}{\mathrel{\vert_\T}}
\newcommand{\incomp}[1]{\mathrel{\vert_{#1}}}
\newcommand{\orr}{\mathrel{\vee}}
\newcommand{\andd}{\quad \& \quad}
\newcommand{\intersect}{\mathrel{\cap}}
\newcommand{\union}{\mathrel{\cup}}
\newcommand{\Union}{\mathop{\bigcup}}
\newcommand{\satisfies}{\mathrel{\models}}
\newcommand{\of}{\mathrel{\circ}}
\newcommand{\bigiff}{\enskip\iff\enskip}

\newcommand{\lcode}{\left\langle}
\newcommand{\rcode}{\right\rangle}
\newcommand{\pair}[1]{{\lcode#1\rcode}}
\newcommand{\size}[1]{{\left|#1\right|}}
\newcommand{\paren}[1]{{\left(#1\right)}}
\newcommand{\compliment}[1]{\overline{#1}}
\newcommand{\abs}[1]{\lvert#1\rvert}

\newcommand{\cons}{\pagebreak[2]\noindent\textit{Construction.}\enskip}

\newenvironment{thmenum}{\begin{enumerate}[\scshape (i)]}{\end{enumerate}}
\newenvironment{defenum}{\begin{enumerate}[\scshape (i)]}{\end{enumerate}}
\newenvironment{step}{\begin{enumerate}[\itshape Step 1.]}{\end{enumerate}}
\newcommand{\nlb}{\nolinebreak[3]}

\theoremstyle{plain}        \newtheorem{thm}{Theorem}[section]
\theoremstyle{definition}   \newtheorem{defn}[thm]{Definition}
\theoremstyle{plain}        \newtheorem{prop}[thm]{Proposition}
\theoremstyle{plain}        \newtheorem{conj}[thm]{Conjecture}
\theoremstyle{plain}        \newtheorem{cor}[thm]{Corollary}
\theoremstyle{plain}        \newtheorem{lemma}[thm]{Lemma}
\theoremstyle{remark}       
\theoremstyle{plain}        \newtheorem{claim}[thm]{Claim}
\theoremstyle{remark}       \newtheorem*{rem}{Remark}
\theoremstyle{definition}   \newtheorem*{nota}{Notation}

\numberwithin{equation}{chapter}
\numberwithin{thm}{section}
\numberwithin{figure}{chapter}

\begin{document}

\begin{dedication}
For my Mama, whose *-minimal index is computable (because it's finite).

\vspace{20ex}

Love the Ma.
\end{dedication}

\begin{acknowledgements}
It seems nothing short of a miracle that I now find myself in a position to
graduate with a doctorate in mathematics.  This thesis is the product of
contributions from a few outstanding individuals that I've known (in the
biblical sense) along my academic path.  Without the personal connections that
I am about to mention, there would indeed be no work here to acknowledge.

In the summer of 2004, I showed up on Lance Fortnow's doorstep in Chicago,
poised to learn something about computer science. Not only did Lance \emph{not}
send me home (which would have been impossible anyway, since I didn't come from
anywhere), he completely encouraged me to learn.  Lance dedicated enormous time
and energy to helping me achieve excellence in my research.  Not only did he
provide me with invaluable guidance and fill me with enthusiasm for my academic
work, but in a formal sense, he allowed me to continue my studies by taking me
on as a PhD student.  When people read this PhD dissertation a hundred years
from now, Lance Fortnow will undoubtedly be remembered for his outstanding
research contributions to the field of theoretical computer science.  But in
2107, this thesis will also remind readers that Lance is an exceptional mentor
who showed exceeding kindness towards strangers. His efforts have made all the
difference. Thanks for taking a chance on me.

I would like to thank the faculty, staff, and administration at University of
Chicago and its Computer Science Department, as well as the affiliates of the
Committee on Institutional Cooperation for providing me with all of the
resources I needed to complete a research dissertation in mathematics during
these past two years.  Not only have these institutions provided me with access
to an exceptional educational opportunities in theoretical computer science,
but they also helped me to formally complete the degree requirements at my home
institution. Even though I was never a degree student at University of Chicago,
I've always felt like family.

Thanks to my mom and dad, who provided me with support when my home institution
wasn't able to.  In 2004, my parents gave me a place to live, a lot of tasty
food, and a bit of encouragement when I really needed it.

Thanks to Professor William Ziemer and Chunlai Zhou, who together taught me the
love of mathematical analysis early on in my graduate career.  Chunlai, in
particular, has been a most insightful study companion.  His contagious
commitment to learning surely contributed to my passing grades on the math
qualifying exams during my first year of graduate school.

Inspiration for this project came from Marcus Schaefer, Robert Soare, and Joe
Miller.  Marcus Schaefer invented the tantalizing $\MIN^*$ problem, the
catalyst for this entire work.  Robert Soare not only invited me to his
enlightening computability theory class this year, but he also shared some
technical suggestions which gave me confidence to prove a couple of the hard
theorems in this thesis.  In 2003, I first experienced computability theory
during Joe Miller's office hours for Recursion Theory.  Although I had a
scheduling conflict which prevented me from ever attending his lectures, I
discovered the beautiful nature of computability theory outside of class time,
thanks to the instructor's patience and enthusiasm.

Thanks to Kevin Zumbrun and Misty Cummings at Indiana University, whose efforts
ultimately assembled the bureaucratic pieces of my academic career and allowed
me to earn a doctorate in mathematics.  For his contributions as a mentor, I'd
also like to thank Larry Moss, who would have been on my dissertation committee
had he not moved to Africa this year.

Finally, I would like to acknowledge my colleagues at University of Chicago,
especially Parinya Chalermsook and Raghav Kulkarni, who introduced me to the
field of computer science.  And Gabriela Turcu, who went ice skating with me
sometimes.

\bigskip

Goodnight, and g-d bless America.
\end{acknowledgements}

\begin{abstract}
It is a basic fact that, given a computer language and a computable integer
function, there exists a shortest program in that language which computes the
desired function. Once a programmer establishes the correctness of a program,
she then need only verify that the program is the shortest possible in order to
declare complete victory. Unfortunately, she can't.  A creature that could
identify minimal programs would not only be able to decide the halting problem,
but she could even decide the halting problem for machines with halting set
oracles. Such a creature exceeds the powers of ordinary machine cognition, and
must therefore be a divine jument-numen.

Suppose, however, the programmer would be satisfied to know just whether or not
her program is minimal up to finitely many errors.  In this case, even the
jument-numen is helpless: something much stronger is needed.  Just as we can
associate equality with the ordinary notion of ``minimal,'' we can associate an
equivalence relation, $=^*$, with the principle of ``minimal up to finitely
many errors.''  This thesis is the first to explore the extensive realm of
minimal indices beyond the $=^*$ relation.  Every equivalence relation gives
rise to a notion of minimality, modulo that relation.  We call the resulting
collection of minimal indices a \emph{spectral set}, because it selects exactly
one function from each equivalence class.  Spectral sets are rare, natural
examples of non-index sets which are neither $\Sigma_n$ nor $\Pi_n$-complete.

In this thesis, we classify spectral sets according to their thinness and
information content.  We give optimal immunity results for the spectral sets
(that is, we identify types of sets which are not contained in spectral sets),
and we place these sets in the arithmetic hierarchy (which quantitatively
measures their information contents).  Some lower bounds in the arithmetic
hierarchy follow from immunity properties alone, but we further refine these
immunity bounds using direct techniques.  We also measure information content
with Turing equivalences. In fact, the spectral sets become canonical
iterations of the halting set when we list our computer programs in the right
order. Regardless of the particular numbering, a reasonable amount of
information is always present in such sets.

We now informally describe the contents of some spectral sets.  The
$\Pi_3$-Separation Theorem says that the spectral sets pertaining to
$\equiv_1$, $=^*$, and $\equiv_\m$ each have the same complexity with respect
to the arithmetic hierarchy, yet each of these sets is immune against a
different level of the arithmetic hierarchy.  We can thus quantitatively
compare the strength of equivalence relations by way of immunity.  We also
prove a result which we call the Forcing Lowness Lemma.  Using this lemma, we
show that $\zero''''$ is decidable in $\MIN^\T$ (the spectral set for
$\equiv_\T$) together with $\zero''$.  This result is probably optimal, and we
apply the Forcing Lowness Lemma again to show that, in some formal sense, this
fact will be difficult to prove.

Armed with this new machinery, we highlight its utility with some new
applications.  First, we prove the Peak Hierarchy Theorem: there exists a set
which neither contains nor is disjoint from any infinite arithmetic set, yet
the set is majorized by a computable function.  Furthermore, the set that we
construct is natural in the sense that it contains a spectral set.  Along the
way, we construct a computable sequence of c.e$.$ sets in which no set can be
computed from the join of the others, for any iteration of the jump operator.

We use the machinery of spectral sets to quantitatively compare the power of
nondeterminism with the power of the jump operator.  We show that in the world
of computably enumerable sets, nondeterminism contributes nothing to immunity.
In this respect, the jump operator outshines nondeterminism.  Nonetheless, we
can ascend naturally from the lowest level of the spectral hierarchy using
nondeterminism.

Finally, we present connections to the Arslanov Completeness Criterion which
stand as immediate consequences of immunity properties for spectral sets.
\end{abstract}

\frontmatter

\maketitle

\signaturepage
\copyrightpage
\makededication
\makeack
\makeabstract

\tableofcontents

\mainmatter

\include{episode_IV}
\include{turing}
\include{immunity}
\include{thickville}
\include{kolmogorov}
\include{peak_hierarchy}

\appendix

\include{open}

\bibliographystyle{plainnat}
\bibliography{spectral}
\printindex
\backmatter

\end{document}

%% file: episode_IV.tex
\chapter{Introduction}

\section{Episode IV: A New Hierarchy}

It all begins with Occam's razor.  Since the fourteenth century, and certainly
not before then, mankind has established a universal preference for simplicity.
For centuries, people have enjoyed the pleasures of short descriptions and the
joys of simple explanations.  In the twentieth century, this predisposition
abundantly manifests itself in computer science: humans love short computer
programs.

We now consider two types of people who are especially keen (resp$.$ not keen)
on short programs.  Note that the length of the shortest program describing an
algorithm is a measure of its complexity.  An output which follows a simple,
constructive pattern can not be seen as random.  In particular, outputs with
short descriptions are not \emph{Kolmogorov random}.  Therefore, we expect that
fans of Kolmogorov random strings will not like shortest programs (unless the
programs are really long). In the tradition of minimal indices, on the other
hand, machine learnists are generally dissatisfied with machines that merely
learn the index of a target function \citep{ACJS04}.  They prefer indices which
are close to minimal. Machine learning thus provides a practical application
for the theory of shortest programs.

Minimal indices, or shortest programs, are generalizations of Kolmogorov random
strings.  According to Kolmogorov complexity, strings which lack short
descriptions contain more information than those that have.  This point-of-view
matches our intuition: a string which is truly random does not follow a simple
pattern and requires many bits to describe it.  The index of a shortest program
is always a Kolmogorov random string because if it were not, then some smaller
(a.k.a shorter) index would compute that program.  The set $\MIN$
(Definition~\ref{defnminsets}) thus contains the shortest descriptions, or
``random strings'' for c.e$.$ sets. Later we shall consider further
generalizations, such as $\MIN^\T$, the shortest descriptions for c.e$.$ Turing
degrees.

If simplicity constitutes our objective, then we have no better place to begin
our study of computability theory than with the following eloquent observation:
\begin{quote}
\textit{The set of shortest programs is not computable.}
\end{quote}
By the ``set of shortest programs,'' we mean the set succinctly characterized
by
\[\fMIN := \{e : (\forall j < e)\: [\phe_j \neq \phe_e]\}.\]
Despite its outwardly congenial appearance, $\fMIN$ exhibits some barbaric
properties. In 1972, Meyer demonstrated that $\fMIN$ admits a neat Turing
characterization, namely $\fMIN \equiv_\T \zero''$ \citep{Mey72}.  Yet it is
difficult, if not impossible, to pin down the degree of $\fMIN$ for stronger
reductions (see Section~\ref{noneffectiveorderings}) \citep{Kin77},
\citep{Sch98}.  Along these lines, Schaefer showed that, unlike the familiar
index sets, a strong reduction will never reduce the halting set to $\fMIN$
(Section~\ref{omega immune}) \citep{Sch98}.

At this point, the reader is likely wondering, ``In the definition of $\fMIN$,
what happens if we replace the functions with sets, and also try to replace
equality with other equivalence relations?''  We trace the scant history of
this problem.  In Spring 1990 (according to the best recollection of the
author), John Case issued a homework assignment with the following definition
\citep{Cas06}:
\[\fMIN^* := \{e : (\forall j < e)\: [\phe_j \neq^* \phe_e]\}, \]
where $=^*$ means equal except for a finite set.  Case notes that $\fMIN^*$ is
$\Sigma_2$-immune, although his assignment exclusively refers to the $\Sigma_2$
sets as ``$\lim$-r.e$.$'' sets:
\begin{defn}
$A$ is \emph{$\lim$-r.e$.$} iff there exists a computable function $g$
satisfying
\[(\forall x)\: [\chi_A(x) = \liminf_{t \to \infty}g(x,t)].\]
\end{defn}
On October 1, 1996, six years after the initial homework assignment, Case
introduced the set $\fMIN^*$ to Marcus Schaefer in an email.

The following year, Schaefer published a Master's Thesis on minimal indices
\citep{Sch98}, which became the first public account of $\fMIN^*$.  In his
survey thesis, Schaefer proved that $\fMIN^* \join \zero' \equiv_\T \zero'''$,
leaving open the tantalizing question of whether or not $\fMIN \equiv_\T
\zero'''$. All that would be required to answer this question affirmatively is
to show that $\fMIN^* \geq_\T \zero'$, care of Schaefer's result.  The reader
is encouraged to attempt this reduction before proceeding. This concludes our
historical remarks.  The remaining scholarship on this subject is probably
contained in this thesis.

In attempt to comprehend the world of minimal indices, we characterize the
spectral sets from Section~\ref{MINdefinitions} in three ways.  First, we
describe the sets in terms of the arithmetic hierarchy.  The arithmetic
hierarchy gives us an idea of the complexity of sets by describing how many
quantifiers are needed to determine membership.  The arithmetic hierarchy does
not, however, tell us exactly which sets are computable from a set in question.
For this reason, we devote Chapter~\ref{Turing characterizations} to a
discussion of Turing degrees for minimal indices.  We discuss several reduction
techniques, and ultimately discover that spectral sets contain as much
information as any set with the same complexity (modulo some nontrivial
iteration of the halting set). Even without an oracle, this can still be true.
In particular, we show that it possible to effectively enumerate the partial
computable functions in such a way that the spectral sets are Turing equivalent
to the familiar sets $\zero'$, $\zero''$, $\zero'''$, $\dotsc$ (see
Chapter~\ref{a kolmogorov numbering}).  Both Chapter~\ref{Turing
characterizations} and Chapter~\ref{a kolmogorov numbering} make use of the
Forcing Lowness Lemma (Lemma~\ref{forcerecursiveorindependent}), an interesting
result in its own right.

In Chapter~\ref{immunitychapter}, we classify spectral sets in terms of
immunity.  Immunity is a measure of the ``thinness'' of a set, and it is an
especially appropriate benchmark for spectral sets.  All of the sets we
consider are $\omega$-immune and not hyperimmune (see Section~\ref{omega
immune} and Section~\ref{properties of MINT}), but another immunity notion is
useful for comparisons.  In particular, we examine immunity with respect to the
arithmetic hierarchy.  Weak equivalence relations give rise to ``thin''
spectral sets which are immune against high levels of the arithmetic hierarchy.
This may be the first time that a class of sets has been characterized in this
manner.  Of note in Chapter~\ref{immunitychapter} is the $\Pi_3$-Separation
Theorem (Theorem~\ref{pi3separation}), which gives a neat, if not surprising,
way for distinguishing between spectral sets in $\Pi_3$.

Chapter~\ref{thickville} crushes false generalizations that one might surmise
after reading just the first three chapters.  At the same time, it provides
additional direction by introducing an operator on equivalence relations.  We
learn that spectral sets occupy every level of the arithmetic hierarchy,
including $\Sigma_n - \Pi_n$, and we gain intuition for why a simple converse
to the immunity-completeness theorems does not appear in
Section~\ref{completeness criterion}.  Furthermore, we observe that the
location of an equivalence relation within the arithmetic hierarchy tells us
little about its immunity.  The operator in this chapter offers a way to
compare the power of nondeterminism against the jump operator within the realm
of the c.e$.$ sets.  If we accept the notion that weaker relations indicate
greater computational power, then the jump operator comes out on top.

Finally, what is the sparsest set you can imagine?  We follow this line of
thought into the last chapter.  In particular, there exists a spectral set so
sparse that it doesn't contain any infinite arithmetic sets.  Most of the work
in Chapter~\ref{The Peak Hierarchy Theorem} is devoted to showing that this
remarkable spectral set is not hyperimmune.  Consequently, we are able to show
that there exists a $\0$-majorized set which takes a bite out of every
arithmetic set but never eats the whole thing.

\section{Preliminaries}
The computability background required for this paper is completely covered in
Soare's book \citep{Soa87}, and we use the standard notation from his book
throughout this thesis.  The other important reference for this thesis is
Schaefer's survey on minimal indices \citep{Sch98}.  Schaefer's work is
approachable and comprehensive.  We will not cover all his results here.

\subsection{Basic computability theory}

p.c$.$ stands for partial computable, and c.e$.$ stands for computably
enumerable \citep{Soa87}.  ``c.e$.$ in $A$'', or equivalently, ``$A$-c.e$.$''
means enumerable with an $A$ oracle.  $A$-computable means computable with an
$A$ oracle.  We say $A$ is \emph{co-c.e$.$}if $A$ is c.e.  Unless otherwise
specified, we assume a fixed enumeration of the partial computable functions,
$\phe_0$, $\phe_1$, $\dotsc$. $W_0$, $W_1$, $\dotsc$ is an enumeration of their
domains. $\dom f$ denotes the domain for a partial function $f$, and $\rng f$
is the range of $f$. $\converge$ is for converge, and $\diverge$ is for
diverge. $\phe_{e,t}(x) \converge$ means that $\phe_e(x)$ converges in $t$
steps.  $(\mu x)$ means ``the smallest $x$ such that $\dotsc$.''  $K := \{x :
\phe_x(x)\converge\}$ denotes the halting set, $'$ denotes the jump operator,
$\null^{(n)}$ is the $n^\text{th}$ iteration of the jump operator, and
$\zero^{(\omega)} := \{\pair{x,n} : x \in \zero^{(n)}\}$. ``$\lim$'' means
limit.

For any set $A$, $\compliment{A}$ denotes the complement of $A$.  $\size{A}$
denotes the cardinality of $A$.  $\chi_A$ is the characteristic function for
$A$, which we sometimes write as just $A$.  $A(n)$ is the $n^\text{th}$ element
of $A$ under the usual ordering.    $\omega$ denotes the set of natural
numbers.  $\pair{\cdot, \cdot} : \omega \times \omega \to \omega$ is a
bijective pairing function.

For any equivalence class $\equiv_\alpha$,the \emph{$\equiv_\alpha$-degree} of
a set $A$ is the class of all sets equivalent to $A$ under $\equiv_\alpha$.  If
no equivalence relation is specified, we mean Turing equivalence. The degree
containing $\zero$ is $\0$, the degree containing $\zero'$ is $\0'$, and the
degree containing $\zero^{(n)}$ is $\0^{(n)}$.  An acquaintance with the
statements of the $s$-$m$-$n$ Theorem, the Recursion Theorem, and the Jump
Theorem \citep[Theorem III.2.3]{Soa87} are reccomended.

\begin{defn}
$\Psi_e^{A} (x)$ denotes the output of oracle Turing machine $e$ with oracle
$A$ on input $x$.  $\psi_e^{A} (x)$ is the corresponding \emph{use function},
the maximum query made to the oracle during computation (if it converges).
$W_e^A$ denotes the domain of $\Psi_e^{A} (x)$.

Occasionally, we will also use $\psi$ to denote a partial computable function.
In this case $\psi$ will not receive an oracle superscript, so as not to be
confused with the use function.
\end{defn}

\begin{defn}
For any set $A$,
\[A \restr n := A \intersect \{0, \ldots, n \}\]
is the $n^\text{th}$ \emph{initial segment} of $A$.
\end{defn}

The join operator allows us to use two or more sets as oracles simultaneously.
\begin{defn}
Let $A$ and $B$ be sets.  We define the \emph{join} of $A$ and $B$, denoted $A
\join B$, to be the set \[A \join B := \{2x: x \in A \} \cup \{2x + 1 : x \in
B\}.\]  For a sequence of sets $\{A_i\}$, define the \emph{infinite join} to be
\[\ijoin_{i \in \omega} A_i := \{\pair{x,i} : x \in A_i \}.\]
\end{defn}

Sometimes we only care about the first number in an ordered pair:

\begin{nota}[projections]
Let $\pi_1: \omega \to \omega$ denote the function which maps pairs to their
first coordinates, i.e.
\[ \pi_1 \left(\pair{x,y} \right) := x.\]
Similarly,
\[ \pi_2 \left(\pair{x,y} \right) := y.\]
\end{nota}

\begin{defn}
An integer $n$ is an \emph{$i^\text{th}$ prime power} if $n = p_i^k$ for some
$k \geq 1$, where $p_i$ is the $i^\text{th}$ prime number.  If $n$ is an
$i^\text{th}$ prime power for some $i$, then we may simply say $n$ is a
\emph{prime power}.
\end{defn}

\begin{defn}
Let $A$ and $B$ be thing.  Wacka wacka.
\end{defn}

\begin{defn}
$A$ is called an \emph{index set} if \[[x \in A \andd \phe_y = \phe_x]\quad
\implies \quad y \in A.\]
\end{defn}

\begin{defn}
A few familiar index sets will come into play.  For $n \geq  0$:
\begin{align*}
\INF &:= \{e : \size{W_e} = \infty\}, \\
\TOT &:= \{e : W_e = \omega\}, \\
\COF &:= \{e : W_e =^* \omega\}, \\
\mCOMP &:= \{e : W_e \equiv_\m K \}, \\
\LOW^n &:= \{e : (W_e)^{(n)} \equiv_\T \zero^{(n)} \}, \\
\HIGH^n &:= \{e : (W_e)^{(n)} \equiv_\T \zero^{(n+1)} \}.
\end{align*}

Note that $$\LOW^0 = \{e : W_e \equiv_\T \zero \},$$ and $$\HIGH^0 = \{e : W_e
\equiv_\T K \}.$$
\end{defn}

We sometimes view a set of natural numbers as a matrix, in which case the rows
may have special meanings:

\begin{defn} For any set $A$, \label{setrowdefn}
\begin{align*}
A^{[y]} &:= \{\pair{x,y} : \pair{x,y} \in A\}, \\
A^{[\widetilde{x}]} &:= \{\pair{x,y} : \pair{x,y} \in A\}.
\end{align*}
Here $A^{[y]}$ is the ``$y^\text{th}$ row of $A$,'' and $A^{[\widetilde{x}]}$
is the ``$x^\text{th}$ column of A,''
\begin{align*}
A^{[\leq y]} &:= \bigcup_{z \leq y} A^{[z]}, \quad \text{and} \\[1ex]
 A^{[> y]} &:= \bigcup_{z > y} A^{[z]}.
\end{align*}
\end{defn}

\begin{defn}
\begin{defenum}
\item A subset $A \subseteq B$ is a \emph{thick} subset if $A^{[y]} =^* B^{[y]}$ for all
$y$.

\item $B$ is \emph{piecewise computable} if $B^{[y]}$ is computable for all $y$.
\end{defenum}
\end{defn}

\begin{defn} Let $\paren{D_e}_{e\in\omega}$ be the canonical numbering of the finite sets.
\begin{defenum}
\item A set is \emph{immune} if it is infinite and contains no infinite c.e$.$
sets.

\item A set $A$ is \emph{hyperimmune} if it is infinite and there is no computable
function $f$ such that:
\begin{enumerate}
\item $\paren{D_{f(i)}}_{i \in \omega}$ is a family of pairwise disjoint sets,
and

\item $D_{f(i)} \intersect A \neq \emptyset$.
\end{enumerate}
\end{defenum}
\end{defn}

\subsection{Reductions and arithmetic hierarchy}
We list the main reductions and equivalence relations.

\begin{defn}[reductions] Let $A$ and $B$ be sets.
\begin{defenum} \label{familiar reductions}
\item  Let $f$ and $g$ be functions.
\[f =^* g \bigiff (\exists N) (\forall x > N)\: [f(x) = g(x)].\]
If $f$ and $g$ are the characteristic functions for $A$ and $B$ respectively,
then \[A =^* B \bigiff f =^* g.\]  Furthermore, $A \subseteq^* B$ when $A =^*
C$ for some $C \subseteq B$.

\item $A \leq_\m B$ if there exists a computable function $f$ such that
\[x \in A \bigiff f(x) \in B.\]

\item $A \leq_1 B$ if $A \leq_m B$ via an \emph{injective} function $f$.

\item Let $\{\sigma_n\}$ be an enumeration of all propositional truth tables with
predicates of the form ``$k \in S$,'' where $k \in \omega$.  We say that a set
$X$ \emph{satisfies} a truth table $\sigma_n$, or $X \satisfies \sigma_n$ if
the proposition $\sigma_n$ is true when ``$S$'' is interpreted as $X$.

$A \leq_\Tt B$ just in case there exists a computable function $f$ such that
\[ x \in A \bigiff B \satisfies \sigma_{f(x)}.\]

\item $A \leq_\btt B$ means that $A \leq_\Tt B$ via an $f$ which requires only
constant many queries to $B$.

\item $A \leq_\T B$ if there exists an index $e$ such that $\chi_A = \Psi_e^B$.
 In general, we write $f \leq_\T B$ for a function $f$ if $f = \Psi_e^B$ for
 some $e$.

\item $A \leq_\bT B$ if $A \leq_\T B$ and the largest query to the $B$ oracle
is computably bounded.  That is, $A \leq_\bT B$ if there exists a computable
function $f$ and an index $e$ such that for all $x$, $$\chi_A (x) = \Psi_e^{B
\restrs f(x)} (x).$$  Alternatively, $A \leq_\bT B$ if there exists a
computable function $f$ such that
\[ x \in A \bigiff B \satisfies \xi_{f(x)},\]
where $\{\xi_n\}$ is an enumeration of p.c$.$ truth tables which converge upon
satisfaction and diverge otherwise.  For this reason, $bT$-reductions are also
called ``weak truth-table'' reductions.

A function $f \leq_\bT B$ if $f \leq_\T B$ and we can computably bound the
largest query to $B$.

\item For every $n < \omega$, let
\begin{align*}
A \leq_{\T^{(n)}} B \enskip &\iff \enskip A^{(n)} \leq_\T B^{(n)},
\intertext{and} A \leq_{\T^{(\omega)}} B \quad &\iff \quad (\exists n)\:
[A^{(n)} \leq_\T B^{(n)}].
\end{align*}

\end{defenum}
If $A \leq_\alpha B$ and $B \leq_\alpha A$ for some partial ordering
$\leq_\alpha$, then we write $A \equiv_\alpha B$.  If $A \leq_\alpha B$ and $B
\not\leq_\alpha A$, then $A <_\alpha B$.  If $A \not\leq_\alpha B$ and $B
\not\leq_\alpha A$, then $A \incomp{\alpha} B$.  This notation applies to all
of the reductions in this definition. Finally, two sets are equal if they are
equal.
\end{defn}

We define the member classes of the arithmetic hierarchy: $\Delta_n$,
$\Sigma_n$, and $\Pi_n$ for $n \geq 0$.
\begin{defn}[arithmetic hierarchy] \label{arithmetichierarchydefn}
Let $n \geq 1$.
\begin{defenum}
\item $\Delta_0 = \Sigma_0 = \Pi_0$ is the class of computable sets.

\item $A \in \Sigma_n$ if there exists a computable relation $R$ such that
\[x \in A \bigiff (\exists y_1)(\forall y_2) (\exists y_3) \dotso (Q y_n)\:
R(x, y_1, \dotsc, y_n),\] where $Q$ is $\forall$ if $n$ is even, and $\exists$
if $n$ is odd.  Similarly, \label{arithmetichierarchydefnsigman}

\item $A \in \Pi_n$ if there exists a computable relation $R$ such that
\[x \in A \bigiff (\forall y_1)(\exists y_2) (\forall y_3) \dotso (Q y_n)\:
R(x, y_1, \dotsc, y_n),\] where $Q$ is $\exists$ if $n$ is even, and $\forall$
if $n$ is odd. \label{arithmetichierarchydefnpin}

\item $\Delta_n := \Sigma_n \intersect \Pi_n$.
\end{defenum}
\end{defn}

We also relativize the arithmetic hierarchy in the following way:

\begin{defn}[relativized arithmetic hierarchy]
Let $S$ be a set, and let $n \geq 1$.
\begin{defenum}
\item $\Delta_0^S = \Sigma_0^S = \Pi_0^S$ is the class of $S$-computable sets.

\item $A \in \Sigma_n^S$ is just as in
Definition~\ref{arithmetichierarchydefn}\nlb(\ref{arithmetichierarchydefnsigman}),
except that ``a computable relation $R$'' is replaced with ``an $S$-computable
relation $R$.''

\item $A \in \Pi_n^S$ is just as in
Definition~\ref{arithmetichierarchydefn}\nlb(\ref{arithmetichierarchydefnpin}),
except that ``a computable relation $R$'' is replaced with ``an $S$-computable
relation $R$.''

\item $\Delta_n^S := \Sigma_n^S \intersect \Pi_n^S$.
\end{defenum}

\end{defn}

\begin{defn}
A set $A$ is:
\begin{defenum}
\item $\Sigma_n$-complete if $A \in \Sigma_n$ and for every $B \in \Sigma_n$, $B \leq_\m
A$.

\item $\Pi_n$-complete if $A \in \Pi_n$ and for every $B \in \Pi_n$, $B
\leq_\m A$.
\end{defenum}
\end{defn}

The reader who is reading about the arithmetic hierarchy for the first time
should familiarize herself with (Relativized) Post's Theorem, the Hierarchy
Theorem, and the Limit Lemma \citep{Soa87}.

\subsection{Minimal indices} \label{MINdefinitions}
We formally define our objects of study.

\begin{defn}
Let $\equiv_\alpha$ be an equivalence relation on sets.  Then
\[ \MIN^{\equiv_\alpha} := \{ e : (\forall j < e)\: [W_j \not \equiv_\alpha W_e]\}.\]
Similarly, for an equivalence relation $\equiv_\beta$ on functions we define,
\[ \fMIN^{\equiv_\beta} := \{ e : (\forall j < e)\: [\phe_j \not \equiv_\beta \phe_e]\}.\]
A set of either form is called a \emph{spectral set}, or, equivalently, a
\emph{$\MIN$-set}.
\end{defn}

 We will refer to certain spectral sets often, and we use the following abbreviations for these
 sets.  We employ equivalence relations from Definition~\ref{familiar
 reductions}.
\begin{defn}\label{defnminsets}
For notational clarity, we sometimes abbreviate the relations $=^*$,
$\equiv_\m$, $\equiv_\T$, and $\equiv_{\T^{(n)}}$ as $*$, $\m$, $\T$ and
$\T^{(n)}$, repectively.  The following are in effect, for $n \geq 0$:
\begin{align*}
\MIN &:= \{e : (\forall j < e)\: [W_j \neq W_e] \}, \\
\MIN^* &:= \{e : (\forall j < e)\: [W_j \neq^* W_e] \}, \\
\MIN^\m &:= \{e : (\forall j < e)\: [W_j \not\equiv_\m W_e] \}, \\
\MIN^\T &:= \{e : (\forall j < e)\: [W_j \not\equiv_\T W_e] \}, \\
\MIN^{\T^{(n)}} &:= \{e : (\forall j < e)\: [W_j \not\equiv_{\T^{(n)}} W_e] \},
\intertext{and}
\MIN^{\T^{(\omega)}} &:= \bigcap_{n \in \omega} \MIN^{\T^{(n)}} \\
&=\{e : (\forall j < e)(\forall n)\: [(W_j)^{(n)} \not\equiv_\T (W_e)^{(n)}]\}.
\end{align*}
\end{defn}

In the case of $\MIN^\m$, we modify the usual definition of $\equiv_\m$ so that
all recursive sets, including $\zero$ and $\omega$, have the same $\m$-degree.
This makes Theorem~\ref{JLSSmain} true without modification.

A similar set of notations applies for indices of minimal functions, but only
for a few specific equivalence relations.

\begin{nota}  We shall consider the following ``function'' spectral sets.
\begin{align*}
\fR &:= \{e : (\forall j < e)\: [\phe_j (0) \neq \phe_e (0) ] \}, \\
\fMIN &:= \{e : (\forall j < e)\: [\varphi_j \neq \varphi_e] \}, \\
\fMIN^* &:= \{e : (\forall j < e)\: [\varphi_j \neq^* \varphi_e] \}.
\end{align*}
\end{nota}

Occasionally, we will want to compute the minimal index for a c.e$.$ set:
\begin{defn}
For every equivalence relation $\equiv_\alpha$, we define a function
$\mn^{\equiv_\alpha}$ by
\[\mn^{\equiv_\alpha}(e) := (\mu x)\: [W_x \equiv_\alpha W_e].\]
If $\equiv_\alpha$ is not specified in the notation, we mean equality.
\end{defn}

The following proposition is easily verified.
\begin{prop} \label{MINinclusionprop}
Let $\equiv_\alpha$ and $\equiv_\beta$ be equivalence relations.  Assume that for all $X,
Y \subseteq \omega$, \[X \equiv_\alpha Y \implies X \equiv_\beta Y.\]  Then
$\MIN^{\equiv_\alpha} \supseteq \MIN^{\equiv_\beta}$.
\end{prop}

\begin{cor} \label{inclusionpropexamples}
\begin{thmenum}
\item $\fMIN \supsetneq \MIN \supsetneq \MIN^*$,

\item $\fMIN \supsetneq \fMIN^* \supsetneq \MIN^*$,

\item $\MIN \supsetneq \MIN^* \supsetneq \MIN^\m \supsetneq \MIN^\T \supsetneq
\MIN^{\T'} \dotsb$.
\end{thmenum}
\end{cor}

In the following proposition, $\equiv_\alpha$ can be taken to be any familiar
intermediate reduction, such as $\equiv_\btt$, $\equiv_\Tt$, or $\equiv_\bT$.
It might appear, in light of Proposition~\ref{minm'isminT}, that the spectral
sets form a simple, linear ordering under reverse inclusion. However, in
Chapter~\ref{thickville} and Section~\ref{almost thickness} we explore a
natural class of equivalence relations which do not fit between
$\MIN^{\T^{(n)}}$ and $\MIN^{\T^{(n+1)}}$ for any $n$.

\begin{prop}\label{minm'isminT}
\begin{thmenum}
\item For every $n \geq 0$, $\MIN^{\m^{(n+1)}} = \MIN^{\T^{(n)}}.$ \label{minm'isminTmain}

\item Let $\equiv_\alpha$ be any equivalence relation which is weaker than $\equiv_1$
and stronger than $\equiv_\T$.  For any sets $A$ and $B$, let
\[A \equiv_{\alpha^{(n)}} B \quad \iff \quad A^{(n)} \equiv_\alpha B^{(n)},\]
and define
\[\MIN^{\alpha^{(n)}} := \{e : (\forall j < e)\ [W_j \not\equiv_{\alpha^{(n)}} W_e]\}.\]
Then for all $n$,
\[\MIN^{\T^{(n)}} \supseteq \MIN^{\alpha^{(n+1)}} \supseteq \MIN^{\T^{(n+1)}}.\]
\label{minm'isminTgen}
\end{thmenum}
\end{prop}

\begin{proof}
\begin{proof}[(\ref{minm'isminTmain})]
It suffices to show that for any sets $A$ and $B$,
\begin{equation}
A^{(n+1)} \equiv_\m B^{(n+1)} \iff A^{(n)} \equiv_\T B^{(n)}.
\label{jumpthmeqn}
\end{equation}
We show (\ref{jumpthmeqn}) by proving the Jump Theorem~\cite[Theorem
III.2.3]{Soa87}, namely:
\[A' \leq_\m B' \iff A \leq_\T B.\]

Assume $A' \leq_\m B'$ via a computable function $h$.  Then
\begin{align*}
A &\leq_\m A' \leq_\m B', \\
\compliment{A} &\leq_\m A' \leq_\m B'.
\end{align*}
Indeed, $A \leq_\m A'$ via the function $f$ defined by
\[ \Psi_{f(x)}^A(n) =
\begin{cases}
1 &\text{if $x \in A$} \\
\diverge &\text{otherwise},
\end{cases}\]
because
\[x \in A \iff \Psi_{f(x)}^A[f(x)] \converge \iff f(x) \in A'.\]
An analogous function $g$ yields $\compliment{A} \leq_\m A'$.

It follows that $A$ and $\compliment{A}$ are c.e$.$ in $B$ via the enumerations
$x \in A_s \iff h \of f(x) \in B'_s$ and $x \in \compliment{A}_s \iff h \of
g(x) \in B'_s$, where $B'_s$ is a $B$-enumeration of $B'$. Therefore $A \leq_\T
B$.

Conversely, assume $A \leq_\T B$.  Since $A'$ is c.e$.$ in $A$, $A'$ must be
c.e$.$ in $B$.  This means that $A' \leq_\m B'$, since $B'$ is $m$-complete
relative to $B$.
\end{proof}

\begin{proof}[(\ref{minm'isminTgen})]
According to the Jump Theorem~\citep[Theorem III.2.3]{Soa87}, \[A \equiv_\T B
\iff A' \equiv_1 B'.\] Therefore,
\[A^{(n)} \equiv_\T B^{(n)} \implies A^{(n+1)} \equiv_1 B^{(n+1)} \implies A^{(n+1)} \equiv_\alpha B^{(n+1)},\]
and more obviously,
\[A^{(n+1)} \equiv_\alpha B^{(n+1)} \implies A^{(n+1)} \equiv_\T B^{(n+1)}.\qedhere\]
\end{proof}
\end{proof}

\section[Complexity]{Complexity of spectral sets} \label{arithmetic properties}

We place spectral sets in the arithmetic hierarchy.  Our lower bounds
immediately show that $\MIN$-sets are not computable, although our laconic
proofs do not involve the familiar technique of reduction to the halting set.
Spectral sets can be found in every level of the arithmetic hierarchy.  Unlike
index sets, which are always $\geq_\m \K$ (Rice's Theorem \citep{Soa87}),
spectral sets never have this property.  In fact, $K$ doesn't even $btt$-reduce
to $\MIN$-sets (Corollary~\ref{notbtttozero'}).

Based on Corollary~\ref{inclusionpropexamples}, it would be reasonable to
extrapolate that $A \supseteq B$ implies that $B$ lies in a higher arithmetic
level than $A$.  This turns out not to be the case when we considered minimal
indices of functions.  Indeed, there is a notable exception:
\begin{defn}[\citet{Sch98}]
We call
\[\fMIN = \{e : (\forall j < e)\, [\varphi_j \neq \varphi_e]\} \]
the \emph{set of minimal indices for functions}, and
\[\fR = \{e : (\forall j < e)\, [ \phe_j (0) \neq \phe_e (0) ] \} \]
denotes the \emph{set of shortest descriptions for nonegative integers}.
\end{defn}
$\fMIN \supseteq \fR$, yet $\fMIN \in \Sigma_2 - \Pi_2$ and $\fR \in \Delta_2$.
$\fR$, which does not appear to have a spectral analogue for sets, highlights a
potential difference between minimal indices for sets and minimal indices for
functions.  We shall exhibit an infinite $\Delta_2$ subset of $\MIN$ in
Section~\ref{Deltaimmune}, however, it is not a spectral set.  We remark that
the results in Sections~\ref{arithupperbounds} and
Sections~\ref{arithlowerbounds} are by-and-large subsumed by
Corollary~\ref{bestTuring}.

\subsection{Upper bounds}\label{arithupperbounds} We reveal upper bounds for
a number of sets, including the following rare example.
\begin{defn} \label{MINtmdefn}
\[\MIN^\tm := \{ e : (\forall j < e)\, (\exists \pair{x,t}) \left[W_{j,t}(x) \neq
W_{e,t}(x)\right]\}.\]
\end{defn}
The $\equiv_\tm$ identifies indices which are not only equal, but their
respective computations converge in exactly the same amount of time.
$\equiv_\tm$ is the only $\MIN$-set in this thesis which contains $\MIN$; the
rest are subsets of either $\MIN$ or $\fMIN$.

In Theorem~\ref{lowerbound}, we will show that the following upper bounds are
optimal (except for part~(\ref{mintminsigma1}), which follows from
Theorem~\ref{MINtmis0'}).  In light of Proposition~\ref{minm'isminT},
Proposition~\ref{basicarithmetic} also shows that $\MIN^{\m^{(n)}} \in
\Pi_{n+3}$.

\begin{prop} \label{basicarithmetic}
Let $n \geq 0$.  Then
\begin{thmenum}
\item $\MIN^\tm \in \Sigma_1$. \label{mintminsigma1}
\item $\fR \in \Delta_2$. \label{fRindelta2}
\item $\MIN, \fMIN \in \Sigma_2$. \label{mininsigma2}
\item $\MIN^*, \fMIN^* \in \Pi_3$. \label{min*inpi3}
\item $\MIN^{\m} \in \Pi_3$. \label{minminpi3}
\item $\MIN^{\equiv_1} \in \Pi_3$. \label{min1inpi3}
\item $\MIN^{\T^{(n)}} \in \Pi_{n+4}$. \label{mintinpi4}
\end{thmenum}
\end{prop}

\begin{proof}

\begin{proof}[(\ref{mintminsigma1})]
Immediate from the definition.
\end{proof}

\begin{proof}[(\ref{fRindelta2})]
$\phe_j(0) = \phe_e(0)$ can be decided with a $\zero'$ oracle.  So $\fR \in
\Delta_2$ by the Limit Lemma \citep{Soa87}.
\end{proof}

\begin{proof}[(\ref{mininsigma2})]
$\{ \pair{j,e} : W_j = W_e \} \in \Pi_2$ \citep{Soa87}.
\end{proof}

\begin{proof}[(\ref{min*inpi3})]
$\{ \pair{j,e} : W_j =^* W_e \} \in \Sigma_3$ \citep{Soa87}.
\end{proof}

\begin{proof}[(\ref{minminpi3})]
For any c.e$.$ sets $A$ and $B$,
\[A \leq_\m B \iff (\exists e) (\forall x)\ \left[\phe_e (x) \converge \quad \& \quad (x \in A \iff \phe_e(x) \in B)
\right],\] which shows that $A \leq_\m B$ is a $\Sigma_2^{\zero'}$ relation. It
follows that $A \equiv_\m B$ is also a $\Sigma_2^{\zero'}$ relation.  In
particular, for
\[C := \left\{ \pair{j,e} : W_j \equiv_\m W_e \right\}, \]
we have
\[C \in \Sigma_2^{\zero'} = \Sigma_3.\]
Hence
\[ e \in \MIN^{\m} \iff (\forall j < e)\ \left[\pair{j,e} \not \in C \right],\]
which places $\MIN^{\m} \in \Pi_3$.
\end{proof}

\begin{proof}[(\ref{min1inpi3})]
The same idea from (\ref{minminpi3}) works because injectivity can be tested
with a $\zero'$ oracle.
\end{proof}

\begin{proof}[(\ref{mintinpi4})]
For any sets $A$ and $B$,
\begin{align*}
A \leq_\T B &\iff (\exists e) \left[ A = \Psi_e^B \right] \\
&\iff (\exists e) (\forall x) \left[ \Psi_e^B (x) \converge \quad \& \quad (x
\in A \iff \Psi_e^B (x) = 1) \right],
\end{align*}
which shows that $A \leq_\T B$ is a \[\Sigma_2^{B' \join \: (A \join B)} =
\Sigma_2^{A \join B'}\] relation, and it follows that $A \equiv_\T B$ is a
$\Sigma_2^{A' \join B'}$ relation.  In particular, for
\[C_n := \left\{ \pair{j,e} : (W_j)^{(n)} \equiv_\T (W_e)^{(n)} \right\}, \]
we have
\[C_n \in \Sigma_2^{\paren{(W_j)^{(n)}}' \join \paren{(W_e)^{(n)}}'} \subseteq \Sigma_2^{\zero^{(n+2)}} = \Sigma_{n+4}.\]
It follows that
\[ e \in \MIN^{\T^{(n)}} \iff (\forall j < e) \left[\pair{j,e} \not \in C_n \right],\]
which makes $\MIN^{\T^{(n)}} \in \Pi_{n+4}$.
\end{proof}
\end{proof}

\subsection{Lower bounds} \label{arithlowerbounds} It's not too hard to show
that $\MIN$-sets are noncomputable (modulo a few well-known theorems), however,
the more familiar method of $m$-reduction to the halting set doesn't work.
Theorem~\ref{lowerbound}(\ref{MINlowerbound}) was known to Meyer \citep{Mey72},
and I attribute Theorem~\ref{lowerbound}\nlb(\ref{fRlowerbound}) to Lance
Fortnow.

\begin{thm} \label{lowerbound}
Let $n \geq 0$.
\begin{thmenum}
\item $\fR \not \in \Sigma_1 \union \Pi_1$. \label{fRlowerbound}

\item $\MIN, \fMIN \not \in \Pi_2$. \label{MINlowerbound}

\item $\MIN^*, \fMIN^* \not \in \Sigma_3$. \label{MIN*lowerbound}

\item $\MIN^\m \not \in \Sigma_3$. \label{MINmlowerbound}

\item $\MIN^{\equiv_1} \not \in \Sigma_3$. \label{MIN1lowerbound}

\item $\MIN^{\T^{(n)}} \not \in \Sigma_{n+4}$. \label{MINTlowerbound}
\end{thmenum}
\end{thm}

\begin{proof}

\begin{proof}[(\ref{fRlowerbound})]
$\fR \not\in \Sigma_1$ follows immediately from the fact that $\fR$ is immune
\citep{Sch98}.  Suppose $\fR \in \Pi_1$.  Let $a$ be the smallest index such
that $\phe_a(0) \diverge$.  Define a computable function $f$ by way of the
$s$-$m$-$n$ Theorem \citep{Soa87} and the following constant function:
\[ \phe_{f(x)} (y) :=
\begin{cases}
(\mu t)\: [\phe_{x,t}(0) \converge] &\text{if $\phe_x (0) \converge$}, \\
\diverge &\text{otherwise.}
\end{cases}\]
Let \[K_0 := \{e : \phe_e(0) \converge \}.\]  $K_0$ is $\Sigma_1$-complete.
 Note that
\begin{align*}
e \in K_0 &\iff \phe_{f(e)}(0) \converge \\
&\iff  \bigl(\exists j \in [\{0, \dotsc, f(e) \} \intersect \fR] - \{a\}\bigr)
\left[\phe_j (0) \converge \andd \phe_{e, \phe_j(0) \converge}(0) \converge \right],\\
&\iff \bigl(\exists j \leq f(e)\bigr) \left[j \in \fR - \{a\} \andd \phe_j (0)
\converge \andd \phe_{e, \phe_j(0)}(0) \converge\right].
\end{align*}
This means that $\compliment{K_0} \in \Sigma_1$, since $j \in \fR - \{a\}
\implies \phe_j(0) \converge$.  But that's a contradiction, because now $K_0$
is computable.
\end{proof}

\begin{proof}[(\ref{MINlowerbound})]  Suppose that $\MIN \in \Pi_2$, let $a$ be the minimal index for $\omega$, and recall that
\[\TOT = \{e : W_e = \omega\}\]
is $\Pi_2$-complete \citep{Soa87}.  Then
\begin{align*}
\TOT &= (\MIN \intersect \TOT) \union (\compliment{\MIN} \intersect \TOT) \\
&= \{a\} \union \left\{e : (\forall j < e)\: \left[j \in \MIN - \{a\}\enskip
\implies \enskip W_j \neq W_e\right]\right\}.
\end{align*}
Now $\TOT \in \Sigma_2$, since $W_j = W_e$ can be decided in $\Pi_2$, and
because $\MIN - \{a\} \in \Pi_2$ by assumption.  This contradicts the fact that
$\TOT$ is $\Pi_2$-complete.
\end{proof}

\begin{proof}[(\ref{MIN*lowerbound})]
We reuse the argument from part (\ref{MINlowerbound}).  Suppose $\MIN^* \in
\Sigma_3$, let $a$ be the *-minimal index for $\omega$, and recall that the set
of cofinite indices \[\COF := \{e : W_e =^* \omega\}\] is $\Sigma_3$-complete
\citep{Soa87}. Then
\begin{align*}
\COF &= (\MIN^* \intersect \COF) \union (\compliment{\MIN^*} \intersect \COF) \\
 &= \{a\} \union \left\{e : (\forall j < e)\: \left[j \in \MIN^* - \{a\}\enskip
\implies \enskip W_j \neq^* W_e\right]\right\}
\end{align*}
Now $\COF \in \Pi_3$, since $W_j =^* W_e$ can be decided in $\Sigma_3$, and
because $\MIN^* - \{a\} \in \Sigma_3$ by assumption.  This contradicts the fact
that $\COF$ is $\Sigma_3$-complete.
\end{proof}

\begin{proof}[(\ref{MINmlowerbound})]
$\{e : W_e \equiv_\m C\}$ is $\Sigma_3$-complete whenever $C \neq \emptyset$,
$C \neq \omega$, and $C$ is c.e. \citep{Yat66b}.
\end{proof}

\begin{proof}[(\ref{MIN1lowerbound})]
$\{e : W_e \equiv_1 C\}$ is $\Sigma_3$-complete whenever $C$ is c.e$.$,
infinite, and coinfinite \citep{Her92}. Since $W_j \equiv_1 W_e$ is decidable
in $\Sigma_3$, the same argument again applies.
\end{proof}

\begin{proof}[(\ref{MINTlowerbound})]
Combining the Yates Index Set Theorem with the Sacks Jump Theorem yields
\begin{equation*}
\HIGH^n = \{e: W_e \equiv_{\T^{(n)}} \zero'\}
\end{equation*}
is $\Sigma_{n+4}$-complete, which is exactly what is needed to prove the
theorem. This fact seems to have been first observed by Schwarz in his PhD
thesis \citep[Theorem 3.3.1]{Sch82}, \citep[Theorem XII.4.4]{Soa87}. He writes
simply,
\begin{quote}
``We discovered the unexpectedly short argument [that $\HIGH_n$ is
$\Sigma_{n+4}$-complete] quite by accident, after having given up on finding
any more direct line of proof.''
\end{quote} \vspace{-5ex}\qedhere
\end{proof}

\end{proof}

\section{Noneffective orderings and other disasters}
\label{noneffectiveorderings}

$\MIN$ is sensitive to the order in which we list the partial computable
functions.  This is exacerbated by the fact that some c.e$.$ classes can be
enumerated without repetition  \citep{Weh95}, \citep{Ers77}.

\begin{defn}
A \emph{numbering} of a set $S$ is a surjective mapping of $\omega$ onto $S$.
If a numbering $\phe$ is computable, we say $\phe$ is a \emph{computable
numbering}.  If $S$ is not specified, we mean the set of partial computable
functions.  A p.c$.$ function $\phe$ is a \emph{p.c$.$ numbering} if  $$e
\mapsto \phe (\pair{e,\cdot})$$ maps onto the partial computable functions. For
any p.c$.$ numbering $\phe$, we denote the function $\phe (\pair{e,\cdot})$ by
$\phe_e$.
\end{defn}

\begin{defn} \label{godelkolmogorovnumbdefn}
A \emph{G\"{o}del numbering} $\phe$ is a p.c$.$ numbering such that if $\psi$
is a p.c$.$ function, then there exists a computable function $f$ satisfying
\[\phe_{f(e)} (x) = \psi (\pair{e,x}).\]
If in addition $f$ is linearly bounded, we say that $\phe$ is a
\emph{Kolmogorov numbering}.
\end{defn}

Since the $\psi$ in Definition~\ref{godelkolmogorovnumbdefn} might itself be a
numbering, we can effectively find a $\phe$-index for any algorithm when $\phe$
is a G\"{o}del numbering.  Furthermore, any reasonably encoded universal Turing
machine is a Kolmogorov numbering \citep{Sch74}.  We use a subscript to
indicate the numbering for a $\MIN$-set, as in $\MIN_\phe$.  If the subscript
is omitted, then we mean an arbitrary G\"{o}del numbering.

\subsection{G\"{o}del numberings}
The degrees for spectral sets are not always invariant with respect to
G\"{o}del numberings. For example, while we do not yet know the truth table
degree of $\fR$ \citep{Sch98}, we do have
Theorem~\ref{numberingsRANDfR}\nlb(\ref{numberingsRANDfRbT}).
Theorem~\ref{numberingsRANDfR}\nlb(\ref{fRpheis0'}) is due to Schaefer
\citep{Sch98}, and we shall revisit this result again in
Theorem~\ref{thickpheTuringcomplete}.  Most surprising, however, are Martin
Kummer's results on the truth-table degree of $\RAND$, the set of
\emph{Kolmogorov random strings}:

\begin{defn}\label{kolmogorovdefn}
For any finite string $x$, the Kolmogorov complexity of $x$ is
\[C(x) := \min\{\size{e} : \phe_e(0) = x \},\]
where $\size{\cdot}$ denotes the length of an integer encoded in binary.  The
set of \emph{Kolmogorov random strings} is
\[\RAND := \{x : C(x) \geq \size{x} \}.\]
\end{defn}

\begin{thm}[\citet{Sch98, Kum96}] \
\begin{thmenum} \label{numberingsRANDfR}
\item For any G\"{o}del numbering $\phe$, $\RAND_\phe \equiv_\bT \zero' \equiv_\bT
\fR_\phe$. \label{numberingsRANDfRbT}

\item There exists a Kolmogorov numbering $\phe$ such that $\fR_\phe \equiv_\Tt
\zero'$. \label{fRpheis0'}

\item For any Kolmogorov numbering $\phe$, $\RAND_\phe \equiv_\Tt \zero'$.
\label{RANDttequiv0'}

\item There exists a G\"{o}del numbering $\phe$ such that $\RAND_\phe \not \geq_\Tt \zero'$.
\end{thmenum}
\end{thm}

We turn to numberings for $\fMIN$.  Using Proposition~\ref{proper subsets of
immune sets}, Young gave a short proof of the following fact.

\begin{thm}[\citet{Mey72}]
For any G\"{o}del numbering $\phe$, there exists a G\"{o}del numbering $\psi$
such that $\fMIN_\psi <_1 \fMIN_\phe$.
\end{thm}

A more complicated argument reveals even more sensitivity.  Kinber was the
first to prove the following two results (both for G\"{o}del numberings),
however Schaefer's proof of (\ref{fMINpheis0''}) is decidely simpler.

\begin{thm}[\citet{Kin77}, \citet{Sch98}] \label{MINbttMIN}
\begin{thmenum}
\item There exist G\"{o}del numberings $\phe$ and $\psi$ such that $\fMIN_\phe \not\equiv_\btt
\fMIN_\psi$.

\item There exists a Kolmogorov numbering $\phe$ such that $\fMIN_\phe \equiv_\Tt
\zero''$. \label{fMINpheis0''}
\end{thmenum}
\end{thm}

The ``closer'' for $\fMIN$ numberings, however, is still at-large. Namely,
Meyer's question from 1972 of whether $\fMIN_\phe \equiv_\Tt \zero''$ \emph{for
all} G\"{o}del numberings $\phe$ remains open \citep{Mey72}.

\subsection{Enumeration without repetition}

\begin{defn}[\cite{Yat66b}]
For any equivalence relation $\equiv_\alpha$, we define \[G_{\equiv_\alpha}(C)
:= \{e : W_e \equiv_\alpha C \}.\]
\end{defn}

Yates proved the following theorem only for c.e$.$ sequences in the case
$\equiv_\alpha$ equal to $\equiv_\T$ \citep{Yat69}, but as we demonstrate here,
his argument easily generalizes to other relations.  In the following theorem,
when we say the $\equiv_\alpha$ degrees can not be enumerated, we mean that it
is impossible to make a list consisting of exactly one index from each
$\equiv_\alpha$ equivalence class.

\begin{thm}\label{enumerateT}
Let $\equiv_\alpha$ be an equivalence relation satisfying
\begin{equation*}
\{ \pair{i,j} : W_i \equiv_\alpha W_j \} \in \Sigma_n.
\end{equation*}
Assume there is some c.e$.$ set $C$ such that $G_{\equiv_\alpha} (C)$ is
$\Sigma_n$-complete.  Then no $\Sigma_n$ sequence of c.e$.$ sets enumerates the
$\equiv_\alpha$-degrees without repetition.
\end{thm}

\begin{proof}
Let $\equiv_\alpha$ be a relation satisfying
\begin{equation*}
\{ \pair{i,j} : W_i \equiv_\alpha W_j \} \in \Sigma_n.
\end{equation*}
Suppose there is some $A \in \Sigma_n$ which contains exactly one index from
each $\equiv_\alpha$ class. Let $c \in A$ be the index such that $W_c
\equiv_\alpha C$. Then
\[G_{\equiv_\alpha}(C) = \left\{ e : (\forall k) \left[k \in A - \{c\} \implies W_e \not\equiv_\alpha W_k \right] \right\} \in \Pi_n,\]
which implies that $G_{\equiv_\alpha}(C)$ is not $\Sigma_n$-complete.
\end{proof}

\begin{cor}
Let $n \geq 0$.  Then
\begin{thmenum}
\item the $=^*$ degrees can not be enumerated without repetition, \label{*enum}

\item the $\equiv_\m$ degrees can not be enumerated without repetition, and
\label{menum}

\item the $\equiv_{\T^{(n)}}$ degrees can not be enumerated without repetition.
\label{Tenum}
\end{thmenum}
\end{cor}

\begin{proof}

\begin{proof}[(\ref{*enum})]
$=^*$ is a $\Sigma_3$ relation on c.e$.$ sets, and $G_{*}(\omega) = \COF$ is
$\Sigma_3$-complete \citep{Soa87}.
\end{proof}

\begin{proof}[(\ref{menum})]
$\equiv_\m$ is a $\Sigma_3$ relation on c.e$.$ sets, and $G_\m(\K)$ is
$\Sigma_3$-complete (see Theorem~\ref{lowerbound}\nlb(\ref{MINmlowerbound})).
\end{proof}

\begin{proof}[(\ref{Tenum})]
$\equiv_\T$ is a $\Sigma_{n+4}$ relation on c.e$.$ sets, and $G_{\T^{(n)}}(\K)$
is $\Sigma_{n+4}$-complete (see
Theorem~\ref{lowerbound}\nlb(\ref{MINTlowerbound})).
\end{proof}
\end{proof}

Theorem~\ref{enumerateT} also eliminates the possibility that the $=^*$ sets
might be enumerable using a $\zero''$ oracle, that the $\equiv_\m$ might be
enumerable using a $\zero''$ oracle, and that the $\equiv_{\T^{(n)}}$ sets
might be enumerable using a $\zero^{(n+3)}$ oracle.

Theorem~\ref{enumerateT} does not apply when $\equiv_\alpha$ is the equality
relation, since $G_= (\omega)  = \TOT \in \Pi_2$. Consequently, Friedberg was
able to prove the following theorem \citep{Fri58} \citep{Gon02}, but we cite
Kummer for the elegance of his later proof.  Kummer's construction of a
Friedberg ordering is an application of the set $\MIN$.

\begin{thm}[\citet{Kum90}]
The c.e$.$ sets can be enumerated without repetition.
\end{thm}

The noneffective Friedberg ordering $\psi$ makes $\MIN_\psi = \omega$.  If we
are willing to entertain arbitrary numberings, then $\MIN$ can be any set we
like (or don't like).  The partial computable functions admit analogous
pathological numberings for $\fMIN$, thus threatening to turn our study of
minimal indices into a triviality.  For this reason, we hereby restrict our
attention to G\"{o}del numberings.  The remaining results in this thesis do not
depend on the particular choice of G\"{o}del numbering.  So from this point
forth, we simply fix an enumeration of the partial computable functions (with
one exception, Chapter~\ref{a kolmogorov numbering}).

%% file: turing.tex
\chapter[Turing characterizations]{Turing characterizations} \label{Turing characterizations}

When squeezed gently, a fair amount of information can be extracted from
spectral sets.  To show that $\zero^{(n)}$ reduces to a spectral set, one first
tries to achieve this (difficult) reduction with the aid of some oracle. By
repeatedly substituting with successively weaker oracles, eventually one
eliminates the oracle entirely (hopefully).  Each time that a weaker oracle is
introduced, a new reduction technique is required. This chapter is organized
according to technique.  Each section describes one or more reduction methods
which pertain to oracles of particular strength.

\section{Generic reductions}

Lemma~\ref{downonelevel} shows how to ``drop'' a $\MIN$-set ``down one level.''
We demonstrate an especially short proof which is peculiar to $\MIN^\m$,
however there is a canonical strategy which works for $\MIN$-sets in general.
The canonical strategy is presented in the proofs of (\ref{MINjoin0'is0''}) and
(\ref{MINTjoin0'''is0''''}).  In each case, we give the reduction in only one
direction because the opposite directions are immediate from our arithmetic
upper bounds (Proposition~\ref{basicarithmetic}). (\ref{MINjoin0'is0''}) and
(\ref{MIN*join0''is0'''}) first appeared in \citep{Sch98} and \cite{Mey72} for
$\fMIN$ and $\fMIN^*$, respectively.

\begin{lemma} For $n \geq 0$, \label{downonelevel}
\begin{thmenum}
\item $\MIN \join \zero' \equiv_\T \zero''$, \label{MINjoin0'is0''}

\item $\MIN^* \join \zero'' \equiv_\T \zero'''$, \label{MIN*join0''is0'''}

\item $\MIN^\m \join \zero'' \equiv_\T \zero'''$, \label{MINmjoin0''is0'''}

\item $\MIN^{\T^{(n)}} \join \zero^{(n+3)} \equiv_\T \zero^{(n+4)}$. \label{MINTjoin0'''is0''''}
\end{thmenum}
\end{lemma}

\begin{proof}
\begin{proof}[(\ref{MINjoin0'is0''})]
Let $a$ be the minimal index for $\TOT$, and let $e$ be any index.  Note that
$W_e = W_x$ for exactly one $x$ in
\[B:= \{0, \dotsc, e\} \intersect \MIN.\]
Since
\[ \{ \pair{j,e} : W_j \neq W_e \} \in \Sigma_2,\]
we can enumerate all the indices $y \in B$ such that $W_y \neq W_e$ using a
$\zero'$ oracle.  Eventually, we enumerate all of the indices except for one.
If the leftover index is $a$, then $W_e = W_a$, so $e \in \TOT$.  Otherwise, $e
\not\in \TOT$.  Thus, we can decide membership for a $\Pi_2$-complete set using
only a $\MIN \join \zero'$ oracle.
\end{proof}

\begin{proof}[(\ref{MIN*join0''is0'''})]
Schaefer's proof of $\fMIN \join \zero'' \equiv_\T \zero'''$ uses the fact that
there is an ordering $\phe$ such that $\fMIN^*_\phe \equiv_\T \zero'''$
\citep{Sch98}. The argument in (\ref{MINTjoin0'''is0''''}) with $\COF$
substituted for $\HIGH^n$ yields an analogous result, without taking into
consideration other G\"{o}del numberings.
\end{proof}

\begin{proof}[(\ref{MINmjoin0''is0'''})]
Define a $\compliment{\MIN^\m}$-computable function $f$ by \[ f(e) := (\mu i)
\left[i \in \MIN^\m \andd i > e \right].\]  Then \[(\forall e) \left[W_e
\not\equiv_\m W_{f(e)}\right].\]  Since $\compliment{\MIN^\m} \in \Sigma_3$, it
follows from the $\equiv_\m$-Completeness Criterion
(Theorem~\ref{jockuschcompleteness}\nlb(\ref{jockuschcompletenessm}),
\citep{JLSS}) that
\[\MIN^\m \join \zero'' \equiv_\T \compliment{\MIN^\m} \join \zero'' \equiv_\T \zero'''. \qedhere\]
\end{proof}

\begin{proof}[(\ref{MINTjoin0'''is0''''})]
Recall that $\mn^{\T^{(n)}}(e)$ denotes the function which computes the
$\equiv_{\T^{(n)}}$-minimal index of $e$.  We claim that
\[\mn^{\T^{(n)}} \leq_\T \MIN^{\T^{(n)}} \join \zero^{(n+3)}.\]
Let $a$ denote the $\T^{(n)}$-minimal index for $\zero^{(n+1)}$.  In
Theorem~\ref{basicarithmetic}\nlb(\ref{mintinpi4}), we showed
\[\left\{ \pair{j,e} : W_j \equiv_{\T^{(n)}} W_e \right\} \in \Sigma_{n+4}, \]
so we can enumerate the pairs of $\equiv_{\T^{(n)}}$-equivalent c.e$.$ sets
using a $\zero^{(n+3)}$ oracle.

For any index $e$, $W_e \equiv_{\T^{(n)}} W_x$ for exactly one $x$ in
\[\{0, \dotsc, e \} \intersect \MIN^{\T^{(n)}}. \]
Since a unique $x$ is guaranteed to exist, we have that $x =
\min^{\T^{(n)}}(e)$ can be computed from a $\MIN^{\T^{(n)}} \join
\zero^{(n+3)}$ oracle.  This proves the claim.

Now since
\[ \HIGH^n = \{ e : W_e \equiv_{\T^{(n)}} \zero' \} \]
is $\Sigma_{n+4}$-complete (see
Theorem~\ref{lowerbound}\nlb(\ref{MINTlowerbound})), it suffices to determine,
using a $\MIN^{\T^{(n)}} \join \zero^{(n+3)}$ oracle, whether a given index $e$
is in $\HIGH^n$.  To do this, just compute $\mn^{\T^{(n)}}(e)$, and check
whether it is equal to $a$.
\end{proof}

\end{proof}

Note that Lemma~\ref{downonelevel}\nlb(\ref{MINmjoin0''is0'''}) gives us
another way of showing that $\MIN^\m \in \Pi_3 - \Sigma_3$.  We obtain similar
arithmetic results from the other parts of Lemma~\ref{downonelevel}.

\section{(Old)-timers} Prior to this work, the only technique
which was successful in reducing a $\MIN$-set by a second ``level'' was to use
$\MIN$ queries to build a ``timer'' for the convergence of some function,
thereby turning an enumerable object into something computable.  Unlike the
technique of Lemma~\ref{downonelevel}, however, the ``timer'' method appears to
be peculiar to the equivalence relation under consideration.  We demonstrate
this method in Lemma~\ref{downtwolevel}.

The following theorem isolates the main idea of Lemma~\ref{downtwolevel}.
Reading this proof may help to remember the proof of $\MIN \geq_\bT \zero'$. We
defined $\MIN^\tm$ in Definition~\ref{MINtmdefn}.
\begin{thm}
$\MIN^\tm \geq_\bT \zero'$. \label{MINtmis0'}
\end{thm}

\begin{proof}
It suffices to determine whether $W_e = \emptyset$ using a $\MIN^\tm$ oracle.
Let $a$ be the minimal index of the function which diverges everywhere. For $j
\in \MIN^\tm - \{a\}$, define a function $s$ by
\[ s(j) := (\mu \pair{t,x}) \left[ \phe_{j,t} (x) \converge \right],\]
and let
\[S(i) := \max_{\substack{ j \leq i \\ (j \in \MIN^\tm - \{a\})}
} s(j).\] $S(e)$ is computable from $\MIN^\tm$, as $s$ converges everywhere on
its domain. Now either $e \in \MIN^\tm$, or else $\phe_e$ duplicates the
computation of some $\phe_j$, $j < e$.  That is, either $\phe_{e,t}(x)$
converges for some $\pair{t,x} \leq S(e)$, or else $W_e = \emptyset$.
\end{proof}
The reverse inequality for Theorem~\ref{MINtmis0'} is immediate, as $\MIN^\tm
\in \Sigma_1$ and hence $\MIN^\tm \leq_\m \zero'$.

Schaefer proved Lemma~\ref{downtwolevel} for $\fMIN$ and $\fMIN^*$, but a
similar proof works for both sets and functions.

\begin{lemma}[\citet{Sch98}] \label{downtwolevel}
\ \begin{thmenum}
\item $\MIN \geq_\bT \zero'$, \label{MINgeq0'}

\item $\MIN^* \join \zero' \geq_\T \zero''$. \label{MIN*join0'geq0''}
\end{thmenum}
\end{lemma}

\begin{proof}
\begin{proof}[(\ref{MINgeq0'})]
Let $e$ be an index.  We show how to decide whether $\phe_e(e) \converge$ with
a $\MIN$ oracle.  Using the $s$-$m$-$n$ Theorem, define a computable function
$f$ by
\begin{equation} \label{MINgeq0'eq}
\phe_{f(i)}(x) :=
\begin{cases}
1 & \text{if $\phe_{i,x}(i)\converge$,} \\
\diverge & \text{otherwise.}
\end{cases}
\end{equation}
 Now $e \in \K$ iff $W_{f(e)} \neq \emptyset$.  Since $\phe_{f(i)}(x)$
effectively counts the steps in computation $\phe_i(x)$, we can now proceed as
in Lemma~\ref{MINtmis0'}.

Let $a$ be the minimal index of the function which diverges everywhere.  Define
a function $s: \MIN - \{a\} \to \omega$ by
\[ s(j) := (\mu x) \left[ \phe_j (x) \converge \right],\]
and let
\[S(i) := \max_{\substack{j \leq i \\ (j \in \MIN - \{a\})}} s(j).\]
Since $\phe_{f(e)}$ agrees with some index in $\MIN \intersect \{0, \dotsc,
f(e)\}$, it must be the case that
\begin{align*}
W_{f(e)} \neq \emptyset &\iff W_{f(e)} \intersect \{0, \dotsc, S[f(e)]\} \neq \emptyset\ \\
 &\iff \phe_{e,S[f(e)]}(e)\converge.
\end{align*}
Since $S$ is computable in $\MIN$, we can decide $W_{f(e)} \neq \emptyset$.
\end{proof}

\begin{proof}[(\ref{MIN*join0'geq0''})]
Recall that $\TOT \equiv_\T \zero''$.  Since $\compliment{\TOT}$ is c.e$.$ in
$\zero'$, it suffices to enumerate $\TOT$ using a $\MIN^* \join \zero'$ oracle.
Define computable functions $f$ and $g$ by
\begin{align*}
\phe_{f(i)}(x) &:=
\begin{cases}
\bigl\langle x, (\mu s)\, (\forall y \leq x) \left[\phe_{i,s} (y) \converge
\right] \bigr\rangle& \text{if
such an $s$ exists,} \\
\diverge & \text{otherwise.}
\end{cases}
\\[2ex]
\phe_{g(i)}(x) &:=
\begin{cases}
\pi_2 [\phe_i(y)] & \text{if $(\exists y) \left[y \geq x \andd \phe_i(y) \converge \right]$} \\
\diverge & \text{otherwise.}
\end{cases}
\end{align*}
Let $a$ be the $=^*$-minimal index for the function which diverges everywhere.
Define
\begin{equation} \label{AisTOTforMIN*}
A := \left\{e : (\exists \pair{j,N})\, \left[j \in [\MIN^* - \{a\}] \intersect
\{0, \dotsc, f(e)\} \andd (\forall x) [\phe_{e, \max\{N, \phe_{g(j)}(x)\}} (x)
\converge]\right]\right\}.
\end{equation}
We claim:
\begin{enumerate}
\item $A$ is enumerable with a $\MIN^* \join \zero'$ oracle, and

\item $A$ = $\TOT$.
\end{enumerate}

Note that $W_j$ is infinite when $j \in \MIN^* -\{a\}$, which makes
$\phe_{g(j)}$ a total function.  The bracketed clause in \eqref{AisTOTforMIN*}
is therefore computable in $\MIN^* \join \zero'$, which proves (1).

If $e \in A$ then the universal clause in \eqref{AisTOTforMIN*} is satisfied,
so $e \in \TOT$.  Conversely, assume $e \in \TOT$.  Then $f(e) \in \INF$, so
$f(e)$'s $=^*$-minimal index is not $a$.  Let $j$ be the $=^*$-minimal index
for $f(e)$, choose $n$ large enough so that \[(\forall x
> n) \left[ W_j (x) = W_{f(e)}(x) \right],\]  and choose $N$ large enough so
that \[(\forall x \leq n) \left[\phe_{e,N}(x) \converge \right].\]  Then for
all $x$,
\[ \max \{N, \phe_{g(j)}(x) \} \geq \pi_2[\phe_{f(e)}(x)], \]
because $\pi_2[\phe_{f(e)}]$ is a nondecreasing function.  Hence \[(\forall
x)\: [\phe_{e, \max\{N, \phe_{g(j)}(x)\}} (x) \converge],\] so our selected
pair $\pair{j,N}$ exhibits that $e \in A$.
\end{proof}
\end{proof}

\section[GBS]{The Forcing Lowness Lemma}

We show how to ``drop'' $\MIN^{\T^{(n)}}$ by a second ``level.''
 Lemma~\ref{forcerecursiveorindependent} is easiest to digest when we recall
that $\LOW^0$ is the set of indices with computable domains. The lemma gives
slightly more than we need to prove the main theorem of this section, which is
Theorem~\ref{MINTjoin0''geq0'''}.  The argument in
Theorem~\ref{MINTjoin0''geq0'''} only depends on knowing the index
$a_\pair{k,n}(0)$, however the entire countable sequence $a_\pair{k,n}(0),
a_\pair{k,n}(1), \dotsc$, as well as uniformity in $n$, will be required for
Theorem~\ref{psimandTttcomplete}.

We state a simple version of \citep[Theorem 6.3]{Sac63b} by Sacks for use in
the next lemma.  Sacks does not explicitly mention uniformity in his original
proof, however Soare does \citep[Theorem VIII.3.1]{Soa87}.

\begin{thm}[Sacks Jump Theorem \citep{Sac63}] \label{sacksjumptheorem}
Let $B$ be any set, and let $S$ be c.e$.$ in $B'$ with $B' \leq_\T S$.  Then
there exists a $B$-c.e$.$ set $A$ with $A' \equiv_\T S$.  Furthermore, an index
for $A$ can be found uniformly from an index for $S$.
\end{thm}

\begin{lemma}[forcing lowness] \label{forcerecursiveorindependent}
There exists a ternary computable function $a_\pair{k,n}(i)$ such that for
every index $k$ and any number $i$, $W_{a_\pair{k,n}(i)} \leq_{\T^{(n)}} W_k$.
In particular, and furthermore:
\begin{thmenum}
\item $k \in \LOW^n \implies (\forall i)\: [a_\pair{k,n}(i) \in \LOW^n]$,
\label{forcercursiveorindependentcase1}

\item $k \not \in \LOW^n \implies (\forall i \neq j) \left[W_{a_\pair{k,n}(i)}
\incomp{\T^{(n)}}  W_{a_\pair{k,n} (j)}\right]$.
\label{forcercursiveorindependentcase2}
\end{thmenum}
In either case, $a_\pair{k,n}(i) \in \LOW^{n+1}$ for all $k$, $n$, and $i$.
\end{lemma}

\begin{proof}
This lemma is secretly \citep[Exercises VII.2.7 and VII.2.3]{Soa87}, in mild
disguise. Indeed, we shall combine finite injury (\citep{Fri57}, \citep{Muc56})
with standard permitting (\citep{Dek54}, \citep{Yat65}) by playing the
Friedberg-Muchnik strategy (\citep{Muc56}, \citep{Fri57}) under
$\paren{W_k}^{(n)}$.  Our construction follows \citep{Soa08}.

Given inputs $n$ and $k$, we show how to effectively find $\zero^{(n)}$-c.e$.$
sets $A_0, A_1, \dotsc$ so that $A_0 = (W_{a_\pair{k,n}(0)})^{(n)}$, $A_1 =
(W_{a_\pair{k,n}(1)})^{(n)}$, $\dotsc$ etc$.$ satisfy the conclusions of the
theorem.  If $n$ is nonzero, then we can subsequently (and uniformly) find
appropriate indices for c.e$.$ sets by iteratively applying the Sacks Jump
Theorem (Theorem~\ref{sacksjumptheorem}). For clarity purposes, we adopt the
following abbreviations:
\begin{align*}
B_i &:= \ijoin_{j \neq i} A_j, \\
\paren{B_i}_s &:= \ijoin_{j \neq i} \paren{A_j}_s,
\end{align*}
where $\paren{A_j}_0 \subseteq \paren{A_j}_1 \subseteq \dotsc$ is a
$\zero^{(n)}$-enumeration for $A_j$.

If $k \in \LOW^n$, our construction will satisfy for all $i$,
\begin{align*}
Q_i &: A_i \equiv_{\T^{(n)}} \zero,
\end{align*}
and if $k \not\in \LOW^n$, our construction will meet the requirements, for all
$i$ and $e$:
\begin{align*}
N_i &: A_i \leq_{\T^{(n)}} W_k, \\
R_\pair{e,i} &: A_i \neq \Psi_e^{B_i}.
\end{align*}

In the following construction, we imagine $Y$ to be the set $\zero^{(n)}$.  We
write $Y$ in place of $\zero^{(n)}$ simply to emphasize that our algorithm is
independent of the choice of oracle.  Furthermore, our construction will be
uniform in $k$.  Let \[ C_k := (W_k)^{(n)} \join \omega.\]  Now $C_k$ is c.e$.$
in $\zero^{(n)}$, and an index for $C_k$ (with $\zero^{(n)}$ oracle) can be
found uniformly from $k$.  'The ``$\omega$'' is added into the definition of
$C_k$ just to ensure that the set is infinite.  Since our construction will no
longer refer to the value $k$, we abbreviate with $C := C_k$.  Using the
$\zero^{(n)}$-index for $C$, we can effectively find a 1:1 function $c \leq_\T
\zero^{(n)}$ such that $c(0), c(1), c(2), \dotsc$ is an enumeration of $C$.

\cons
\begin{enumerate}
\item[\textit{Stage $s=0$}.] Define $r\left(\pair{e,i},0\right) =
-1$ for all $\pair{e,i}$.  Set $\paren{A_i}_0 = \zero \join Y$ for all~$i$.

\item[\textit{Stage $s+1$ ($s+1$ is an $i^\text{th}$ prime power)}.]
Choose the least $e$ such that
\begin{multline}\label{recursiveorindependenteq}
r\left(\pair{e,i},s\right) = -1 \andd (\exists \text{ even } x)\Bigl[x
\in\omega^{\left[\pair{e,i}\right]} - \left( A_i \right)_s  \andd
\Psi_{e,s}^{\paren{B_i}_s}(x)\converge = 0  \\
\andd \left(\forall \pair{z,j} < \pair{e,i}
\right)\left[r\left(\pair{z,j},s\right) < x \right] \andd c(s) \leq x \Bigr].
\end{multline}

If there is no such $e$, then do nothing and go to stage $s+2$.  If $e$ exists,
then we say $R_\pair{e,i}$ \emph{acts} at stage $s+1$, Perform the following
steps.
\begin{step}

\item Enumerate $x$ in $A_i$.

\item Define $r \left(\pair{e,i}, s+1 \right) = s + 1$.

\item For all $\pair{z,j} > \pair{e,i}$, define $r \left(\pair{z,j}, s+1 \right) =
-1$.

\item For all $\pair{z,j} < \pair{e,i}$, define $r\left(\pair{z,j}, s+1 \right) =
r\left(\pair{z,j}, s \right)$.
\end{step}
When $r(\pair{z,j},s+1)$ is reset to $-1$. we say that requirement
$R_\pair{z,j}$ is \emph{injured}.

\item[\textit{Stage $s+1$ ($s+1$ is not a prime power)}.] Do nothing.  Get some
coffee.
\end{enumerate}

\begin{claim} \label{recursiveorindependentclaim1}
For all $i$, $A_i \leq_\T C$.
\end{claim}

\begin{proof}
To decide whether $x \in A_i$, wait for a stage $s$ such that all the elements
of $C$ below $x+1$ have been enumerated into $C$, i.e., $$C \restr x \subseteq
\{c(0), c(1), \dotsc, c(s)\}.$$  Such a stage $s$ is guaranteed to exist, and
the oracle $C$ lets us identify when this occurs.  The final clause of
\eqref{recursiveorindependenteq}, ``$c(s) \leq x$,'' ensures that no element
$\leq x$ get enumerated into $A_i$ after stage $s$.  Hence \[x \in A_i \iff x
\in \paren{A_i}_{s+1}. \qedhere\]
\end{proof}

If $C \leq_\T \zero^{(n)}$, then by Claim~\ref{recursiveorindependentclaim1},
$A_i$ is $\zero^{(n)}$-computable for every $i$. This proves
case~(\ref{forcercursiveorindependentcase1}). It remains to consider
case~(\ref{forcercursiveorindependentcase2}).

\begin{claim} \label{recursiveorindependentclaim2}
If requirement $R_\pair{e,i}$ acts at some stage $s+1$ and is never later
injured, then requirement $R_\pair{e,i}$ is met and $r \left(\pair{e,i}, t
\right) = s + 1$ for all $t \geq s + 1$.
\end{claim}
\begin{proof}
Suppose $R_\pair{e,i}$ acts at stage $s+1$ and say $e$ is an $i^\text{th}$
prime power. Then \[\Psi^{\paren{B_i}_s}_e (x)\converge = 0\] for some $x \in
(A_i)_{s+1}$.  Since no $R_\pair{z,j}, \pair{z,j} < \pair{e,i}$ ever acts after
stage $s + 1$, it follows by induction on $t > s$ that $R_\pair{e,i}$ never
acts again and $r \left(\pair{e,i}, t \right) = s + 1$ for all $t > s$.  Hence
no $R_\pair{z,j}$, $\pair{z,j} > \pair{e,i}$, enumerates any $x \leq s$ into
any $A_j$ $(j \neq i)$ after stage $s + 1$. Therefore,
\[B_i \restr s = \paren{B_i}_s \restr s\] and
\[\Psi^{B_i}_e (x)\converge = 0 \neq A_i (x).\qedhere\]
\end{proof}

\begin{claim} \label{recursiveorindependentclaim3}
Assume $C >_\T \zero^{(n)}$. Then for every $\pair{e,i}$, requirement
$R_\pair{e,i}$ is met, acts at most finitely often, and
$r\left(\pair{e,i}\right) := \lim_s r\left( \pair{e,i}, s \right)$ exists.
\end{claim}
\begin{proof}
Fix $\pair{e,i}$ and assume the statement holds for all $R_\pair{z,j}$,
$\pair{z,j} < \pair{e,i}$.  Let $v$ be the greatest stage when some such
$R_\pair{z,j}$ acts, if ever, and $v = 0$ if none exists.  Then $r \left(
\pair{e,i}, v \right) = -1$, and this persists until some stage $s + 1 > v$ (if
ever) when $R_\pair{e,i}$ acts.  If $R_\pair{e,i}$ acts at some stage $s + 1$,
then $R_\pair{e,i}$ becomes satisfied and never acts again.  It then follows
from Claim~\ref{recursiveorindependentclaim2} that $r \left(\pair{e,i}, t
\right) = s + 1$ for all $t \geq s + 1$.

Either way, $r \left(\pair{e,i}\right)$ exists and $R_\pair{e,i}$ acts at most
finitely often.  Now suppose that $R_\pair{e,i}$ is not met.  Then \[A_i =
\Psi^{B_i}_e.\] By stage $v$, at most finitely many elements $x \in
\omega^{\left[\pair{e,i}\right]}$ have been enumerated in $A_i$.  No further
elements are enumerated from $\omega^{[\pair{e,i}]}$ because only requirement
$R_\pair{e,i}$ can enumerate in this row.  Let $x \in
\omega^{\left[\pair{e,i}\right]} - \paren{A_i}_v$ be such that $x
> v$. Eventually there will be a stage $s$ such that
\[\Psi_{e,s}^{\paren{B_i}_{s}}(x)\converge = 0,\] because $x \not\in A_i$.
Since $x$ never becomes a witness that $R_\pair{e,i}$ is satisfied, it must be
the permitting clause ``$c(s) \leq x$'' in \eqref{recursiveorindependenteq}
which prevents this from happening. Therefore
\[C \restr x = \{c(0), \dotsc, c(s)\} \restr x.\]
Since $x$ was chosen arbitrarily, we now have an algorithm to compute any
finite initial segment of $C$.  Our algorithm used only a $\zero^{(n)}$ oracle
to compute the function $c$.  Therefore $C \leq_\T \zero^{(n)}$, contrary to
assumption. So requirement $R_\pair{e,i}$ must be met.
\end{proof}

Case~(\ref{forcercursiveorindependentcase2}) is now satisfied because the
requirements $R_\pair{e,i}$ are met.  Finally,
\begin{claim}[\citet{Soa72}] \label{forcerecursiveorindependentclaimsoare}
For every $k$, $n$, and $i$, we have $a_\pair{k,n}(i) \in \LOW^{n+1}$.
\end{claim}
\begin{proof}
We may assume $C >_\T \zero^{(n)}$ because otherwise the result follows
immediately from Claim~\ref{recursiveorindependentclaim1}.  Using the
relativized $s$-$m$-$n$ theorem, define a computable function $f$ such that for
all $Y \subseteq \omega$,
\[\Psi_{f(e)}^Y (x) :=
\begin{cases}
0 &\text{if $\Psi_e^Y (e) \converge$,} \\
\diverge &\text{otherwise}.
\end{cases}\]
$\Psi_{f(e)}^Y$ is either the constant zero function or diverges everywhere,
depending on $Y$.  Define a computable ``witness'' function $w$ by
\[w(\pair{e,i},s) :=
\begin{cases}
\text{most recent member of $A_i \intersect \omega^{[\pair{e,i}]}$ after stage $s$, or} \\
\text{$\left\langle 0, \pair{e,i} \right\rangle$ if none exists.}
\end{cases}\]
Since each requirement acts only finitely often
(Claim~\ref{recursiveorindependentclaim3}), the limit
$$\hat{w}(e,i) := \lim_{s} w
\paren{\pair{e,i},s}$$ exists and witnesses
$\Psi_e^{B_i}[\hat{w}(e,i)] \neq A_i [\hat{w}(e,i)]$. Finally, define a
sequence of functions $g_i \leq_{\T^{(n)}} \zero$ by
\[g_i(e,s) :=
\begin{cases}
1 &\text{if $\Psi_{f(e),s}^{\paren{B_i}_s} \Bigl(w
\left[\pair{f(e),i},s\right] \Bigr)\converge = 0$,} \\
0 &\text{otherwise.}
\end{cases}\]
We show that
\begin{equation}
\hat{g}_i(e) := \lim_s g_i(e,s) \label{recursiveorindependenteq2}
\end{equation}
is the characteristic function for $(B_i)'$, which implies that $(B_i)'
\leq_{\T^{(n)}} \zero'$ by the Limit Lemma.

Let $t$ be a large enough stage so that $R_\pair{f(e),i}$ never gets injured
after stage $t$, and large enough so that $w(\pair{f(e),i}, \cdot)$ has
settled, i.e.
\[(\forall s > t)\: (w[\pair{f(e),i},s] = w[\pair{f(e),i},t] = \hat{w}[f(e),i]).\]
For clarity, let $\tilde{w}$ denote the value $\hat{w}[f(e),i]$, and let $v_s$
denote the function
$$v_s (x) := \Psi_{f(e),t}^{\paren{B_i}_s} (x).$$
Now for all $s > t$, $g_i(e,s) = g_i(e,t)$, so the limit in
\eqref{recursiveorindependenteq2} exists.  Indeed, if $v_t
\paren{\tilde{w}} \converge = 0$, and at some later stage $s$, $\neg \left[v_s
\paren{\tilde{w}} \converge = 0 \right]$, this would force our construction to find
a new witness for $R_\pair{e,i}$, contradicting the fact that $\tilde{w}$ is
the final witness.  If, on the other hand, $\neg \left[v_t
\paren{\tilde{w}} \converge = 0 \right]$, then this computation on $\tilde{w}$
must be preserved forever, lest $\R_\pair{e,i}$ acts again.

Since $\hat{g}_i(e) = g_i(e,t)$, it follows that
\begin{align*}
\hat{g}_i(e) = 1 &\bigiff \Psi_{f(e),t}^{\paren{B_i}_t}
\Bigl(\hat{w}[f(e),i] \Bigr) \converge = 0 \\
& \bigiff \Psi_{f(e)}^{B_i} \Bigl(\hat{w}[f(e),i] \Bigr) \converge = 0 \\
& \bigiff \Psi_e^{B_i}(e) \converge.
\end{align*}
Therefore $\hat{g}_i$ is the characteristic function for $(B_i)'$.  This proves
$a_\pair{k,n}(j) \in \LOW^{n+1}$ for all $j \neq i$, as $a_\pair{k,n}(j)$ is
the $\zero^{(n)}$-index for $A_j \leq_\T B_i$.  Since $i$ was chosen
arbitrarily, we conclude that, in fact, $a_\pair{k,n}(i) \in \LOW^{n+1}$ for
all $i \in \omega$.
\end{proof}
\end{proof}

Our first application of Lemma~\ref{forcerecursiveorindependent} is the
following theorem:

\begin{thm} \label{MINTjoin0''geq0'''}
$\MIN^{\T^{(n)}} \join \zero^{(n+2)} \geq_\T \zero^{(n+3)}$.
\end{thm}

\begin{proof}
Since $\LOW^n$ is $\Sigma_{n+3}$-complete, it suffices to determine membership
in $\LOW^n$ using a $\MIN^{\T^{(n)}} \join \zero''$ oracle.  On input $k$,
first compute $a_\pair{k,n}(0)$, where $a_\pair{k,n}$ is the computable
function defined in Lemma~\ref{forcerecursiveorindependent}, and let $c$ be the
least index such that
$$W_c \equiv_{\T^ {(n)}} \zero,$$
(i.e$.$, $c \in \LOW^n$).  We would like to know whether $\mn^{\T^{(n)}}(k) =
c$.

Let
\[e := a_\pair{k,n}(0),\]
and
\[S_e := \{0, \dotsc, e\} \intersect \MIN^{\T^{(n)}}.\]
There exists a unique $x \in S_e$ satisfying $W_{x} \equiv_{\T^{(n)}} W_e$,
however unlike in Theorem~\ref{downonelevel}\nlb(\ref{MINTjoin0'''is0''''}), we
can not discover which one it is by direct enumeration because we are now
missing the $\zero^{(n+3)}$ oracle.  So we use ``double enumeration'' instead.
Since $e \in \LOW^{n+1}$, the set
\[Y_e := S_e \intersect \left\{y : W_y \leq_{\T^{(n)}} W_e \right\}\]
is c.e$.$ in $\MIN^{\T^{(n)}} \join \zero^{(n+2)}$
(Proposition~\ref{basicarithmetic}\nlb(\ref{mintinpi4})). Let $Y_{e,t}$ denote
the elements which have been added into $Y_e$ after $t$ steps of this
enumeration.  We remark that $Y_{e,t} \leq_\T \MIN^{\T^{(n)}} \join
\zero^{(n+2)}$.

\begin{claim} \label{MINTdowntwolevelsclaim}
Define a function $Z$ from $\rng [a_\pair{\cdot,n}(0)]$ to finite sets by
\[Z(e) := Y_e \intersect \left\{y : W_e \leq_{\T^{(n)}} W_y \right\}.\]
Then
\begin{thmenum}
\item $Z \leq_\T \MIN^{\T^{(n)}} \join \zero^{(n+2)}$, and
\label{MINTdowntwolevelsclaim1}

\item $Z(e) = \{\mn^{\T^{(n)}}(e)\}$. \label{MINTdowntwolevelsclaim2}
\end{thmenum}
\end{claim}
\begin{proof}
\eqref{MINTdowntwolevelsclaim2} is immediate because $z \in Z(e)$ implies $W_z
\equiv_{\T^{(n)}} W_e$, and $\mn^{\T^{(n)}}(e)$ is the unique member of $S_e$
with this property. It remains to compute $Z(e)$ with a $\MIN^{\T^{(n)}} \join
\zero^{(n+2)}$ oracle. Note that when $y \in Y_{e,t}$, the relation
\begin{equation} \label{MINTdowntwolevelseq}
(\exists i \leq t)\: (\forall x) \left[\Psi_i^{(W_y)^{(n)}}(x)\converge \andd
\left( x \in (W_e)^{(n)} \iff \Psi_i^{\paren{W_y}^{(n)}}(x) = 1 \right) \right]
\end{equation}
is in $\Pi_1^{\zero^{(n+1)}} = \Pi_{n+2}$ because $y \in \LOW^{n+1}$. Therefore
knowing \emph{a priori} that we are considering only members of $Y_{e,t}$, we
can decide membership in \eqref{MINTdowntwolevelseq} using the $\zero^{(n+2)}$
oracle.

The algorithm for $Z$ is as follows.  Assume that we have not yet converged by
stage $t$.  For each $y \in Y_{e,t}$, we check using $\zero^{(n+2)}$ whether
$y$ satisfies \eqref{MINTdowntwolevelseq}.  If we find a $y \in Y_{e,t}$
satisfying \eqref{MINTdowntwolevelseq}, then we know $W_e \leq_{\T^{(n)}} W_y$,
hence $Z(e) = \{y\}$, so the algorithm terminates. Otherwise we proceed
similarly in stage $t+1$.  Eventually we will discover a $y \in Y_e$ satisfying
\eqref{MINTdowntwolevelseq}, namely $y = \mn^{\T^{(n)}}(e)$.

We have glossed over one important detail of our algorithm, namely whether or
not we can check for membership in \eqref{MINTdowntwolevelseq} \emph{uniformly
in $e$}.  In fact, we can.  In order to make the algorithm uniform in $e$, we
not only need to know that $(W_y)' \leq_{\T^{(n)}} \zero'$, but we also need to
know explicitly what the reduction is so that we can make the correct queries
to $\zero''$ (regarding \eqref{MINTdowntwolevelseq}).

Here are the missing details.  When we enumerate $y$ into $Y_e$, we
automatically obtain a witness for $W_y \leq_{\T^{(n)}} W_e$, namely the index
of this reduction. Using this witness, we can effectively find a second index
witnessing $(W_y)' \leq_{\T^{(n)}} (W_e)'$.  Finally, $e$ is a special set of
the form $a_\pair{\cdot,n}(0)$, and so
Claim~\ref{forcerecursiveorindependentclaimsoare} gives a recipe for deciding
membership in $(W_e)^{(n+1)}$ given $\zero^{(n+1)}$.
\end{proof}

By Lemma~\ref{forcerecursiveorindependent},
\begin{align*}
Z(e) = \{c\} &\bigiff \mn^{\T^{(n)}}(e) = c \\
&\bigiff a_\pair{k,n}(0) = e \in \LOW_n \\
&\bigiff k \in \LOW^n.
\end{align*}
Thus, membership in $\LOW^n$ is decidable in $\zero^{(n+2)} \join
\MIN^{\T^{(n)}}$.
\end{proof}

\section{Conclusion}

We summarize the main results of this chapter.

\begin{cor} \label{bestTuring}
\begin{thmenum}
\item $\fR \equiv_\bT \zero'$. \label{fRis0}

\item $\MIN \equiv_\T \zero''$, \label{MINis0''}

\item $\MIN^* \join \zero' \equiv_\T \zero'''$. \label{MIN*join0'is0'''}

\item $\MIN^\m \join \zero'' \equiv_\T \zero'''$. \label{MINmjoin0'is0'''}

\item $\MIN^{\T^{(n)}} \join \zero^{(n+2)} \equiv_\T \zero^{(n+4)}$. \label{MINTjoin0'is0''''}
\end{thmenum}
\end{cor}

\begin{proof}
\begin{proof}[(\ref{fRis0})]
Use the proof from Lemma~\ref{downtwolevel}\nlb(\ref{MINgeq0'}), but in
\eqref{MINgeq0'eq} make $f$ check for convergence on $0$  rather than $i$.
\end{proof}

\begin{proof}[(\ref{MINis0''}), (\ref{MIN*join0'is0'''})]
Combine Lemma~\ref{downonelevel} with Lemma~\ref{downtwolevel}.
\end{proof}

\begin{proof}[(\ref{MINmjoin0'is0'''})]
Lemma~\ref{downonelevel}.
\end{proof}

\begin{proof}[(\ref{MINTjoin0'is0''''})]
Combine Lemma~\ref{downonelevel} with Theorem~\ref{MINTjoin0''geq0'''}.
\end{proof}

\end{proof}

It would be interesting to know whether or not the $\zero'$, $\zero''$, or
$\zero^{(n+2)}$ oracle is necessary in any of the above reductions.
Theorem~\ref{thickpheTuringcomplete} shows, in a formal sense, that a positive
answer to this question will be difficult to prove.

%% file: immunity.tex
\chapter[Immunity]{Immunity and fixed points} \label{immunitychapter}
We discuss ``thinness'' of spectral sets.  Spectral sets are naturally sparse,
as weaker relations give rise to thinner spectral sets (for example, $\MIN^\T
\subsetneq \MIN$).  The notion of ``thinness'' is formally captured by
immunity. Based on the examples of $\MIN$, $\MIN^*$, and $\MIN^\T$, one might
be tempted to extrapolate that $\MIN$-sets which are higher in the arithmetic
hierarchy are also more immune. In general, however, arithmetic level turns out
to be a crude and inaccurate indicator of thinness.  It is even possible to
find a pair of $\MIN$-sets where the arithmetic level is higher in one set and
immunity is greater in the other. For example, $\MIN^\m \in \Pi_3$ and
$\MIN^{\thick\equiv_*} \in \Sigma_4 - \Pi_4$ (see
Section~\ref{thickarithmetics}), but the first set is $\Sigma_3$-immune while
the latter is only $\Sigma_2$-immune.

The theorems in this chapter provide an alternative method to
Chapter~\ref{arithmetic properties} for showing that $\MIN$-sets are
noncomputable.  We illustrate a connection between these methods and
generalizations of the Arslanov completeness criterion.

\begin{defn}
Let $\mathcal{C}$ be a family of sets.  A set is \emph{$\mathcal{C}$-immune} if
it is infinite and contains no infinite members of $\mathcal{C}$.  If
$\mathcal{C}$ is the class of c.e$.$ sets, then we write \emph{immune} in place
of $\mathcal{C}$-immune.
\end{defn}

For example, the $\Pi_1$ set of Kolmogorov random strings, $\RAND$
(Definition~\ref{kolmogorovdefn}), is immune \citep[Corollary 2.7.1]{LV97}.  In
fact, $\compliment{\RAND}$ is a natural example of a \emph{simple} set, being
infinite, c.e$.$, and having a complement which is immune.  Simple sets were
first invented by Emil Post in attempt to exhibit a c.e$.$ set $A$ satisfying
$\zero < A < \zero'$ \citep{Pos44}.  Post's program ultimately failed, however
Post's problem (and consequently his notion of immunity) shaped the focus of
computability theory in the 1940's and 1950's \cite{Soa87}.

\section{Lower spectral sets}
Marcus Schaefer \citep{Sch98} made the following observations with regards to
minimal functions, but the results translate easily into sets.   He attributes
the main idea of (\ref{MINimmune}) to Blum \citep[Theorem 3]{Blu67} and
(\ref{MIN*immune}) to John Case:
\begin{thm}[\citet{Sch98}]\label{RMINMIN*immune}
\begin{thmenum}
\item $\fR$ is immune. \label{Rimmune}

\item $\MIN$ is immune. \label{MINimmune}

\item $\MIN^*$ is $\Sigma_2$-immune. \label{MIN*immune}
\end{thmenum}
\end{thm}

The ideas from Theorem~\ref{RMINMIN*immune} will come in handy when we prove
the $\Pi_3$-Separation Theorem (Theorem~\ref{pi3separation}).

First, we consider the problematic relation $\equiv_1$.  One might be tempted
to modify this relation by identifying finite sets, but we decline to do this
here.  Consequently, $\MIN^{\equiv_1}$ contains a representative of each finite
size.  The reason for doing this is not only that this finite property makes
the $\Pi_3$-Separation Theorem sparkle, but also because ``finiteness'' is
essentially an unavoidable aspect of 1:1 equivalence on computable sets:

\begin{prop}[\citet{DM60}]\label{proper subsets of immune sets}
Let $A$ be an immune set of nonegative integers.  Then
\[A <_1 A \union \{-1\} <_1 A \union \{-1, -2\} <_1 A \union \{-1, -2, -3\} <_1 \dotsc\]
\end{prop}

\begin{proof}
Let
\begin{align*}
A_0 &:= A, \\
A_1 &:= A \union \{-1\}, \\
A_2 &:= A \union \{-1, -2\}, \\
\vdots
\end{align*}
Clearly, $A_n \leq_1 A_{n+1}$ via the identity function.  Suppose towards a
contradiction that $A_{n+1} \leq_1 A_n$ for some $n$, and let $f$ be the
computable function that witnesses this relation.  Let $x \in A_{n+1} - A_n$.
Then the sequence
\begin{align*}
f(x) \\
f \of f(x) \\
f \of f \of f(x) \\
\vdots
\end{align*}
has no repetitions, and is therefore an infinite c.e$.$ subset of $A_n$.
Indeed, $x$ is not in the range of $f$, so a repetition of values would
indicate that $f$ is not injective.  This means that $A_n$ is not immune, a
contradiction.
\end{proof}

The following theorem shows that immunity can be used to distinguish between
certain $\MIN$-sets, even when the arithmetic hierarchy can not.

\begin{thm}[$\Pi_3$-Separation] \label{pi3separation}
$\MIN^\m$, $\MIN^*$, and $\MIN^{\equiv_1}$ are all in $\Pi_3 - \Sigma_3$, but
\begin{thmenum}
\item $\MIN^\m$ is $\Sigma_3$-immune, whereas \label{pi3separationm}

\item $\MIN^*$ contains an infinite $\Sigma_3$ set, and \label{pi3separation*}

\item $\MIN^{\equiv_1}$ contains an infinite $\Sigma_2$ set. \label{pi3separation1}
\end{thmenum}
\end{thm}

\begin{proof}
We already showed $\MIN^\m, \MIN^*, \MIN^{\equiv_1} \in \Pi_3 - \Sigma_3$ in
Theorem~\ref{lowerbound}.

\begin{proof}[(\ref{pi3separationm})]
$\MIN^\m$ is infinite because it's noncomputable
(Theorem~\ref{lowerbound}\nlb(\ref{MINmlowerbound})).  Let $A$ be an infinite,
$\Sigma_3$ set, and suppose $A \subseteq \MIN^\m$. Since $A$ is infinite and
c.e$.$ in $\zero''$, we can define a $\zero''$-computable function $g$ by
\[g(e) = \pi_1 \Bigl( \left( \mu \pair{i,t} \right)\: [i > e \andd i \in A_t] \Bigr), \]
where $\{A_t\}$ is a $\zero''$-enumeration of $A$.

Now for all $e$, $g(e) > e$ and $g(e) \in \MIN^\m$.  Therefore
\[ (\forall e)\: [W_e \not\equiv_\m W_{g(e)}],\]
contradicting a theorem of Jockusch et al. (Theorem~\ref{JLSSmain}): for every
$f \leq_\T \zero''$,
\[(\exists e)\: [W_e \equiv_\m W_{f(e)}]. \qedhere\]
\end{proof}

\begin{proof}[(\ref{pi3separation*})]
For every $k$, let
\begin{align*}
P_k &:= \{ n : n \ \text{is a $k^\text{th}$ prime power} \}, \\
A_k  &:= \{e : W_e \subseteq^* P_k \} \intersect \INF, \\
A &:= \{e : (\exists k)\ (\forall j < e)\ [e \in A_k \enskip \& \enskip j
\not\in A_k ]\}.
\end{align*}
Now $A \subseteq \MIN^*$, as $e \in A$ implies $W_j \not =^* W_e$ for all $j <
e$. Since the $A_k$'s are disjoint, any infinite $B$ satisfies $B \subseteq^*
A_k$ for at most one $k$.  Moreover, each $A_k$ contributes a distinct element
to $A$, hence $A$ is infinite. Finally,
\begin{align*}
W_e \subseteq^* P_k &\iff (\exists N)\: (\forall x \geq N)\: [x \in W_e \implies x \in P_k] \\
&\iff (\exists N)\: (\forall x \geq N)\: [x \not\in W_e \; \orr \; x \in P_k] \\
&\iff (\exists N)\: (\forall x \geq N)\: (\forall t)\ [x \not\in W_{e,t} \;
\orr \; x \in P_k],
\end{align*}
which makes $A_k \in \Delta_3$, on account of $\INF \in \Pi_2$.  It follows
that $A \in \Sigma_3$.
\end{proof}

\begin{proof}[(\ref{pi3separation1})]
Define a sequence of finite sets by
\[A_k := \{x : 0 \leq x \leq k \}.\]
Furthermore, define
\[B_k := \{e : \text{$W_e$ has \emph{at least} $k$ elements}\} \in \Sigma_1,\]
which means that
\begin{equation*}
C_k := \{e : \text{$W_e$ has \emph{exactly} $k$ elements}\} = B_k \intersect
\compliment{B_{k+1}} \in \Delta_2.
\end{equation*}
It follows from the Pigeonhole Principle that
\[W_e \equiv_1 A_k \iff e \in C_k,\]
and therefore
\[\{ \pair{e,k} : W_e \equiv_1 A_k \} \in \Delta_2.\]
Now
\[A := \left\{e : (\exists k)\, (\forall j < e) \left[ W_j \not\equiv_1 A_k \andd W_e \equiv_1 A_k\right]\right\}\]
is a $\Sigma_2$ set.  Moreover, $A$ is infinite because each $A_k$ represents a
distinct $\equiv_1$ class.  Since $A \subseteq \MIN^{\equiv_1}$, it follows
that $\MIN^{\equiv_1}$ is not $\Sigma_2$-immune.
\end{proof}
\end{proof}

\begin{rem}
It is worth noting that $\MIN^{\equiv_1}$ is immune (simply because it is a
subset of $\MIN$).
\end{rem}
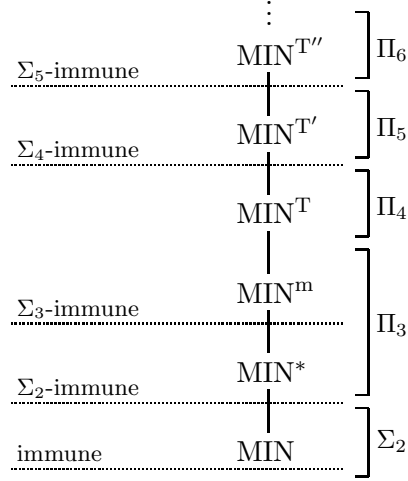
\begin{figure}[h] \label{naivefig}
\begin{picture}(100,165)
\jput(51,165){\vdots}
\jput(40,150){$\MIN^{\T''}$}
\jput(40,120){$\MIN^{\T'}$}
\jput(40,90){$\MIN^\T$}
\jput(40,60){$\MIN^\m$}
\jput(40,30){$\MIN^*$}
\jput(40,0){$\MIN$}

\drawline(52,131)(52,147)
\drawline(52,101)(52,117)
\drawline(52,71)(52,87)
\drawline(52,41)(52,57)
\drawline(52,11)(52,27)

\jput(-43,145){\footnotesize{$\Sigma_5$-immune}}
\dottedline{1.5}(-45,142)(80,142)
\jput(-43,115){\footnotesize{$\Sigma_4$-immune}}
\dottedline{1.5}(-45,112)(80,112)
\jput(-43,55){\footnotesize{$\Sigma_3$-immune}}
\dottedline{1.5}(-45,52)(80,52)
\jput(-43,25){\footnotesize{$\Sigma_2$-immune}}
\dottedline{1.5}(-45,22)(80,22)
\jput(-43,0){\footnotesize{immune}}
\dottedline{1.5}(-45,-3)(80,-3)

\jput(93,153){\small{$\Pi_6$}}
\drawline(85,145)(90,145)(90,170)(85,170)
\jput(93,123){\small{$\Pi_5$}}
\drawline(85,115)(90,115)(90,140)(85,140)
\jput(93,93){\small{$\Pi_4$}}
\drawline(85,85)(90,85)(90,110)(85,110)
\jput(93,49){\small{$\Pi_3$}}
\drawline(85,25)(90,25)(90,80)(85,80)
\jput(93,5){\small{$\Sigma_2$}}
\drawline(85,-6)(90,-6)(90,20)(85,20)
\end{picture}
\caption{A n\"{a}ive approach to spectral sets, by reverse inclusion.}
\end{figure}
The set inclusions and relations in Figure~\ref{naivefig} becomes nonlinear
when we add in spectral sets such as $\fR$, $\fMIN$, $\MIN^{\equiv_1}$,
$\MIN^{\thick*}$ and $\MIN^{\aev\T}$ (see Chapter~\ref{thickville} and
Appendix~\ref{almost thickness}).  $\MIN^{\thick*}$, in particular, does not
seem to fit into this picture at all.  Indeed, $\MIN^{\thick*} \in \Sigma_4 -
\Pi_4$ (Theorem~\ref{thicklowerbound}) but is only $\Sigma_2$-immune
(Theorem~\ref{thickimmune}). A simple, general pattern connecting arithmetics
and immunity does not seem to exist.


\section{Upper spectral sets}

The goal of this section is to determine the immunity of $\MIN^{\T^{(n)}}$.  In
the following theorem, the cases $=^*$ and $\equiv_\T$ were first proved by
Arslanov, and the remaining cases are due to Jockusch et al.
\begin{thm}[generalized fixed points, \citet{Ars81}, \citet{JLSS}]\label{JLSSmain}
For every $n \leq \omega$,
\begin{thmenum}
\item $f \leq_\T \zero' \implies (\exists e)\ [W_e =^* W_{f(e)}]$,

\item $f \leq_\T \zero'' \implies (\exists e)\ [W_e \equiv_\m W_{f(e)}]$,
\label{JLSSmainm}

\item $f \leq_\T \zero^{(n+2)} \implies (\exists e)\ [W_e \equiv_{\T^{(n)}} W_{f(e)}]$.
\end{thmenum}
Furthermore, $e$ can be found effectively from $n$ and an index for $f$.
\end{thm}

\begin{cor}\label{mintnimmune}
For all $n < \omega$, $\MIN^{\T^{(n)}}$ is $\Sigma_{n + 3}$-immune.
\end{cor}
\begin{proof}
We follow the proof of the $\Pi_3$-Separation Theorem
(Theorem~\ref{pi3separation}\nlb(\ref{pi3separationm})), and as before,
$\MIN^{\T^{(n)}}$ is infinite (by
Theorem~\ref{lowerbound}\nlb(\ref{MINTlowerbound})).

Let $n \geq 0$, and let $A$ be an infinite, $\Sigma_{n+3}$ set.  Suppose $A \subseteq
\MIN^{\T^{(n)}}$.  Since $A$ is infinite and c.e$.$ in $\zero^{(n+2)}$, we can define a
$\zero^{(n+2)}$-computable function $g$ by
\[g(e) = \pi_1 \Bigl( \left( \mu \pair{i,t} \right)\: [i > e \andd i \in A_t] \Bigr), \]
where $\{A_t\}$ is a $\zero^{(n+2)}$-enumeration of $A$.

Now for all $e$, $g(e) > e$ and $g(e) \in \MIN^{\T^{(n)}}$.  Therefore
\[ (\forall e)\: [W_e \not\equiv_{\T^{(n)}} W_{g(e)}],\]
contradicting Theorem~\ref{JLSSmain}.
\end{proof}

We now show that Corollary~\ref{mintnimmune} is optimal.

\begin{defn}
A set $A$ is called \emph{low} if $A' \equiv_\T \zero'$.
\end{defn}

\begin{defn}
Let $A$ and $B$ be c.e$.$ sets.  $A$ and $B$  are \emph{pairwise minimal} if
\begin{defenum}
\item $A, B >_\T \zero$, and
\item For every c.e$.$ set $C$, \[[C \leq_\T A \andd C \leq_\T B]\enskip \implies \enskip C
\text{ is computable}.\]
\end{defenum}
\end{defn}

The original minimal pairs construction is due to Lachlan \citep{Lac66} and
Yates \citep{Yat66}.  We generalize their result as follows:

\begin{thm}\label{lowminpairs}
There exists a computable sequence of c.e$.$ sets $A_0, A_1, \dotsc$ which are
low and pairwise minimal.
\end{thm}

\begin{proof}
We shall assume familiarity with Lachlan's tree construction for minimal pairs
\citep{Lac73} as given in \citep{DH06}, since only a minor modification is
needed to prove the theorem.  Lachlan's tree construction meets the following
requirements:
\begin{align*}
R_e &:\ \compliment{A} \neq W_e, \\
Q_e &:\ \compliment{B} \neq W_e, \\
N_\pair{i,j} &:\ [\Psi_i^A = \Psi_j^B = f \andd f \in \TOT]\enskip \implies
\enskip f \text{ is computable}.
\end{align*}

Since we are constructing a sequence of noncomputable sets, we replace $R_e$
and $Q_e$ with an appropriate requirement $R_\pair{e,k}$.  We also add the
lowness requirement $L_\pair{k,e}$ from \citep[Theorem VII.1.1]{Soa87}:
\begin{align*}
R_\pair{k,e} &:\ \compliment{A_k} \neq W_e, \\
L_\pair{k,e} &:\ (\exists^\infty s) \left[\Psi_{e,s}^{\paren{A_k}_s}(e) \converge \right] \enskip \implies \enskip \Psi_e^{A_k}(e) \converge, \\
N_\pair{i,j,m,n} &:\ [\Psi_i^{A_m} = \Psi_j^{A_n} = f \andd f \in \TOT]\enskip
\implies \enskip f \text{ is computable}.
\end{align*}

Requirements $R_\pair{k,e}$ and $N_\pair{i,j,m,n}$ are satisfied in exactly the
same way as the original construction, once we place these requirements on the
analogous levels of the tree.  The lowness requirement, $L_\pair{k,e}$,
combines easily with the the $N_\pair{i,j,m,n}$ requirement because both are
negative requirements which only try to protect existing computations.  We
satisfy $L_\pair{k,e}$ simply by adding an extra constraint on the witnesses
chosen to satisfy $R_\pair{k,e}$: define a computable restraint function $r$ by
\[r(k,e,s) := \psi_{e,s}^{\paren{A_k}_s}(e), \]
where $\psi$ denotes the use function.  Then restrain, in stage $s$ with
priority $\pair{k,e}$ (lower numbers having higher priority), any element less
than $r(k,e,s)$ from entering $\paren{A_k}_{s+1}$.  In some stage t, after
finite many injuries, $r$ eventually protects the computation on $e$ with
oracle $\paren{A_k}_t$ (whether or not the computation converges), thereby
satisfying $L_\pair{k,e}$.
\end{proof}

Theorem~\ref{lowminpairs} easily relativizes:

\begin{thm} \label{genminpair}
For every $n$, there exists a computable sequence of c.e$.$ sets $A_0, A_1,
\dotsc$ such that for all $C$ c.e$.$ in $\zero^{(n)}$ and $i \neq j$,
\begin{thmenum}
\item $\zero^{(n)} <_\T (A_i)^{(n)}$.

\item $(A_i)^{(n+1)}  \equiv_\T \zero^{(n+1)}$,

\item $\left[C \leq_\T (A_i)^{(n)} \quad \& \quad C \leq_\T (A_j)^{(n)} \right] \implies C \leq_\T
\zero^{(n)}$.\label{genminpairmincond}
\end{thmenum}
\end{thm}

\begin{proof}
We perform the construction from Theorem~\ref{lowminpairs} above $\zero^{(n)}$,
but with each c.e$.$ set $A_k$ replaced with a respective set $B_k$ which is
c.e$.$ in $\zero^{(n)}$. $n$ applications of the Sacks Jump Theorem
(Theorem~\ref{sacksjumptheorem}) then gives the desired reduction:
${(A_k)}^{(n)} \equiv_\T B_k$.
\end{proof}

Theorem~\ref{genminpair} will be useful in the proof of
Theorem~\ref{minTnotimmune}.

\begin{thm}\label{minTnotimmune}
For all $n \geq 0$, $\MIN^{\T^{(n)}}$ is not $\Sigma_{n+4}$-immune.
\end{thm}
\begin{proof}
Let $n \geq 0$, and let $A_0, A_1, \dotsc$ be the corresponding sequence of
sets obtained from Theorem~\ref{genminpair}.  Recall that
\[\LOW^n := \{e : (W_e)^{(n)} \equiv_\T \zero^{(n)} \},\]
and define
\begin{gather*}
B_k := \biggl[ \{ x : W_x \leq_{\T^{(n)}} A_k \}\ \intersect\ \compliment{\LOW^n} \biggr], \\
B := \left\{ e : (\exists k)\:(\forall j < e)\ [e \in B_k\quad \&\quad j \not\in B_k ]
\right\}.
\end{gather*}
The proof of Proposition~\ref{basicarithmetic}(\ref{mintinpi4}) mentions that
$B \leq_\T A$ is a $\Sigma_2^{B \join A'}$ relation.  Since, for any $x$, both
$(W_x)^{(n)} \leq_\T \zero^{(n+1)}$ and $(A_k)^{(n+1)} \leq_\T \zero^{(n+1)}$,
it follows that \[\left\{ x : (W_x)^{(n)} \leq_\T (A_k)^{(n)} \right\} \in
\Sigma_2^{\zero^{(n+1)}} = \Sigma_{n+3}.\] This places $B_k \in \Delta_{n+4}$,
on account of $\compliment{\LOW^n} \in \Pi_{n+3}$.  Therefore $B \in
\Sigma_{n+4}$.

It remains to show that $B$ is an infinite subset of $\MIN^{\T^{(n)}}$.  Note
that $B_i \intersect B_j = \zero$ for $i \neq j$.  Indeed, if $e \in B_i
\intersect B_j$, then
\[W_e \leq_{\T^{(n)}} A_i \quad \& \quad W_e \leq_{\T^{(n)}} A_j \andd e \not\in \LOW^n,\]
contradicting Property~(\ref{genminpairmincond}) of Theorem~\ref{genminpair}.
Now since $B_k \neq \emptyset$ and each $B_k$ contributes exactly one element
to $B$, $B$ must be infinite.

Finally, assume $e \in B$, and let $k$ be such that $e \in B_k$ and  $j \not \in B_k$ for
all $j < e$.  Then for $j < e$,
\[W_e \leq_{\T^{(n)}} A_k \andd W_j \not\leq_{\T^{(n)}} A_k,\]
which implies $W_e \not\equiv_{\T^{(n)}} W_j$.  So $e \in \MIN^{\T^{(n)}}$.  That
is, $B \subseteq \MIN^{\T^{(n)}}$.
\end{proof}

\section{Completeness criterion} \label{completeness criterion}

As a application of our immunity results, we obtain generalizations of the
Arslanov Completeness Criterion.  The classical theorem and Arslanov's original
generalization can be stated as follows:

\begin{thm}[$=$ and $=^*$-Completeness Criterion, {\citet{Ars85}}]
\label{arslanovcompleteness}
\begin{thmenum}
\item Let $A$ be c.e.  Then
\[A \equiv_\T \zero' \iff (\exists f \leq_\T A)\, (\forall e) \left[ W_e \neq W_{f(e)} \right]. \]

\item Let $A \in \Sigma_2$ and $\zero' \leq_\T A$.  Then
\[A \equiv_\T \zero'' \iff (\exists f \leq_\T A)\, (\forall e) \left[W_e \neq^*
W_{f(e)}\right].\]
\end{thmenum}
\end{thm}
The forward directions of Theorem~\ref{arslanovcompleteness} follow immediately
from the fact that $\MIN$ is not $\Sigma_2$-immune and $\MIN^*$ is not
$\Sigma_3$-immune (recall $\MIN \in \Sigma_2$, and
Theorem~\ref{pi3separation}\nlb(\ref{pi3separation*})). Proofs are analogous to
Corollary~\ref{equiv1completeness} below.  According to \citep{JLSS}, the
hypothesis ``A is c.e$.$ (resp$.$ $\Sigma_2$)'' in
Theorem~\ref{arslanovcompleteness} can be strengthened to ``$A$ is $k$-c.e$.$
(resp$.$ $k$-c.e$.$ in $\zero'$).''  A \emph{$k$-c.e$.$ set} is a limit
computable set in which the function witnessing this fact never changes its
mind more than $k$ times.  By ``$A$ is $k$-c.e$.$,'' we mean that the theorem
holds for any $k$ (in fact, we can assume only $k$-$\REA$, see \citep{JS84} for
a definition).

Using immunity properties of $\MIN$-sets, we are able to give a completeness
criterion for 1:1 equivalence:

\begin{cor}[$\equiv_1$-Completeness Criterion] \label{equiv1completeness}
Let $A$ be c.e.  Then
\[A \equiv_\T \zero' \iff (\exists f \leq_\T A)\, (\forall e) \left[W_e
\not\equiv_1 W_{f(e)}\right].\]
\end{cor}

\begin{proof}
Let $A \equiv_\T \zero'$, and suppose
\[(\forall f \leq_\T A)\, (\exists e) \left[W_e \equiv_1 W_{f(e)}\right].\]
This contradicts the $\Pi_3$-Separation Theorem (Theorem~\ref{pi3separation}):
under this assumption, the immunity argument in
Theorem~\ref{pi3separation}\nlb(\ref{pi3separationm}) shows that
$\MIN^{\equiv_1}$ is $\Sigma_2$ immune, but part~(\ref{pi3separation1}) of that
theorem says that it isn't.

The reverse direction is a direct application of the $=$-Completeness Criterion
(Theorem~\ref{arslanovcompleteness}).  Let $A$ be c.e$.$, and assume
\[(\exists f \leq_\T A)\, (\forall e) \left[W_e \not\equiv_1 W_{f(e)}\right].\]
Then clearly this same assertion holds for equality:
\[(\exists f \leq_\T A)\, (\forall e) \left[W_e \neq W_{f(e)}\right],\]
which implies that $A \equiv_\T \zero'$.
\end{proof}

The following result was first, along with its converse (modulo an appropriate
assumption about $A$), was first proved for $n = 0$ by Arslanov
\citep[Corollary 2.3]{Ars85} and for $n
> 0$ by Jockush et al$.$ \citep[ Corollary~5.17]{JLSS}.

\begin{cor}[necessary generalized completeness]\label{genarslanov}
Let $A$ be a set.  Then for any $n$,
\[ A \equiv_\T \zero^{(n+3)} \implies (\exists f \leq_\T A)\, (\forall e) \left[ W_e \not\equiv_{\T^{(n)}} W_{f(e)} \right].\]
\end{cor}

\begin{proof}
Let $A \equiv_\T \zero^{(n+3)}$, and suppose
\[(\forall f \leq_\T A)\, (\exists e) \left[ W_e \equiv_{\T^{(n)}} W_{f(e)} \right].\]
The argument in Corollary~\ref{mintnimmune} shows that $\MIN^{\T^{(n)}}$ is
$\Sigma_{n+4}$-immune for all $n$, contradicting Theorem~\ref{minTnotimmune}.
\end{proof}

Because $\MIN^{\m^{(n+1)}} = \MIN^{\T^{(n)}}$ (Proposition~\ref{minm'isminT}),
Corollary~\ref{genarslanov} shows that Theorem~\ref{JLSSmain} is optimal in the
sense that $\equiv_{\T^{(n+1)}}$ cannot be replaced with $\equiv_{\m^{(n+1)}}$.
In contrast, the authors of \citep{JLSS} note that $\equiv_\T$ \emph{can} be
substituted with $\equiv_\m$.  The converse for Corollary~\ref{genarslanov} is
known to hold when $A$ is $k$-c.e$.$ in $\zero^{(n+2)}$ (or $k$-REA in
$\zero^{(n+2)}$) for some $k$, and $\zero^{(n+2)} \leq_\T A$. Furthermore, this
additional condition is necessary \citep{JLSS}.  For ``completeness,'' we state
these theorems explicitly for the $\Sigma_n$ sets:

\begin{thm}[$\equiv_m$ and $\equiv_{\T^{(n)}}$-Completeness Criterion, {\citet{JLSS}}]
\label{jockuschcompleteness} \
\begin{thmenum}
\item Let $A \in \Sigma_3$ and $\zero'' \leq_\T A$.  Then
\[A \equiv_\T \zero''' \iff (\exists f \leq_\T A)\, (\forall e) \left[ W_e \not \equiv_\m W_{f(e)} \right].
\] \label{jockuschcompletenessm}

\item Let $A \in \Sigma_{n+3}$ and $\zero^{(n+2)} \leq_\T A$.  Then
\[A \equiv_{\T^{(n)}} \zero^{(n+3)} \iff (\exists f \leq_\T A)\, (\forall e) \left[W_e \not \equiv_{\T^{(n)}}
W_{f(e)}\right].\] \label{jockuschcompletenessT}
\end{thmenum}
\end{thm}

In summary, fixed points give us immunity:

\begin{thm} Let $n \geq 0$. Then \label{immunityfixedpoints}
\[ \left (\forall f \leq_\T \zero^{(n)} \right) (\exists e) \left[ W_e \equiv_\alpha W_{f(e)} \right]
\implies \MIN^{\equiv_\alpha}\ \text{is $\Sigma_{n+1}$-immune}.\]
\end{thm}

A simple converse to Theorem~\ref{immunityfixedpoints}, however, may not be
forthcoming. In Lemma~\ref{semi-fixed}, we develop the notion of
\emph{semi-fixed points}, which are sufficient to ensure certain immunity
properties (see Theorem~\ref{thickimmune}\nlb(\ref{thickMINmimmune})). In fact,
Theorem~\ref{thickimmune}\nlb(\ref{thickMINmimmune}) holds even if $\nu$ were
not computable, but merely $\nu \leq_\T \zero''$.  Thus it appears that fixed
points are a strictly stronger notion than immunity for $\MIN$-sets.

\section{Refinements}
\subsection{$\omega$-immunity} \label{omega immune}

Let $D_0$, $D_1$, $\dotsc$ be a computable numbering of the finite sets.

\begin{defn}[\citet{FS99}]
\begin{thmenum}
\item A set $A$ is called \emph{$k$-immune} if it is infinite and there is no computable function
$f$ such that
\begin{enumerate}
\item $\paren{D_{f(n)}}_{n \in \omega}$ is a family of pairwise disjoint sets,

\item $D_{f(n)} \cap A \neq \emptyset$, and

\item $\abs{D_{f(n)}} \leq k$.
\end{enumerate}

\item A set is called \emph{$\omega$-immune} if is $k$-immune for every $k$.
\end{thmenum}
\end{defn}

\begin{thm}[\citet{FS99}]
Let $A$ be a set. \label{allomegaimmunemain}
\begin{thmenum}
\item $\MIN$ is $\omega$-immune.

\item $\RAND$ is $\omega$-immune.

\item If $A$ is $\omega$-immune, then $\zero'\not\leq_\btt A$. \label{allomegaimmune}
\end{thmenum}
\end{thm}
Schaefer actually proves that $\fMIN$ is immune, but the same proof works for
$\MIN$ \citep{Sch98}.

\begin{cor}\label{notbtttozero'}
Let $A \subseteq \MIN$.  Then $\zero' \not\leq_\btt A$.  This includes all the
$\MIN$-sets  mentioned in this paper.
\end{cor}
\begin{proof}
Any subset of $\MIN$ is also $\omega$-immune.
\end{proof}

\subsection{$\Pi_n$-immunity}
The $\Pi_3$-Separation Theorem~\ref{pi3separation} gives us an optimal immunity
result for $\MIN^\m$, but the analogous theorems for $\MIN$, $\MIN^*$, and
$\MIN^\T$ leave room for improvement.  We can say a bit about $\Pi_n$ subsets
in general, if they exist.

\begin{thm}Let $n \geq 0$. \label{Pisubset}
\begin{thmenum}
\item Let $A$ be an infinite $\Pi_1$ subset of $\MIN$.  Then $A \equiv_\bT
\zero'$, $A$ is not hyperimmune, and $A \not\equiv_\btt \zero'$.
\label{MINPisubset}

\item Let $A$ be an infinite $\Pi_2$ subset of $\MIN^*$ such that $A \geq_\T \zero'$.  Then $A \equiv_\T \zero''$, but $A \not\equiv_\btt
\zero''$. \label{MIN*Pisubset}

\item Let $A$ be an infinite $\Pi_{n+3}$ subset of $\MIN^{\T^{(n)}}$ such that
 $A \geq_\T \zero^{(n+2)}$.
 Then $A \equiv_\T \zero^{(n+3)}$, but $A \not\equiv_\btt \zero^{(n+3)}$.
 \label{MINTPisubset}
\end{thmenum}
\end{thm}

\begin{proof}
\begin{proof}[(\ref{MINPisubset})]
Suppose $\MIN$ has an infinite subset $A \in \Pi_1$.  Since $\MIN$ is strongly
effectively immune \citep{Owi86}, $A$ must also be strongly effectively immune.
Thus $\compliment{A}$ is effectively simple, and it follows immediately that $A
\equiv_\bT \zero'$ \citep{Mar66}.  Furthermore, since a hypersimple set can
never be $\bT$-complete \citep{FR59}, $A$ must not be hyperimmune.  Finally,
$\compliment{A} \not\equiv_\btt \zero'$ follows from the fact that $A$ is
simple \citep[Theorem III.8.8]{Odi89}, \citep{pos44}.
\end{proof}

\begin{proof}[(\ref{MIN*Pisubset})]
The argument is quite the same as part~(\ref{MINPisubset}), but now we turn to
our $\omega$-immunity results.  Suppose $\MIN^*$ has an infinite subset $A \in
\Pi_2$.  Define $f \leq_\T \compliment{A}$ by
\[f(e) := (\mu x)\: [x \in A \andd x > e]. \]
Then
\[ (\forall e)\: [W_{f(e)} \neq^* W_e]. \]
Therefore $\compliment{A} \equiv_\T \zero''$ by the $=^*$ Completeness
Criterion (Theorem \ref{arslanovcompleteness}).

In place of simplicity, $A \subseteq \MIN^*$ implies that $A$ is
$\omega$-immune, and therefore $\compliment{A} \not\geq_\btt \zero'$  by
Theorem~\ref{allomegaimmunemain}.
\end{proof}

\begin{proof}[(\ref{MINTPisubset})]
same as part~(\ref{MIN*Pisubset}).
\end{proof}
\end{proof}

\subsection{$\Delta_n$-immunity}\label{Deltaimmune}

We observed in the introduction to Chapter~\ref{arithmetic properties} that
$\fR \subseteq \fMIN$, which proves $\fMIN$ is not $\Delta_2$-immune.  We now
show that $\MIN$ is not $\Delta_2$-immune either, although our witness to this
fact will not be a spectral set.

\begin{defn}
\begin{defenum}
\item Let $A \subseteq \omega$.  We say that a set $A$ is a \emph{partial function} if
\[(\forall x) \left[\pair{x,y_1}, \pair{x,y_2} \in A \enskip \implies \enskip y_1 = y_2 \right]. \]

\item Define the $\Pi_1$ set
\[ \FUN := \{e : W_e\ \text{is a partial function}\}. \]
\end{defenum}
\end{defn}

\begin{prop}
$\FUN$ is $\Pi_1$-complete.
\end{prop}

\begin{proof}
Using the $s$-$m$-$n$ Theorem, define a computable function $f$ by
\[\phe_{f(e)}(x) :=
\begin{cases}
1 & \text{if $x = e = \pair{\pi_1(e), \pi_2(e)}$ and $\phe_e (e) \converge$,} \\
1 & \text{if $x = \pair{\pi_1(e), \pi_2(e) + 1}$,} \\
\diverge &\text{otherwise.}
\end{cases} \]
Then \[e \in K \iff f(e) \not\in \FUN. \qedhere\]
\end{proof}

\begin{nota}
If $e \in \FUN$, we use $\hat{e}$ to denote the function represented by $W_e$.
In more detail, $\hat{e}(x)$ is the unique integer $y$ satisfying $\pair{x,y}
\in W_e$ if such a $y$ exists, and $\hat{e}(x)$ diverges otherwise.
\end{nota}

\begin{defn}
\[\FUNfR := \{e: e \in \FUN \andd \left( \forall j \in \FUN \intersect \{0, \dotsc, e \} \right) [\hat{j}(0) \neq \hat{e} (0)]\}.\]
\end{defn}

In the following Corollary, the general techniques from
parts~(\ref{MIN*notDelta3immune}) and (\ref{MINTnotDelta4immune}) can be
applied to part~(\ref{MINnotDelta2immune}).  We choose to have $\FUN$, however,
because $\FUN$ hypostatizes an essential connection between sets and functions.

\begin{cor}Let $n \geq 0$. \label{notDeltaimmune}
\begin{thmenum}
\item $\MIN$ is not $\Delta_2$-immune. \label{MINnotDelta2immune}

\item $\MIN^*$ is not $\Delta_3$-immune. \label{MIN*notDelta3immune}

\item $\MIN^{\T^{(n)}}$ is not $\Delta_{n+4}$-immune.
\label{MINTnotDelta4immune}
\end{thmenum}
\end{cor}

\begin{proof}
\begin{proof}[(\ref{MINnotDelta2immune})]
We show that $\FUNfR \subseteq \MIN$.  First, note that $\FUNfR \in \Delta_2$.
Indeed, $\FUN \leq_\T \zero'$, and then convergence of $\hat{j}(0)$ and
$\hat{e}(0)$ can be decided by asking $\zero'$ whether the following sets are
nonempty:
\begin{gather*}
\{y : \pair{0,y} \in W_j \}, \\
\{y : \pair{0,y} \in W_e \}.
\end{gather*}
Since there are infinitely many possible values for $\hat{e}(0)$, $\FUNfR$ must
be infinite.

Finally, it is clear that $\MIN$ contains $\FUNfR$:
\begin{align*}
e \in \FUNfR &\implies (\forall j \in \FUN \intersect \{0, \ldots, e\})
\left[\hat{j}(0) \neq \hat{e}(0) \right] \\
&\implies (\forall j < e)\: [W_j \neq W_e]. \qedhere
\end{align*}
\end{proof}

\begin{proof}[(\ref{MIN*notDelta3immune})]
By Theorem~\ref{pi3separation}, $\MIN^*$ contains an infinite $\Sigma_3$ set
$A$.  Since $A$ is c.e$.$ in $\zero''$, $A$ contains a set $B \leq_\T \zero''$,
namely $B = \{b_0 < b_1 < b_2 < \dotsb \}$ where
\begin{align*}
b_0 &:= \text{any member of A, and} \\
b_n &:= \pi_1 \left((\mu \pair{x,t})\: [x \in A_t \andd x > b_{n-1}] \right).
\end{align*}
Thus $B$ is an infinite, $\Delta_3$ subset of $\MIN^*$.
\end{proof}

\begin{proof}[(\ref{MINTnotDelta4immune})]
Similar to part~(\ref{MIN*notDelta3immune}), Theorem~\ref{minTnotimmune}
provides an infinite $\Sigma_{n+4}$ subset of $\MIN^{\T^{(n)}}$, which in turn
contains a $\Delta_{n+4}$ subset.
\end{proof}
\end{proof}

%% file: thickville.tex
\chapter[Thickville]{Thickville: nonuniformity vs. the jump operator} \label{thickville}

\begin{quote}
``$\dotsc$Oriental onlookers are dubbed with pungent comments such as `He's
roasting King Kong' $\ldots$''
\begin{flushright}
--\emph{The New York Times}, 6/27/1963
\end{flushright}
\end{quote}

We provide intuition for the fact that $\MIN$ is a $\Sigma$-set while $\MIN^*$,
$\MIN^\m$, and $\MIN^\T$ are all $\Pi$-sets.

\section{Intuition}

\begin{defn}
Let $\equiv_\alpha$ be an equivalence relation, and let $A, B \subseteq
\omega.$  Define the relation
\[ A \equiv_{\thick{\equiv_\alpha}} B \iff (\forall n) \left[A^{[n]} \equiv_\alpha B^{[n]}\right]. \]
Similarly,
\[ A \leq_{\thick{\equiv_\alpha}} B \iff (\forall n) \left[A^{[n]} \leq_\alpha B^{[n]}\right]. \]
\end{defn}
Note that for any equivalence relation $\equiv_\alpha$ on $\omega$,
$\equiv_{\thick{\equiv_\alpha}}$ is also an equivalence relation.  Informally,
$\thick\equiv_\alpha$ requires agreement on every row.

\begin{defn} \label{thickMINdef}
Let $\equiv_\alpha$ be an equivalence relation.  Then
\[\thick\MIN^{\equiv_\alpha} := \MIN^{\equiv_\alpha} \intersect \MIN^{\thick{\equiv_\alpha}}\]
\end{defn}
This definition is justified by the fact that $\thick{\equiv_\alpha}$ is
intuitively a stronger notion than $\equiv_\alpha$.  Indeed, for any $A, B
\subseteq \omega$, define $m$-equivalent sets $X \equiv_\m A$ and $Y \equiv_\m
B$ by
\begin{align*}
X^{[0]} &:= A & Y^{[0]} &:= B\\
X^{[n+1]} &:= \emptyset & Y^{[n+1]} &:= \emptyset.
\end{align*}
Then
\[X \equiv_{\thick{\equiv_\alpha}} Y \implies A \equiv_\alpha B.\]
Note also that $\MIN^* = \MIN^* \intersect \MIN^{\thick*} \supseteq
\MIN^{\thick*}$, eliminating the need for the $\equiv_\m$-equivalent sets $X$
and $Y$ at the $*$-level:

\begin{prop} \label{thickprop}
For any equivalence relations $\equiv_\alpha$ and $\equiv_\beta$,
\begin{thmenum}
\item $\MIN^{\equiv_\alpha} \supseteq \thick\MIN^{\equiv_\alpha}$, and
\label{thickpropmininthick}

\item If $(\forall A, B \subseteq \omega) \left[ A \equiv_\alpha B \implies A
\equiv_\beta B \right]$ then $\thick\MIN^{\equiv_\alpha} \supseteq
\thick\MIN^{\equiv_\beta}$. \label{thickpropminandthick}

\end{thmenum}
\end{prop}

Proposition~\ref{thickprop}\nlb(\ref{thickpropmininthick}) formally insinuates
that ``$\thick{\equiv_\alpha}$'' is a stronger relation than $\equiv_\alpha$,
simply because it is always possible to move all the information encoded in a
set into a single row.

Theorem~\ref{thicksubset}\nlb(\ref{thickandjump}) gives basic set-theoretic
properties of the modified thick operator from Definition~\ref{thickMINdef}. If
one wishes to deal strictly with $\MIN$-sets, however, a brief inspection of
Theorem~\ref{thicksubset}\nlb(\ref{thickandjump}) also reveals that
$\MIN^{\thick\T^{(n)}} \supseteq \MIN^{\T^{(n+1)}}$.  On the other hand,
$\MIN^{\T^{(n)}} \supseteq \MIN^{\thick{\T^{(n)}}}$ is not true in general
because two sets which are Turing equivalent need not contain their respective
``information'' in identical rows.

We illustrate this last point with an example.  Define sets $A$ and $B$ by
\begin{align*}
A^{[0]} &:= \K & B^{[0]} &:= \emptyset \\
A^{[i+1]} &:= \emptyset & B^{[1]} &:= \K \\
&& B^{[i+2]} &:= \emptyset.
\end{align*}
Then clearly $A \equiv_{\T^{(n)}} B$, but $A \mathrel{\vert_{\thick{\T^{(n)}}}}
B$.

\begin{thm} \label{thicksubset}
Let $n \geq 0$.  Then:
\begin{thmenum}
\item $\MIN^* \supseteq \MIN^{\thick *}$. \label{MIN*thickMIN*}

\item $\thick\fR = \fMIN$. \label{thickMIN}

\item $\MIN^{\thick=} = \MIN$. \label{thickMIN2}

\item $\MIN^{\T^{(n)}} \supseteq \thick\MIN^{\T^{(n)}} \supseteq \MIN^{\T^{(n+1)}}$. \label{thickandjump}
\end{thmenum}
\end{thm}

\begin{proof}
\begin{proof}[(\ref{MIN*thickMIN*})]
$A =^* B$ implies $A \equiv_{\thick*} B$.
\end{proof}

\begin{proof}[(\ref{thickMIN})]
We interpret $\fR$ as a single row in $\fMIN$.  Define
\begin{align*}
\thick\fR :&= \left\{e :  (\forall j < e)\: (\forall n) \left[\phe_e (n) \neq \phe_j(n)\right]  \right\} \\
&= \fMIN. \qedhere
\end{align*}
\end{proof}

\begin{proof}[(\ref{thickMIN2})]
Two sets are equal if they agree on all rows.
\end{proof}

\begin{proof}[(\ref{thickandjump})]  Assume $A^{(n)} \leq_{\thick\T} B^{(n)}$.
Then $A^{(n+1)} \leq_{\thick\m} B^{(n+1)}$ by the Jump Theorem \citep{Soa87}.
For every $i$, let $f_i$ be the computable function that witnesses
$\paren{A^{(n+1)}}^{[i]} \leq_\m \paren{B^{(n+1)}}^{[i]}$.  We create a
function $h$ which captures the values for all the $f_i$'s.

Define a computable function $h$ by
\[ h (\pair{i,x},s) :=
\begin{cases}
f_i(x) & \text{if $\left[x \in \omega^{[i]} \andd i \leq s \right]$}, \\
0 & \text{otherwise.}
\end{cases}\]
Let \[\hat{h}(\pair{i,x}) := \lim_s h(\pair{i,x},s).\]  By the Limit Lemma,
$\hat{h} \leq_\T \zero'$. Since $\hat{h}(\pair{i,x}) = f_i (x)$ for all $i$, we
have that \[A^{(n+1)} \leq_\T B^{(n+1)} \join \zero' \equiv_\T B^{(n+1)}.\]
Hence
\[A^{(n)} \equiv_{\thick\T} B^{(n)} \implies A^{(n+1)} \equiv_\T B^{(n+1)},\]
which gives the second inclusion (by Proposition~\ref{MINinclusionprop}).  The
first inclusion follows immediately from
Proposition~\ref{thickprop}\nlb(\ref{thickpropmininthick}).
\end{proof}
\end{proof}

Theorem~\ref{thicksubset}\nlb(\ref{thickMIN}) gives some intuition why $\MIN
\in \Sigma_2$.  $\MIN$ trivially equals $\MIN^{\thick=}$, so unlike the other
``natural'' $\MIN$-sets, $\MIN$ actually doubles as a ``thick'' set.  In the
case of functions, we note that $\fMIN$ doubles as a nontrivial thick set.

\begin{rem}
Theorem~\ref{thicksubset}\nlb(\ref{thickandjump}) unambiguously shows that the
jump operator defeats nondeterminacy on c.e$.$ sets when we interpret rows as
nondeterministic enumerations of c.e$.$ sets.
\end{rem}

\section{Arithmetics} \label{thickarithmetics}
We now redeem the thick operator by showing that, just like the jump operator,
thick ``kicks'' complete sets up one level in the arithmetic hierarchy.  Still,
jump and thick are juxtaposed here as antithetical: thick ``kicks'' relations
into $\Pi_n$ whenever the jump operator ``kicks'' them into $\Sigma_n$.

\begin{prop}\label{thickarithmetic}
 Let $n \geq 0$.  Then
\begin{thmenum}

\item $\MIN^{\thick *} \in \Sigma_4$. \label{MINthick*insigma4}

\item $\MIN^{\thick\m} \in \Sigma_4$. \label{MINthickminsigma4}

\item $\MIN^{\thick\T^{(n)}} \in \Sigma_{n+5}$. \label{MINthickTinsigma5}
\end{thmenum}
\end{prop}

\begin{proof}
\begin{proof}[(\ref{MINthick*insigma4})]
$\{ \pair{j,e} : W_j =^* W_e \} \in \Sigma_3$, so $\{ \pair{j,e} : W_j
\equiv_{\thick *} W_e \} \in \Pi_4$. \qedhere
\end{proof}

\begin{proof}[(\ref{MINthickminsigma4})]
$\{ \pair{j,e} : W_j \equiv_{\thick\m} W_e \} \in \Pi_4$.
\end{proof}

\begin{proof}[(\ref{MINthickTinsigma5})]
$\{ \pair{j,e} : W_j \equiv_{\thick{\T^{(n)}}} W_e \} \in \Pi_{n+5}$.
\end{proof}
\end{proof}

Theorem~\ref{thicklowerbound} gives the lower bounds.

\begin{thm} \label{thicklowerbound}
Let $n \geq 0$.  Then
\begin{thmenum}
\item $\MIN^{\thick *} \not\in \Pi_4$. \label{thick*lowerbound}

\item $\MIN^{\thick\m} \not\in \Pi_4$. \label{thickmlowerbound}

\item $\MIN^{\thick{\T^{(n)}}} \not\in \Pi_{n+5}$. \label{thickTlowerbound}
\end{thmenum}
\end{thm}

\begin{proof}

\begin{proof}[(\ref{thick*lowerbound})]
Let $A \in \Pi_4$.  Then there exists a relation $R \in \Sigma_3$ such that
\[x \in A \iff (\forall y)\, R(x,y).\]
Since $\COF$ is $\Sigma_3$-complete \citep{Soa87}, there exists a computable
function $g$ such that $R(x,y)$ iff $W_{g(x,y)}$ is cofinite.  Therefore
\[x \in A \iff (\forall y) \left[W_{g(x,y)} =^* \omega \right].\]
Define a computable function $f$ by \[\phe_{f(x)}^{[y]} := \phe_{g(x,y)}.\]
Then
\begin{align*}
W_{f(x)} \equiv_{\thick*} \omega &\iff (\forall y) \left[W_{g(x,y)} =^* \omega
\right] \\
&\iff x \in A,
\end{align*}
which makes \[\thick\COF := \{ e : W_e \equiv_{\thick*} \omega \}\]
$\Pi_4$-complete.

Now, as in Theorem~\ref{lowerbound}, suppose towards a contradiction that
$\MIN^{\thick*} \in \Pi_4$, and let $a$ be the $\equiv_{\thick*}$-minimal index
for $\omega$.  Then
\begin{align*}
\thick\COF &= \{e : W_e \equiv_{\thick*} \omega \} \\
&= \{a\} \union \left\{e : (\forall j < e) \left[ j \in \MIN^{\thick*} - \{a\}
\implies W_j \not \equiv_{\thick*} W_e \right] \right\}.
\end{align*}
Now $\thick\COF \in \Sigma_4$, since $W_j \equiv_{\thick*} W_e$ can be decided
in $\Pi_4$, and because \[\MIN^{\thick*} - \{a\} \in \Pi_4\] by assumption.
This contradicts the fact that $\thick\COF$ is $\Pi_4$-complete.
\end{proof}

\begin{proof}[(\ref{thickmlowerbound})]
Let $\K^\omega$ be the c.e$.$ set in which each row is the halting set; for all
$k$, \[ (\K^\omega)^{[k]} := \K,\] and recall that
\[ \mCOMP = \{ e : W_e \equiv_\m K \} \]
is $\Sigma_3$-complete (Theorem~\ref{lowerbound}\nlb(\ref{MINmlowerbound})). By
an argument analogous to part~(\ref{thick*lowerbound}), we have that
\[\thick\mCOMP := \left\{ e : W_e \equiv_{\thick\m} \K^\omega \right\}\]
is $\Pi_4$-complete.

Suppose $\MIN^{\thick\m} \in \Pi_4$.  Following the same line of reasoning as
before, and noting that $W_j \equiv_{\thick\m } W_e$ can be decided in $\Pi_4$,
we obtain a contradiction.
\end{proof}

\begin{proof}[(\ref{thickTlowerbound})]
We use the same reasoning a third time.  Let $\K^{(n)^\omega}$ be the c.e$.$
set given by \[ \left(\K^{(n)^\omega}\right)^{[i]} := K^{(n)},\] for all $i$,
and recall that
\[ \HIGH^n = \left\{ e : (W_e)^{(n)} \equiv_\T K^{(n)} \right\} \]
is $\Sigma_{n+4}$-complete (see Theorem~\ref{lowerbound}\nlb(\ref{MINTlowerbound})).  By an argument analogous to
part~(\ref{thick*lowerbound}), we have that
\[\thick\HIGH^n := \left\{ e : (W_e)^{(n)} \equiv_{\thick{\T}} \K^{(n)^\omega} \right\}\]
is $\Pi_{n+5}$-complete.

Suppose $\MIN^{\thick\m} \in \Pi_{n+5}$.  Following the same line of reasoning
as before, and noting that $W_j \equiv_{\thick\T^{(n)}} W_e$ can be decided in
$\Pi_{n+5}$, we obtain a contradiction.
\end{proof}
\end{proof}

\section{Immunity}

Thickness contributes nothing to immunity, as evidenced by
Theorem~\ref{thickimmune}.

\begin{lemma}[semi-fixed points]  There exists a computable function $\nu$ such
that \label{semi-fixed}
\begin{thmenum}
\item $f \leq_\T \zero' \enskip \implies \enskip (\exists e) \left[W_{\nu(e)}
\equiv_{\thick*} W_{f(e)} \right]$, \label{semi-fixed*}

\item $f \leq_\T \zero'' \enskip \implies \enskip (\exists e) \left[W_{\nu(e)} \equiv_{\thick\m}
W_{f(e)} \right]$, \label{semi-fixedm}

\item $f \leq_\T \zero^{(n+2)} \enskip \implies \enskip (\exists e) \left[W_{\nu(e)} \equiv_{\thick{\T^{(n)}}}
W_{f(e)} \right]$. \label{semi-fixedT}
\end{thmenum}
\end{lemma}

\begin{proof}

\begin{proof}[(\ref{semi-fixed*})]
The proof for part~(\ref{semi-fixedm}) will work.
\end{proof}

\begin{proof}[(\ref{semi-fixedm})]
Using the $s$-$m$-$n$ Theorem, define a computable function $\nu$ by
\[\phe_{\nu(x)}\paren{\pair{z,n}} :=
\begin{cases}
\phe_{\phe_{x}(n)}(z) &\text{if $\phe_x (n) \converge$} \\
\diverge &\text{otherwise}.
\end{cases}
\] so that for any $x \in \TOT$, \[W_{\nu(x)}^{[n]} = W_{\phe_x (n)}. \]

Let $f \leq_\T \zero''$, and define,  again using the $s$-$m$-$n$ Theorem, a
computable sequence of $\zero''$-computable functions $\{f_n\}$ by
\[\phe_{f_n(x)} (z) := \phe_{f (x)} \paren{\pair{z,n}} \] so that \[ W_{f_n(x)} = W^{[n]}_{f (x)}.\]

By the Generalized Fixed Point Theorem~\ref{JLSSmain}, we can uniformly find a
computable sequence $\{e_n\}$ such that for all $n$, \[W_{e_n} \equiv_\m
W_{f_n(e)}.\]   Let $e$ be an index so that \[ \phe_e (n) := e_n.\]   Then for
all $n$,
\[W^{[n]}_{\nu(e)} =  W_{\phe_e(n)}  =  W_{e_n} \equiv_\m W_{f_n(e)} = W^{[n]}_{f (e)}. \]
This means that
\begin{equation}\label{thickmfixedpoint}
(\forall f \leq_\T \zero'')\: (\exists e) \left[ W_{\nu(e)} \equiv_{\thick\m}
W_{f(e)} \right],
\end{equation}
which is what we intended to show.
\end{proof}

\begin{proof}[(\ref{semi-fixedT})]
Our proof of (\ref{semi-fixedm}) used no specific properties of $\equiv_\m$
except that this relation satisfies Generalized Fixed Point
Theorem~\ref{mintnimmune}. Since $\equiv_{\T^{(n)}}$ satisfies analogous fixed
point properties, the same argument will work.
\end{proof}
\end{proof}

Comparing Theorem~\ref{thickimmune} with the results from
Chapter~\ref{immunitychapter}, we note that the thick operator does not at all
affect the immunity of our main equivalence relations:

\begin{thm} Let $n \geq 0$.  Then \label{thickimmune}
\begin{thmenum}
\item $\MIN^{\thick *}$ is $\Sigma_2$-immune but not $\Sigma_3$-immune.
\label{thickMIN*immune}
\item $\MIN^{\thick\m}$ is $\Sigma_3$-immune (and not $\Sigma_4$-immune). \label{thickMINmimmune}
\item $\MIN^{\thick\T^{(n)}}$ is $\Sigma_{n+3}$-immune but not $\Sigma_{n+4}$-immune. \label{thickMINTimmune}
\end{thmenum}
\end{thm}

\begin{proof}

\begin{proof}[(\ref{thickMIN*immune})]
$\MIN^{\thick*}$ is $\Sigma_2$-immune follows immediately from
Theorem~\ref{thicksubset}(\ref{MIN*thickMIN*}) and
Theorem~\ref{RMINMIN*immune}\nlb(\ref{MIN*immune}). We show $\MIN^{\thick*}$ is
not $\Sigma_3$ immune by modifying the proof of
Theorem~\ref{pi3separation}\nlb(\ref{pi3separation*}).  All that is needed is
to change the definition of $A_k$ so that it only applies to the first row of
each c.e$.$ set:
\[A_k := \left\{e : W_e^{[0]} \subseteq^* P_k \right\} \intersect \INF.\]
The rest of the proof is the same.
\end{proof}

\begin{proof}[(\ref{thickMINmimmune})]
We tweak the proof of $\Pi_3$-Separation
Theorem~\ref{pi3separation}\nlb(\ref{pi3separationm}).  Let $A$ be an infinite,
$\Sigma_3$ set. Suppose $A \subseteq \MIN^{\thick\m}$. Since $A$ is infinite
and c.e$.$ in $\zero''$, we can define a $\zero''$-computable function $f$ by
\[f(x) = \pi_1 \Bigl( \left( \mu \pair{i,t} \right)\: [i > \nu(x) \quad \& \quad i \in A_t] \Bigr), \]
where $\{A_t\}$ is a $\zero''$-enumeration of $A$ and $\nu$ is the computable
function from Lemma~\ref{semi-fixed}.

Now for all $x$, $f(x) > \nu(x)$ and $f(x) \in \MIN^{\thick\m}$.  Therefore
\[ (\forall x)\ [W_{\nu(x)} \not\equiv_{\thick\m} W_{f(x)}],\]
contradicting Lemma~\ref{semi-fixed}.
\end{proof}

\begin{proof}[(\ref{thickMINTimmune})]
A rehash of ideas from parts (\ref{thickMIN*immune}) and
(\ref{thickMINmimmune}).  To show that $\MIN^{\thick\T^{(n)}}$ is not
$\Sigma_4$-immune, we use the proof from Theorem~\ref{minTnotimmune} in place
of Theorem~\ref{pi3separation}(\ref{pi3separation*}): just redefine
\[B_k := \biggl[ \Bigl\{ x : \left((W_x)^{(n)}\right)^{[0]} \leq_\T (A_k)^{(n)} \Bigr\}\ \intersect\ \compliment{\LOW^n}
\biggr],\] and then the proof from Theorem~\ref{minTnotimmune} works.

The argument in (\ref{thickMINmimmune}) shows that $\MIN^{\T^{(n)}}$ is
$\Sigma_3$ immune.
\end{proof}
\end{proof}

%% file: kolmogorov.tex
\chapter{A Kolmogorov numbering} \label{a kolmogorov numbering}
For certain G\"{o}del numberings, we can exactly determine the truth-table
degree of $\MIN$, $\MIN^*$, and $\MIN^\m$ as well as the Turing degrees of
$\MIN^{\T^{(n)}}$, and $\MIN^{\thick*}$. The main result of this chapter,
Theorem~\ref{thickpheTuringcomplete}, provides a Kolmogorov numbering in which
$\MIN$-sets exactly characterize the Turing degrees $\0$, $\0'$, $\0''$,
$\dotsc$.

\section{Numberings I \& II}

\subsection{Numbering I}

Theorem~\ref{thickorderlemma}, restricted to $\fMIN_\psi$ and $\fMIN^*_\psi$,
was first proved by Schaefer \citep{Sch98}.  He also mentions a G\"{o}del
ordering satisfying (\ref{fRpheis0'thick}) (see
Theorem~\ref{numberingsRANDfR}).  The majority of constructions here are
inspired by ~\citep[Theorem 2.17]{Sch98}.

\begin{thm} \label{thickorderlemma}
There exists a Kolmogorov numbering $\psi$ simultaneously satisfying:
\begin{thmenum}
\item $\fR_\psi \geq_\Tt \zero'$, \label{fRpheis0'thick}

\item $\MIN, \fMIN_\psi \geq_\Tt \zero''$, \label{MINpheis0''}

\item $\MIN^*, \fMIN^*_\psi \geq_\Tt \zero'''$, \label{MIN*pheis0'''}

\item $\MIN_\psi^{\thick *} \geq_\Tt \zero'''$,

\item $\MIN_\psi^{\thick \m} \geq_\Tt \zero'''$, and \label{thickphemis0''''}

\item $\MIN_\psi^{\thick \T^{(n)}} \geq_\Tt \zero^{(n+4)}$.
\label{thickpheTis0'''''}
\end{thmenum}
\end{thm}

\begin{proof}
We first construct a G\"{o}del numbering $\psi$ satisfying
(\ref{thickpheTis0'''''}).  We later argue that our construction can be
modified to produce a Kolmogorov numbering satisfying all six parts of the
lemma.

Let $\phe$ be any G\"{o}del numbering, and let $n \geq 0$.  We define the
numbering $\psi$ as follows.  Define an increasing, computable function $f$ by
\begin{align*}
f(0) &:= 0, \\
f(k+1) &:= 4[f(k) + 1] + 1,
\end{align*}
Let $i \geq 0$.  If $i = f(k)$ for some $k$, then we define $\psi_i := \phe_k$.
This makes $\psi$ an effective ordering.  Otherwise, for some $k$, $f(k) < i <
f(k+1)$. In this case we define
\begin{equation} \label{thickpheequation1}
\psi_i (\pair{x,y}) :=
\begin{cases}
1 &\text{if} \enskip \left[y - f(k) \text{ is odd} \andd y=i \andd \phe_x(x)\converge \right], \\
1 &\text{if} \enskip \left[y - f(k) \text{ is even} \andd y = i - 1 \andd \phe_k (x) \converge\right], \\
\diverge &\text{otherwise.}
\end{cases}
\end{equation}
The functions $\psi_{f(k)+1}, \psi_{f(k)+3} \dotsc, \psi_{4[f(k)+1] - 1}$ code
the halting set into distinct rows, and the remaining functions between $f(k)$
and $f(k+1)$ are used for comparisons.

It remains now only to show that
\[\HIGH^n_\phe \leq_\Tt \MIN_\psi^{\thick{\T^{(n)}}},\]
because $\HIGH^n_\phe$ is $\Sigma_{n+4}$ complete
(Theorem~\ref{lowerbound}\nlb(\ref{MINTlowerbound})).  Here we use the
subscript ``$\phe$'' to emphasize that we are considering $\HIGH^n$ with
respect to the numbering $\phe$.

We claim that
\begin{equation} \label{thickpheequation2}
e \in \HIGH^n_\phe \bigiff \left[\MIN_\psi^{\thick{\T^{(n)}}} \intersect \{f(k)
+ 2, f(k) + 4, \ldots, 4f(k) + 4\}\right] = \emptyset,
\end{equation}
where $k$ is such that $f(k) \leq e < f(k+1)$.  The claim follows by inspecting
pairs of functions $\{\psi_i, \psi_{i+1}\}$.  Indeed, assume $e \in
\HIGH^n_\phe$. Then for all $y$, including $y = f(k) + 1$,
$$\left(\dom \psi_{f(k)+1} \right)^{[y]} \equiv_{\T^{(n)}} \left(\dom \psi_{f(k)+2} \right)^{[y]}.$$
Therefore
$$\dom \psi_{f(k) + 1} \equiv_{\thick\T^{(n)}} \dom \psi_{f(k) + 2},$$
which means that
$$f(k) + 2 \not\in \MIN_{\psi}^{\thick{\T^{(n)}}}.$$
Similarly,
\[f(k)+4, f(k)+6 \dotsc, 4f(k) + 4 \not\in \MIN_{\psi}^{\thick{\T^{(n)}}},\]
which proves the first direction.

Conversely, assume that $e \not\in \HIGH_\phe^n$.  Then for all $i \neq j$,
with
\[i, j \in \{f(k) + 1, f(k) + 2, \ldots, 4f(k) + 4\},\]
we have
\[\psi_i \not\equiv_{\thick{\T^{(n)}}} \psi_j.\]
This means that for $k \geq 1$, $$[4f(k) + 4] - f(k) = 3f(k) + 4$$ distinct
$\equiv_{\thick\T^{(n)}}$-equivalence classes are represented in
\begin{equation}
\{\psi_{f(k)+1}, \psi_{f(k)+2}, \dotsc \psi_{4f(k) + 4}\}.
\label{thickTpheenumeration}
\end{equation}
It follows that at least
$$[3f(k) + 4] - (f(k) + 1) = 2f(k) + 3$$
of the indices from \eqref{thickTpheenumeration} are
$\equiv_{\thick\T^{(n)}}$-minimal, since only those classes also represented in
$\{\psi_1, \dotsc, \psi_{f(k)}\}$ could be
$\equiv_{\thick\T^{(n)}}$-nonminimal.  Thus, any subset from
\[\{f(k)+1, f(k)+2, \dotsc, 4f(k)+4\}\]
with cardinality at least $f(k) + 2$ must contain a
$\equiv_{\thick\T^{(n)}}$-minimal index.  In particular,
\[\left[\MIN_\psi^{\thick{\T^{(n)}}} \intersect \{f(k) + 2, f(k) + 4,
\ldots, 4f(k) + 4\}\right] \neq \emptyset.\] Hence we conclude that
$$\MIN_\psi^{\thick \T^{(n)}} \geq_\Tt \zero^{(n+4)}.$$

 We now describe separate orderings satisfying (\ref{fRpheis0'thick}) --
(\ref{thickphemis0''''}), and then we show that all six numberings can be
combined together into a single  G\"{o}del numbering.  Finally, we argue that
this G\"{o}del numbering can be made into an Kolmogorov numbering by
ambiguously appealing to \citep[Theorem 2.17]{Sch98}.

The remaining, individual numberings are either identical or similar to the
numbering $\psi$ which we just constructed. For instance, the same $\psi$
satisfies
\[\MIN_\psi^{\thick\m} \geq_\Tt \zero'''.\]
In fact, we need only change $\HIGH^n_\phe$ to $\mCOMP_\phe$ in the
verification \eqref{thickpheequation2}, and then the same proof works.  For
$\equiv_{\thick*}$, $=^*$, and $=$, we use a different numbering, say $\nu$,
which is exactly like $\psi$ except the condition ``$\phe_x(x)\converge$'' is
omitted from \eqref{thickpheequation1}.  To verify this numbering works, we
swap either $\COF_\phe$ or $\TOT_\phe$ for $\HIGH^n_\phe$ in
\eqref{thickpheequation2}. For $\fR$, we substitute \eqref{thickpheequation1}
with
\[\xi_i (x) :=
\begin{cases}
\pair{i,1} &\text{if \enskip $i$ is odd,} \\
\pair{i-1,1} &\text{if \enskip [$i$ is even \enskip \& \enskip $\phe_k(k)\converge$]}, \\
\diverge &\text{otherwise.}
\end{cases}\]
In the verification for $\fR$, we replace $\HIGH^n_\phe$ in
\eqref{thickpheequation2} with the halting set complement,
$\compliment{K_\phe}$.

We now merge the numberings $\psi$, $\nu$, and $\xi$ into a single G\"{o}del
numbering $\rho$ satisfying (\ref{fRpheis0'thick}) --
(\ref{thickpheTis0'''''}). All we do is change the p.c$.$ functions filling the
coding ``gap'' between $f(k)$ and $f(k+1)$, so that $\psi$ fills the first gap,
$\nu$ fills the second gap, $\xi$ fills the third gap, $\psi$ again fills the
fourth, etc.  Furthermore, we must repeat each $\phe_k$ function three times,
so that each of numbering strategies may ask questions to it.  For this reason,
we let $\phe$ be a Kolmogorov numbering such that $\phe_k = \phe_{k+1} =
\phe_{k+2}$ whenever $k \equiv 0\ (\mathrm{mod}\ 3)$.  We could settle for a
G\"{o}del numbering for the moment, but we'll need $\phe$ to be a Kolmogorov
numbering anyway after the next paragraph.

We define
\[\rho_i := \phe_k \quad \text{when $i = f(k)$ for some $k$}.\]
Otherwise, $f(k) < i < f(k+1)$ for some $k$.  If $k \equiv 0\ (\mathrm{mod}\
3)$ then we use the $\psi$ strategy for $i$, if $k \equiv 1\ (\mathrm{mod}\ 3)$
we use the $\nu$ strategy for $i$, and if $k \equiv 2\ (\mathrm{mod}\ 3)$ we
use the $\xi$ strategy for $k$. So, for example, if $i = 3 \cdot 4567 + 1$,
then
\[\rho_i (\pair{x,y}) :=
\begin{cases}
1 &\text{if} \enskip \left[y - f(k) \text{ is odd} \andd y=i \right], \\
1 &\text{if} \enskip \left[y - f(k) \text{ is even} \andd y = i - 1 \andd \phe_k (x) \converge\right], \\
\diverge &\text{otherwise.}
\end{cases}\]
We can now make truth-table queries to the appropriate spectral sets, just as
before.

Finally, we transform $\rho$ into a Kolmogorov numbering.  The idea is to
enumerate a large number of $\phe_k$'s between each coding ``gap'' instead of
just the one $k$ from $f(k)$. In the $s^{\text{th}}$ gap, we code a crib for
$\phe_s$ in the same manner as we did with $\rho$.  More formally we define, by
induction,
\begin{align} \label{thickpheequation3}
g(0) &:= 0, \\
h(0) &:= 0, \\
g(k+1) &:= g(k) + h(k) + 2[g(k) + 1], \\
h(k+1) &:= 2(h(k) + 2[g(k) + 1]).
\end{align}
Our new numbering is split into blocks $h(k) \leq i < h(k+1)$ rather than $f(k)
\leq i < f(k+1)$ as before.  For $i$ with $$h(k) \leq i < h(k) + 2[g(k) + 1],$$
we apply the familiar coding scheme from $\rho$ (on $\phe_k$), and for $i$ with
$$h(k) + 2[g(k) + 1] \leq i < h(k+1),$$ we simply enumerate $\phe_{g(k)}$ up to $\phe_{g(k+1)-1}$.  This
construction is a Kolmogorov numbering by \cite[Theorem 2.17]{Sch98}, where
this same induction appears.
\end{proof}

\subsection{Numbering II}

\begin{thm}There exists a Kolmogorov numbering $\psi$ such that for all $n \geq 0$:
\begin{thmenum} \label{psimandTttcomplete}
\item $\MIN^\m_\psi \geq_\Tt \zero'''$. \label{mpsittcomplete}

\item $\MIN^{\T^{(n)}}_\psi \geq_\Tt \zero^{(n+3)}$. \label{Tpsittcomplete}
\end{thmenum}
\end{thm}

\begin{proof}
As in Theorem~\ref{thickorderlemma}, we shall first construct a G\"{o}del
numbering $\psi$ satisfying (\ref{mpsittcomplete}) and (\ref{Tpsittcomplete}),
and we later argue that the construction can be modified so as to achieve a
single Kolmogorov numbering.

Let $\phe$ be an arbitrary G\"{o}del numbering, and assume $\pair{\cdot,\cdot}$
is a bijective pairing function satisfying $\pair{0,0} = 0$.  Let $a$ be the
computable function from Lemma~\ref{forcerecursiveorindependent}, defined in
terms of this ordering. Define a computable function $f$ by
\begin{align*}
f(0) &:= 0, \\
f(k+1) &:= 2f(k) + 3.
\end{align*}  The numbering $\psi$ is defined as follows.
Let $C$ be an arbitrary computable set, and let $\psi_0$ be such that \[\dom
\psi_0 := C.\] Let $i \geq 1$.  If $i = f(\pair{k,n})$ for some pair
$\pair{k,n}$, then $\psi_i := \phe_\pair{k,n}$. Otherwise, $f(\pair{k,n}) < i <
f(\pair{k,n}+1)$ for some $\pair{k,n}$.  In this case,
\[\psi_i := \phe_{a_\pair{k,n}(i)}.\]
Let $\LOW^n_\phe$ and $\LOW^n_\psi$ denote the $\LOW^n$ indices in terms of
$\phe$-indices and $\psi$-indices, respectively.

We claim, for $\pair{k,n} > 0$,
\begin{equation*}
\MIN_\psi^{\T^{(n)}} \intersect \{f(\pair{k,n})+1, f(\pair{k,n})+2, \dotsc,
2f(\pair{k,n})+2 \} \neq \emptyset \bigiff k \in \compliment{\LOW_\phe^n}.
\end{equation*}
Indeed, if $k \in \LOW^n_\phe$, then $a_\pair{k,n}(i) \in \LOW^n_\phe$ for all
$i$, hence
\[\{f(\pair{k,n})+1,  \dotsc, 2f(\pair{k,n})+2\} \subseteq \LOW^n_\psi,\]
and so
\[\MIN_\psi^{\T^{(n)}} \intersect \{f(\pair{k,n})+1,  \dotsc, 2f(\pair{k,n})+2 \} = \emptyset.\]

Conversely, if $k \in \compliment{\LOW_\phe^n}$, then by definition of $a$,
each of the $\psi$-indices
\begin{equation}
f(\pair{k,n})+1,  \dotsc, 2f(\pair{k,n})+2 \label{psimandTttcompleteeq}
\end{equation}
represents a distinct $\T^{(n)}$-degree.  At most $f(\pair{k,n}) + 1$ degrees
are represented with smaller indices, so at least one of the $f(\pair{k,n}) +
2$ degrees in \eqref{psimandTttcompleteeq} must be minimal. That is,
\[\MIN_\psi^{\T^{(n)}} \intersect \{f(\pair{k,n})+1,  \dotsc, 2f(\pair{k,n})+2 \} \neq \emptyset.\]
Since $\LOW^n$ is $\Sigma_{n+3}$-complete, this proves that $\psi$ satisfies
(\ref{Tpsittcomplete}).

Similarly, for $k > 0$,
\begin{equation*}
\MIN_\psi^\m \intersect \{f(\pair{k,0})+1, \dotsc, 2f(\pair{k,0})+2 \} \neq
\emptyset \bigiff k \in \compliment{\LOW_\phe^0},
\end{equation*}
which shows that $\psi$ satisfies (\ref{mpsittcomplete}).  One can now
transform $\phe$ into a Kolmogorov numbering by following the familiar
procedure from Theorem~\ref{thickorderlemma}, starting from
\eqref{thickpheequation3}.
\end{proof}

\section{Truth-table apogee}
We present a Kolmogorov numbering for which $\MIN$-sets achieve maximal
truth-table and Turing degrees.

\begin{lemma} \label{thickdownonelevel}
Let $n \geq 0$.
\begin{thmenum}
\item $\MIN^{\thick *} \join \zero''' \equiv_\bT \zero''''$,

\item $\MIN^{\thick \m} \join \zero''' \equiv_\bT \zero''''$,

\item $\MIN^{\thick \T^{(n)}} \join \zero^{(n+4)} \equiv_\bT \zero^{(n+5)}$.
\end{thmenum}
\end{lemma}

\begin{proof}
The same proof from Lemma~\ref{downonelevel}\nlb(\ref{MINjoin0'is0''}) works
here when we substitute the fact that either $\thick\COF$ is $\Pi_4$-complete,
$\thick\mCOMP$ is $\Pi_4$-complete, or $\thick\HIGH^n$ is $\Pi_{n+5}$-complete
for the fact that $\TOT$ is $\Pi_2$-complete.  Definitions for $\thick\COF$ and
$\thick\mCOMP$ appear in the proof of Theorem~\ref{thicklowerbound}.
\end{proof}

Combining the orderings from Lemma~\ref{thickorderlemma} and
Lemma~\ref{psimandTttcomplete} (using techniques from these lemmas), we obtain:

\begin{thm} \label{ttapogeethm}
There exists a Kolmogorov numbering $\psi$ satisfying
\begin{thmenum}
\item $\fR_\psi \geq_\Tt \zero'$,

\item $\MIN_\psi, \fMIN_\psi \geq_\Tt \zero''$,

\item $\MIN^*_\psi, \fMIN^*_\psi \geq_\Tt \zero'''$,

\item $\MIN^\m_\psi \geq_\Tt \zero'''$,

\item $\MIN^{\T^{(n)}}_\psi \geq_\Tt \zero^{(n+3)}$,

\item $\MIN_\psi^{\thick *} \geq_\Tt \zero'''$,

\item $\MIN_\psi^{\thick \m} \geq_\Tt \zero'''$,

\item $\MIN_\psi^{\thick \T^{(n)}} \geq_\Tt \zero^{(n+4)}$.
\end{thmenum}
\end{thm}

Using the numbering from Theorem~\ref{ttapogeethm}, together with
Lemma~\ref{thickdownonelevel} and Lemma~\ref{downonelevel}, we can conclude the
following.

\begin{cor}There exists a Kolmogorov numbering $\psi$ simultaneously satisfying:
\begin{thmenum} \label{thickpheTuringcomplete}
\item $\fR_\psi \equiv_\Tt \zero'$,

\item $\MIN_\psi \equiv_\Tt \fMIN_\psi \equiv_\Tt \zero''$,
\label{thickpheTuringcompleteMIN}

\item $\MIN^*_\psi \equiv_\Tt \fMIN^*_\psi \equiv_\Tt \zero'''$,

\item $\MIN^\m_\psi \equiv_\Tt \zero'''$,

\item $\MIN^{\T^{(n)}}_\psi \equiv_\T \zero^{(n+4)}$,

\item $\MIN_\psi^{\thick *} \equiv_\T \zero''''$,

\item $\MIN_\psi^{\thick \m} \equiv_\T \zero''''$, and

\item $\MIN_\psi^{\thick \T^{(n)}} \equiv_\T \zero^{(n+5)}$.
\end{thmenum}
\end{cor}

%% file: peak_hierarchy.tex
\chapter[Peak Hierarchy]{Hyperimmunity and the Peak Hierarchy Theorem}
\label{The Peak Hierarchy Theorem}

In Corollary~\ref{containingMINTomega}, we exhibit an infinite sequence of
indices which is common to all spectral sets herein.  We conclude that spectral
sets are not hyperimmune, and we use this fact to build a special ``cutting
set'' in the last section.

\section{A computable sequence of intermediate degrees}
The main goal of this section is to prove Theorem~\ref{intseq}.

\subsection{Main theorem}

\begin{thm}\label{intseq}
There exists a computable sequence $\{x_k\}$ such that for all $n$ and $i$,
\[\left( W_{x_i} \right)^{(n)} \not\leq_\T  \ijoin_{j \neq i}  \left(W_{x_j} \right)^{(n)}.\]
\end{thm}
In particular, $(W_{x_i})^{(n)} \Tincomp (W_{x_j})^{(n)}$ whenever $i \neq j$.
\begin{proof}
The proof uses three lemmata.  Lemma~\ref{antichain} creates a computable
sequence ``upstairs,'' above $\zero'$. Lemma~\ref{sacksjumpseq} brings that
sequence ``downstairs'' using infinite injury from the Sacks Jump Theorem.
Finally, we ``take the elevator to the top'' by way of
Lemma~\ref{2infrecursion}. We now prove Theorem~\ref{intseq}, assuming these
lemmas.

Let $\{a_k\}$ be the sequence of computable functions guaranteed by
Lemma~\ref{sacksjumpseq}.  Then by \eqref{sacksjumpnoteq}, for any $s$ and $i
\neq j$,
\begin{equation}\label{intseqweak}
W_{a_i(s)}^Y \join Y \not\equiv_\T W_{a_j(s)}^Y \join Y.
\end{equation}
By Lemma~\ref{2infrecursion}, let $x$ be the fixed point for the sequence
$\{a_k\}$ satisfying \[\Psi_{a_k(x)}^Y = \Psi_{x_k}^Y,\] so that $\{x_k\}$ is
computable.  Let $i \neq j$.  Taking $Y$ to be the empty set and applying
\eqref{sacksjumpeq} and \eqref{intseqweak} yields
\begin{align*}
\left( W_{x_i} \right)' &\equiv_\T \left( W_{x_i}^{\zero} \join \zero \right)'= \left(
W_{a_i(x)}^\zero \join \zero \right)' \equiv_\T W_{a_i(x)}^{\zero'} \join \zero'
\\[0.5ex]
&\not\equiv_\T W_{a_j(x)}^{\zero'} \join \zero' \equiv_\T \left( W_{a_j(x)}^\zero \join
\zero \right)' = \left(W_{x_j}^\zero \join \zero \right)' \equiv_\T \left( W_{x_j}
\right)'.
\end{align*}
In general, applying \eqref{sacksjumpeq} and \eqref{sacksjumpnoteq} 
from Lemma~\ref{sacksjumpseq} yields
\begin{align*}
\left( W_{x_i} \right)^{(n)} &\equiv_\T \left( W_{x_i}^\zero \join \zero \right)^{(n)} =
\left( W_{a_i(x)}^\zero \join \zero \right)^{(n)} \equiv_\T \left( W_{a_i(x)}^{\zero'} \join \zero' \right)^{(n-1)} \\
&\enskip \vdots \\
&\equiv_\T W_{a_i(x)}^{\zero^{(n)}} \join \zero^{(n)} \not\leq_\T \ijoin_{j \neq i} \left( W_{a_j(x)}^{\zero^{(n)}} \join \zero^{(n)} \right) \\[0.5ex]
&\phantom{\equiv_\T W_{a_i(x)}^{\zero^{(n)}} \join \zero^{(n)}\;}\equiv_\T \ijoin_{j \neq i} \left( W_{a_j(x)}^{\zero^{(n-1)}} \join \zero^{(n-1)} \right)' \\[0.5ex]
&\phantom{\equiv_\T W_{a_i(x)}^{\zero^{(n)}} \join \zero^{(n)}\;}\enskip \vdots \\
&\phantom{\equiv_\T W_{a_i(x)}^{\zero^{(n)}} \join \zero^{(n)}\;}\equiv_\T \ijoin_{j \neq i} \left( W_{a_j(x)}^\zero \join \zero \right)^{(n)} \\[0.5ex]
&\phantom{\equiv_\T W_{a_i(x)}^{\zero^{(n)}} \join \zero^{(n)}\;}= \ijoin_{j \neq i} \left( W_{x_j}^\zero \join \zero \right)^{(n)} \\[0.5ex]
&\phantom{\equiv_\T W_{a_i(x)}^{\zero^{(n)}} \join \zero^{(n)}\;}\equiv_\T
\ijoin_{j \neq i}  \left( W_{x_j} \right)^{(n)}. \qedhere
\end{align*}
\end{proof}

\subsection{Three lemmas}

Lemma~\ref{antichain} was first proved by Kleene and Post in 1954
\citep[Theorem 3.3.1]{KP54}. A more recent, nonrelativized exposition appears
in Odifreddi's book \citep[Proposition V.2.7]{Odi89}.  We isolate this proof in
order to clarify intuition for Lemma~\ref{sacksjumpseq}.

\begin{lemma}\label{antichain}
There exists a computable sequence $\{a_i\}$ such that for any $Y \subseteq \omega$ and
$i \in \omega$,
\[W_{a_i}^Y \join Y \not\leq_\T \ijoin_{j \neq i} \left( W_{a_j}^Y \join Y \right).\]
\end{lemma}
\begin{proof}
We reuse the proof of Lemma~\ref{forcerecursiveorindependent}. Our construction
here is exactly the same as before, but without permitting (we omit ``$c(s)
\leq x$'' from \eqref{recursiveorindependenteq}).  The only other difference is
that we use $A_i$ to denote the relevant set $W_{a_i}^Y \join Y$ rather than
the irrelevant set, $(W_{a_\pair{k,n}(i)})^{(n)}$. The claims from before then
follow verbatim:

\begin{claim}\label{FMmet}
If requirement $R_\pair{e,i}$ acts at some stage $s+1$ and is never later
injured, then requirement $R_\pair{e,i}$ is met and $r \left(\pair{e,i}, t
\right) = s + 1$ for all $t \geq s + 1$.
\end{claim}

\begin{claim}\label{FMmet2}
For every $\pair{e,i}$, requirement $R_\pair{e,i}$ is met, acts at most
finitely often, and $r\left(\pair{e,i}\right) := \lim_s r\left( \pair{e,i}, s
\right)$ exists.
\end{claim}

Recall that $R_\pair{e,i}$ was the requirement
\[R_\pair{e,i} : A_i \neq \Psi_e^{\ijoin_{j \neq i} A_j}.\]
Thus the lemma follows immediately from Claim~\ref{FMmet2}.
\end{proof}

\begin{lemma}\label{sacksjumpseq}
There exists a computable sequence of computable functions $\{a_i\}$ such that for any
computable sequence $\{s_i\}$ and $Y\subseteq \omega$,
\begin{equation}\label{sacksjumpeq}
(W_{a_i(s)}^Y \join Y)' \equiv_\T W_{s_i}^{Y'} \join Y' \quad\text{and}
\end{equation}
\begin{equation}\label{sacksjumpnoteq}
W_{a_i(s)}^Y \join Y \not\leq_\T \ijoin_{j \neq i} \left( W_{a_j(s)}^Y \join Y \right),
\end{equation}
where $s$ is such that $\phe_s(i) = s_i$.  
\end{lemma}

\begin{proof}
We mix the ``true stages'' proof of the Sacks Jump Theorem \citep{Soa87}
together with Lemma~\ref{antichain}, however we omit the ``avoid the cone''
strategy.  The main idea is as follows. In the proof of the Sacks Jump Theorem
(Theorem~\ref{sacksjumptheorem}), one constructs a single set $A$ satisfying
the thickness requirement
\[P_\pair{e,0}: A_0^{[e]} =^* B_0^{[e]},\] where $A_0$ and $B_0$ are analogous
to sets defined below.  In this proof, we will simultaneously construct
infinitely many sets $\{A_i\}$ satisfying requirement $R_\pair{e,i}$ from
Lemma~\ref{antichain} by playing the strategy of Lemma~\ref{antichain} on the
first column of each matrix, $A_0^{[\widetilde{0}]}$, $A_1^{[\widetilde{0}]},
\cdots$.  This strategy won't interfere with the corresponding $P_\pair{e,i}$
requirements because we're changing just one point in each row.

We now prove the lemma.  We construct $\{a_i\}$ uniformly in $s$ and
independent of $Y$. For convenience, let $S_i$ denote the set $W_{s_i}^{Y'}
\join Y'$, and let $\{A_i\}$ denote the sequence of sets we wish to construct,
namely $\{W_{a_i(s)}^Y \join Y\}$. Thus $\{a_i\}$ will be defined implicitly.

Each $S_i \in \Sigma_2^Y$, (by \citep[Theorem~IV.3.2]{Soa87}, relativized to
$Y$) gives us a computable function $h_i$ such that
\begin{align*}
z \in S_i &\implies \size{W_{h_i(z)}^Y} < \infty, \text{ and } \\
z \not\in S_i &\implies W_{h_i(z)}^Y = \omega.
\end{align*}
Define a $Y$-c.e$.$ set $B_i$ by
\begin{equation*}
B_i^{[z]} := \left\{\pair{x,z}: x \in W_{h_i(z)}^Y\right\}
\end{equation*}
so that for all $z$ and $i$,
\begin{equation}\label{sacksjumpS}
\begin{aligned}
z \in S_i &\implies \size{B_i^{[z]}} < \infty, \text{ and } \\
z \not\in S_i &\implies B_i^{[z]} = \omega.
\end{aligned}
\end{equation}
Given any computable enumeration $\{C_s\}$ of a c.e. set $C$, define, for $s >
0$,
\begin{align*}
\hat{u}(C_s) &:=
\begin{cases}
(\mu x)\: [x \in C_s - C_{s - 1}] & \text{if } C_s - C_{s - 1} \neq \zero, \\
\max (C_s \cup \{s\}) & \text{otherwise};
\end{cases} \\[2ex]
\hat{\Psi}_{e,s}^{C_s}(x) &:=
\begin{cases}
\Psi_{e,s}^{C_s}(x) & \text{if defined and } \psi_{e,s}^{C_s}(x) < \hat{u}(C_s), \\
\diverge & \text{otherwise;}
\end{cases} \\[2ex]
\hat{\psi}_{e,s}^{C_s}(x) &:=
\begin{cases}
\psi_{e,s}^{C_s}(x) &\text{if } \hat{\Psi}_{e,s}^{ C_s}(x) \converge, \\
-1 &\text{otherwise;}
\end{cases} \\[2ex]
T(C) &:= \left\{s : C_s \restr \hat{u}(C_s) = C \restr
\hat{u} (C_s) \right\}, \\
\intertext{and according to Definition~\ref{setrowdefn},}
\omega^{\left[\widetilde{x_k}\right]} &:= \left\{ \bigl< x, \pair{y,k} \bigr> :
y \in \omega \right\}.
\end{align*}
is the $k^\text{th}$ row of the $x^\text{th}$ column.

We now construct $\{A_i\}$ uniformly in $\{S_i\}$.  In order to ensure that $A_i
\not\leq_\T \ijoin_{j \neq i} A_j$ and $(A_i)' \equiv_\T S_i$, we meet the following
requirements:
\begin{align*}
R_\pair{e,i} &: A_i \neq \Psi_e^{\ijoin_{j \neq i} A_j}, \\
P_\pair{e,i} &: A_i^{[e]} =^* B_i^{[e]},
\end{align*}
and we attempt to meet the ``pseudo-requirement''
\begin{align*}
Q_\pair{e,i} : (\exists^\infty s)\ [\hat{\Psi}_{e,s}^{\left( A_i \right)_s}(e) \converge]
\implies \Psi_e^{A_i}(e)\converge.
\end{align*}
$R_\pair{e,i}$ will make $\{A_i\}$ computably independent.  $P_\pair{e,i}$
guarantees that $S_i \leq_\T (A_i)'$ because \[F_i(z) := \lim_x A_i
\left(\pair{x,z}\right)\] exists for all $z$ and is the characteristic function
of $\compliment{S_i}$ by (\ref{sacksjumpS}), and because $F_i \leq_\T (A_i)'$
by the Limit Lemma.  We don't actually meet $Q_\pair{e,i}$, since that would
force $(A_i)'$ to be limit computable in $Y$ and hence $(A_i)' \leq_\T Y'$, but
we do meet $Q_\pair{e,i}$ well enough to ensure $(A_i)'
\leq_\T S_i$.  

\bigskip

\cons Let \[\hat{q}_i(e,s) := \hat{\psi}_{e,s}^{\paren{A_i}_s}(e).\] For each
$i$, fix a computable sequence $\{\left( B_i \right)_s\}_{s \in \omega}$ such
that $B_i = \bigcup_s (B_i)_s$. In this construction, $\hat{q}_i(e,\cdot)$ will
be the restraint function for $Q_\pair{e,i}$, and $r\left(\pair{e,i}, \cdot
\right)$ will be the restraint function from Lemma~\ref{antichain}.  Define
$Y$-computable functions $\hat{Q_i}$ by
\[\hat{Q}_i (e, s) := \max_{j \leq e} \hat{q}_i (j,s).\]

\begin{enumerate}
\item[\textit{Stage $s=0$}.] Let $r\left(\pair{e,i},0\right) =
-1$ for all $\pair{e,i}$.  Set $\left( A_i \right)_0 = \zero \join Y$ for
all~$i$.

\item[\textit{Stage $s+1$.}] Do the following when $s+1$ is an $i^\text{th}$ prime power.  If $s+1$ is not a prime
power, do nothing.

\begin{step}
\item For every even $x \not \in \omega^{[\widetilde{0}]}$ and every $e$ such that
\[x \in \left(B_i^{[e]} \right)_{s+1} \quad \text{and} \quad x > \max \left\{\hat{Q}_i \left(e,s \right)
,\ r \left( \pair{e,i}, s \right) \right\},\] enumerate $x$ into
$\left(A_i^{[e]} \right)_{s+1}$.

\item\label{FMmix}
We modify Stage $s+1$ from Lemma~\ref{antichain} so that the
Lemma~\ref{antichain} strategy happens on the first column of each $A_i$.
Choose the least $e$ such that:
\begin{multline*}
r\left(\pair{e,i},s\right) = -1 \quad \& \quad (\exists \text{ even } x)\Bigl[x
\in\omega^{\left[\widetilde{0_\pair{e,i}}\right]} - \left( A_i \right)_s \\
\& \quad \Psi_{e,s}^{\ijoin_{j \neq i} \left( A_j \right)_s}(x)\converge = 0
\quad \& \quad \left(\forall \pair{z,j} <
\pair{e,i}\right)\left[r\left(\pair{z,j},s\right) < x \right]\Bigr].
\end{multline*}
Here $\omega^{\left[\widetilde{0_\pair{e,i}}\right]}$ is the $\pair{e,i}^\text{th}$ row
of the $0$ column.  If there is no such $e$, then do nothing and go to stage $s+2$. If
$e$ exists, then $R_\pair{e,i}$ acts at stage $s+1$. Perform the following steps.
\begin{enumerate}[(a)]

\item Enumerate $x$ into $\left(A_i^{[\widetilde{0}]}\right)_{s+1}$.

\item Define $r \left(\pair{e,i}, s+1 \right) = s + 1$.

\item For all $\pair{z,j} > \pair{e,i}$, define $r \left(\pair{z,j}, s+1 \right) =
-1$.

\item For all $\pair{z,j} < \pair{e,i}$, define $r\left(\pair{z,j}, s+1 \right) =
r\left(\pair{z,j}, s \right)$.
\end{enumerate}
\end{step}
\end{enumerate}
Finally, we have $A_i = \bigcup_s \paren{A_i}_s$.

\begin{claim} \label{column0FMmet}
For every $\pair{e,i}$, $R_\pair{e,i}$ is satisfied.
\end{claim}
\begin{proof}
Claims~\ref{FMmet} and~\ref{FMmet2} from Lemma~\ref{antichain} each hold in
this construction (with the same proofs) because the restraint function $r$
protects the computations in the same way as before. $\hat{q}_i(e,s)$ does not
at all restrict enumeration into $A_i^{[\widetilde{0}]}$.  Since Step~2 only
enumerates in the first column of each $A_i$, the construction in fact
satisfies the stronger relation
\[A_i^{[\widetilde{0}]} \neq \Psi_e^{\ijoin_{j \neq i} A_j}.\qedhere\]
\end{proof}

\eqref{sacksjumpnoteq} follows immediately from Claim~\ref{column0FMmet}.

\begin{claim}\label{sacksjumpreqP}
For every $\pair{e,i}$, $P_\pair{e,i}$ is satisfied.
\end{claim}

\begin{proof}
Since Step~2 of the construction affects only points in
$A_i^{[\widetilde{0}]}$, the points added in Step~2 have no bearing on the
satisfaction of $P_\pair{e,i}$.

We would like to show that for every $\pair{e,i}$,
\begin{align*}
\hat{L}\left(\pair{e,i}\right) &:= \lim_s \max \left\{\hat{Q}_i (e, s),\
r\left( \pair{e,i}, s \right) \right\}
\end{align*}
is finite, because if this is true, then for all $x
>\hat{L}\left(\pair{e,i}\right)$, $x$ is enumerated into $A_i^{[e]}$ iff $x \in
B_i^{[e]}$.  That is, $A^{[e]} =^* B^{[e]}$.  Claim~\ref{FMmet2} shows that
\[r \left(\pair{e,i}\right) := \lim_s \left[r\left( \pair{e,i}, s \right)\right] <
\infty,\] so it suffices to prove
\begin{equation}\label{sackslim}
\lim_s \hat{Q}_i (e, s) < \infty.
\end{equation}
We won't be able to show this, however, because the limit in \eqref{sackslim}
probably doesn't exist. Nevertheless, if we restrict ourselves to the set of
``true stages,'' we can not only guarantee that the restricted limit exists,
but also that it's finite.  Let
\[T_i := T(A_i).\] $T_i$ is the set of \emph{true stages} in the enumeration $\left\{
\paren{A_i}_s \right\}_{s \in \omega}$.  A true stage $t \in T_i$ guarantees that all
nonzero computations below $\hat{u} \left[\paren{A_i}_t\right]$ are correct.
Thus any \emph{apparent} computation $\hat{\Psi}_{e,t}^{\paren{A_i}_t}(x) = y$
is, in fact, a \emph{true} computation $\Psi_{e}^{A_i}(x) = y$, for any $x$.

$T_i$ is infinite because $\left\{ \paren{A_i}_s \right\}$ is an enumeration of
$A_i$. Indeed, it must happen infinitely often that we enumerate an element $x$
into $\left( A_i \right)_t$ so that $\paren{A_i}_t \restr x = A_i \restr x$.
Furthermore,
\begin{equation}\label{sackslim2}
\lim_{t \in T_i} \hat{Q}_i (e, t)
\end{equation}
is finite because once $\hat{q}_i(j, t)$ converges at some true stage $t$
(which must happen), $\hat{q}_i(j ,\cdot)$ remains unchanged through all
subsequent stages.

Finally, why is it sufficient for the limit \eqref{sackslim2} to be finite
\emph{only on the true stages}?  Because $A_i^{[e]}$ enumerates, in
sufficiently late true stages, all elements in $B_i^{[e]}$ which are greater
than both \eqref{sackslim2} and $r \left( \pair{e,i} \right)$.  This makes
$A_i^{[e]} =^* B_i^{[e]}$,  as desired.
\end{proof}

\begin{claim}\label{sacksjumpleq}
For every $i$, $(A_i)' \leq_\T S_i$.
\end{claim}
\begin{proof}
We determine membership in $(A_i)'$ using an $S_i \join Y' \equiv_\T S_i$ oracle.  Let
\[T_i^e := T\paren{A_i^{[\leq e]}},\] and observe that
\begin{equation}\label{sacksjumpingeq}
e \in (A_i)' \iff \Psi_e^{A_i} (e) \converge \iff (\exists t)\ \left[ t \in T_i^e \quad
\& \quad \hat{\Psi}_{e,t}^{\paren{A_i}_t} (e) \converge \right].
\end{equation}
Indeed, the second $\iff$ in (\ref{sacksjumpingeq}) follows from the fact that
enumerations into $A_i^{[>e]}$ cannot, by definition of $\hat{q}_i$, disturb
settled computations on $e$ with oracle $A_i^{[\leq e]}$.

It remains to show that the main predicate in~(\ref{sacksjumpingeq}) is $Y$-computable
and can be constructed from $S_i$, because then membership in $A_i$ is determined by a
formula decidable in $Y'$.  Note that $A_i^{[\leq e]}$ is computable because $A_i^{[\leq
e]} =^*B_i^{[\leq e]}$, and because $B$ is piecewise computable.  Hence $T_i^e \equiv_\T
A_i^{[\leq e]}$ is computable.  Moreover, the index for $T_i^e$ can be found from the
following series of \emph{uniform} reductions:
\[T_i^e \leq_\T A_i^{[\leq e]} \leq_\T B_i^{[\leq e]} \join \left( Y' \join \{i\} \right) \leq_\T S_i \join \{i\},\]
where $Y' \join \{i\}$ is used to compute the finite, $Y$-c.e$.$ set
$A_i^{[\widetilde{0}]} \cap \omega^{[\leq e]}$. This means that the computable predicate
$t \in T_i^e$ can be constructed from $S_i \join \{i\}$.  Furthermore,
$\hat{\Psi}_{e,t}^{\paren{A_i}_t} (e) \converge$ is a $Y$-computable predicate, which
makes
\[(\exists t) \left[ t \in T_i^e \quad \& \quad \hat{\Psi}_{e,t}^{\paren{A_i}_t} (e) \converge \right] \]
from \eqref{sacksjumpingeq} decidable in $Y'$.  Thus \[\paren{A_i}' \leq_\T
(S_i \mathrel{\join} \{i\}) \join Y' \equiv_\T S_i. \qedhere\]
\end{proof}

Claim~\ref{sacksjumpreqP} and Claim~\ref{sacksjumpleq} together now prove
\eqref{sacksjumpeq}.
\end{proof}

\begin{lemma}\label{2infrecursion}
Let $\{a_k\}$ be a computable sequence of partial computable functions.  Then there
exists a computable sequence $\{x_k\}$ such that for all $Y \subseteq \omega$,
\begin{gather*}
\Psi^Y_{a_0(x)} = \Psi_{x_0}^Y \\
\Psi^Y_{a_1(x)} = \Psi_{x_1}^Y \\
\vdots
\end{gather*}
where $x$ is such that $\phe_x(k) = x_k$ for all $k$.
\end{lemma}
\begin{proof}
Let $\{a_k\}$ be a computable sequence of p.c$.$ functions. Define a partial
computable function $f$ by
\[f(z,k) := a_k(z).\]
By the $s$-$m$-$n$ Theorem, there is a computable function $g$ such that
\[\phe_{g(z)}(k) = f(z,k).\]
By the Recursion Theorem, there exists a fixed point $x$ satisfying
\[\phe_x(k) = \phe_{g(x)}(k) = f(x,k).\]
Thus
\[a_k(x) = f(x,k) = \phe_x(k) := x_k.\]
It follows that for all $Y \subseteq \omega$,
\[\Psi_{a_k(x)}^Y = \Psi_{x_k}^Y.\qedhere\]
\end{proof}

\section{Properties of $\MIN^{\T^{(\omega)}}$} \label{properties of MINT}

\begin{defn}
Let $f$ and $g$ be total functions, and let $A = \{ a_0 < a_1 < \dotsb \}$ be a set.
\begin{defenum}
\item $f$ \emph{majorizes} $g$ if $(\forall n)\: [ f(n) > g(n)]$.

\item $f$ \emph{dominates} $g$ if $(\forall^\infty n)\: [ f(n) > g(n)]$, where
$\forall^\infty$ means ``for all but finitely many.''

\item The function $p_A(n) := a_n$ is called the \emph{principal function} of $A$.

\item A function $f$ \emph{majorizes} a set $A$ if $(\forall n)\: [f(n) > p_A(n)].$

\item Let $\mathbf{a}$ be a Turing degree.  A set $A$ is \emph{$\mathbf{a}$-dominated (resp$.$ $\mathbf{a}$-majorized)} if there
exists an $\mathbf{a}$-computable function $f$ which dominates (resp$.$
majorizes) $A$.
\end{defenum}
\end{defn}

\begin{thm}[\cite{KMU55}, \cite{Soa08}]\label{hyperimmuneiffdominated}
An infinite set $A$ is hyperimmune iff $A$ is not $\0$-dominated.
\end{thm}

We obtain the following paradoxical result:

\begin{thm}[peak hierarchy]\label{peak}
$\MIN^{\T^{(\omega)}}$
\begin{thmenum}
\item is infinite, \label{mintomegaisinf2}

\item contains no infinite arithmetic sets, and \label{peaknotcontain}

\item is not hyperimmune. \label{peakhyperimmune}
\end{thmenum}
\end{thm}

\begin{proof}

\begin{proof}[(\ref{mintomegaisinf2})]
Theorem~\ref{intseq} provides a denumerable list of distinct
$\equiv_{\T^{(\omega)}}$ classes.
\end{proof}

\begin{proof}[(\ref{peaknotcontain})]
Follows from Corollary~\ref{mintnimmune}, because $\MIN^{\T^{(\omega)}} \subseteq
\MIN^{\T^{(n)}}$ for every $n$.
\end{proof}

\begin{proof}[(\ref{peakhyperimmune})]
We verify that $\MIN^{\T^{(\omega)}}$ gets majorized.  Let $\{x_k\}$ be as in
Theorem~\ref{intseq}. Then for all $n$ and $i \neq j$,
\[W_i \not\equiv_{\T^{(n)}} W_j.\]
Without loss of generality, $x_0 < x_1 < \dotsb$ since $\{x_k\}$ is computable. Define
the computable function
\begin{align*}
f(0) &:= x_1 \\
f(n+1) &:= x_{[2f(n)]},
\end{align*}
 and let $p$ be the principal function
of $\MIN^{\T^{(\omega)}}$. Note that $f(0) > 0 = p(0)$, and assume for the purposes of
induction that $f(n) > p(n)$.  Note that
\[ p(n) \leq x_{p(n)} < x_{f(n)} < x_{f(n) + 1} < \dotsb < x_{2f(n)} = f(n+1), \]
so at least $f(n)$ $x_k$'s lie strictly between $p(n)$ and $f(n+1)$, namely
\[\{x_{f(n)}, x_{f(n)+1}, \dotsc, x_{2f(n) - 1} \}.\]
Hence, at least $f(n)$ distinct $\equiv_{\T^{(\omega)}}$-equivalence classes
are represented by indices strictly between $p(n)$ and $f(n+1)$.  Since less
than $f(n)$ classes are represented in indices up to $p(n)$, there necessarily
must be a new $\equiv_{\T^{(\omega)}}$-class introduced strictly between $p(n)$
and $f(n+1)$.  This forces $p(n+1) < f(n+1)$.  Hence $f$ majorizes
$\MIN^{\T^{(\omega)}}$.   The result now follows immediately from
Theorem~\ref{hyperimmuneiffdominated}.
\end{proof}
\end{proof}

Consequently, the other $\MIN$-sets in this thesis share properties
(\ref{mintomegaisinf2}) and (\ref{peakhyperimmune}):
\begin{cor} \label{containingMINTomega}
Every set containing $\MIN^{\T^{(\omega)}}$ is infinite but not hyperimmune.
\end{cor}

\begin{rem}
If we only wanted to prove that $\MIN^\T$ is $\0$-dominated, we could have simplified the
proof of Theorem~\ref{peak} by omitting Lemmas~\ref{sacksjumpseq}
and~\ref{2infrecursion}, thereby avoiding infinite injury.
\end{rem}

$\zero^{(\omega)}$ is another familiar set which is hyperarithmetic and
$\0$-dominated. However, unlike $\MIN^{\T^{(\omega)}}$, $\zero^{(\omega)}$
contains a copy of $\zero'$.  This means that $\zero^{(\omega)}$ is not at all
immune.

\section{A strange ``cutting'' set}
Lusin once constructed a set of reals which neither contains nor is disjoint
from any perfect set  \citep{Lus21}, \citep[Theorem 2.25]{MW85}.  By modifying
Lusin's construction and gently expanding $\MIN^{\T^{(\omega)}}$, we obtain an
analogous construction for the arithmetic hierarchy which is remarkably
well-behaved.

\begin{cor} \label{cuttingset}
There exists a set $X \supseteq \MIN^{\T^{(\omega)}}$ such that $X$:
\begin{thmenum}
\item contains no infinite arithmetic sets, \label{cuttingarithmetic}

\item is not disjoint from any infinite arithmetic set, and
\label{cuttingdisjoint}

\item is $\0$-majorized. \label{cuttingmajorized}
\end{thmenum}
\end{cor}

\begin{proof}
We simultaneously enumerate disjoint sets $X$ and $Y$ so that both $X$ and $Y$
intersect every infinite arithmetic set.  Let $A_0, A_1, A_2, \dotsc$ be an
enumeration of the arithmetic sets.  Let
\begin{align*}
X_0 &:= \MIN^{\T^{{(\omega)}}}, \\
Y_0 &:= \emptyset.
\end{align*}
Now assume that $X_n$ and $Y_n$ have already been constructed, with
\begin{align*}
X_n &=^* \MIN^{\T^{(\omega)}}, \\
Y_n &=^* \emptyset,
\end{align*}
$X_n \intersect Y_n = \emptyset$, and for all $k < n$,
\begin{align*}
X_n \intersect A_k &\neq \emptyset, \\
Y_n \intersect A_k &\neq \emptyset.
\end{align*}

Since $\MIN^{\T^{(\omega)}}$ does not contain any infinite arithmetic sets
(Theorem~\ref{peak}), it follows that $A_n \intersect
\compliment{\MIN^{\T^{(\omega)}}}$ is infinite, which means that
\begin{equation*}
R_n := \paren{A_n \intersect \compliment{X_n}} - Y_n
\end{equation*}
is infinite.  Enumerate the least element in $R_n$ into $X_n$, and call this
expanded set $X_{n+1}$.  Now
\begin{equation*}
S_n := \paren{A_n \intersect \compliment{Y_n}} - X_{n+1}
\end{equation*}
is also infinite.  Enumerate the least element in $S_n$ into $Y_n$. Finally,
let
\begin{align*}
X &:= \Union_{n \in \omega} X_n, \\
Y &:= \Union_{n \in \omega} Y_n.
\end{align*}
This concludes the construction.

It is clear that $X_{n+1}$ satisfies the same inductive hypotheses as $X_n$.
Consequently $X$ and $Y$ are disjoint, because otherwise a common element would
have been introduced after finitely many stages.  Furthermore,
\begin{align*}
X_n \intersect A_n &\neq \emptyset, \\
Y_n \intersect A_n &\neq \emptyset,
\end{align*}
for every $n$, which proves (\ref{cuttingdisjoint}) and
(\ref{cuttingarithmetic}).  Finally, (\ref{cuttingmajorized}) follows because
$\MIN^{\T^{(\omega)}}$ is an infinite subset of $X$, by
Corollary~\ref{containingMINTomega} and Theorem~\ref{hyperimmuneiffdominated}.
\end{proof}

\begin{rem}
It is straightforward to construct a set which satisfies properties
(\ref{cuttingarithmetic}) and (\ref{cuttingdisjoint}) from
Corollary~\ref{cuttingset} but not (\ref{cuttingmajorized}).  In the
construction of Corollary~\ref{cuttingset}, one can achieve this by first
replacing $X_0 = \MIN^{\T^{(\omega)}}$ with $X_0 = \emptyset$, and then, during
enumeration into $X_n$, by choosing an element so that $X(n) > A_n(n)$ rather
than just enumerating the least element in $R_n$.
\end{rem}

%% file: open.tex
\chapter{Open problems} \label{open problems}

\section{Truth table degrees}

Meyer's original question from 1972 remains open: is $\fMIN \equiv_\Tt \zero''$
\citep{Mey72}?   A reduction $\fMIN \geq_\bT \zero''$ would suffice to show
$\fMIN \equiv_\Tt \zero''$, if it were the case that $\zero' \leq_\Tt \MIN$
\citep[Section 8]{Sch98}.  Similarly, Schaefer asks, is $\fR \equiv_\Tt \zero'$
\citep{Sch98}?  The fact that we know $\RAND_\phe \equiv_\Tt \zero'$ for any
Kolmogorov numbering $\phe$
(Theorem~\ref{numberingsRANDfR}\nlb(\ref{RANDttequiv0'})) but we don't know the
truth-table degree of its cousin $\fR$ indicates that there is still much to
learn about similarites between randomness and minimal indices.

\section{Is $\MIN^\T \equiv_\T \zero''''$?}

We conjecture that Corollary~\ref{thickpheTuringcomplete} does not hold for
arbitrary G\"{o}del numberings.  In particular, we conjecture that that
Corollary~\ref{bestTuring} is optimal in the following sense:
\begin{conj} Let $n\geq 0$. \label{Tincompletenessconjecture}
\begin{thmenum}
\item There exists a G\"{o}del numbering $\phe$ such that $\MIN^*_\phe \not
\geq_\T \zero'$.

\item There exists a G\"{o}del numbering $\phe$ such that $\MIN^\m_\phe \join \zero' \not
\geq_\T \zero''$.

\item There exists a G\"{o}del numbering $\phe$ such that $\MIN^{\T^{(n)}}_\phe
\join \zero^{(n+1)} \not \geq_\T \zero^{(n+2)}$.
\end{thmenum}
\end{conj}
Even showing $\MIN^\m_\phe \not \geq_\T \zero''$ or $\MIN^\T_\phe \not \geq_\T
\zero''$ for some G\"{o}del numbering $\phe$ would be enough to resolve the
Turing degree of $\MIN^\m$ or $\MIN^\T$.

All of the initial information in a $=^*$ set can be faulty \citep{Sch98}, so
intuitively one needs a halting set oracle to extract useful information from
$\MIN^*$. Similarly, $\MIN^\m$ and $\MIN^\T$ presume knowledge of total
functions, making $\zero'' \equiv_\T \TOT$ undecidable relative to these sets.
The difficulty in constructing the necessary numberings for
Conjecture~\ref{Tincompletenessconjecture} is revealed by considering a simpler
problem where we try to find \emph{any} $A \in \Sigma_4$ satisfying:
\begin{gather*}
A \join \zero'' \equiv_\T \zero'''', \\
A \not \geq_\T \zero''.
\end{gather*}
A set $A$ with these properties can be constructed using two iterations of the
Sacks Jump Theorem (Theorem~\ref{sacksjumptheorem}), however this is already a
nontrivial, infinite injury construction.  Making this construction work with
$A = \MIN_\phe^\T$ for some \emph{G\"{o}del} numbering $\phe$ can only be more
complicated.

If Conjecture~\ref{Tincompletenessconjecture} holds, then spectral sets are
(possibly the first) natural examples of sets which are not Turing equivalent
to any of the canonical $\Sigma_n$-complete sets.  If
Conjecture~\ref{Tincompletenessconjecture} fails, then spectral sets are a new
and remarkable characterizations of the Turing degrees $\0'$, $\0''$, $\0'''$,
$\dotsc$.

One approach to solving the $\MIN^*$ problem is to look first at the related
problem of $\MIN^\m$.  This approach is promising because it has not received
much attention.  It is also promising for mathematical reasons.  We now sing
praises of $\MIN^\m$.  If indeed $\MIN^\m \join \zero'' \equiv_\T \zero'''$ and
$\MIN^* \join \zero' \equiv_\T \zero'''$ are both optimal results (in the sense
of Conjecture~\ref{Tincompletenessconjecture}), then it seems easy to find a
numbering $\phe$ in which $\MIN^\m_\phe$ avoids (merely) the cone of degrees
above $\zero''$, when compared to the (daunting) task of forcing $\MIN^*_\phe$
to avoid the the cone above $\zero'$.  The second reason to take up $\MIN^\m$
is for the elegance and brevity of results given in this thesis which are
unique to $\MIN^\m$. The Generalized Fixed Point
Theorem~\ref{JLSSmain}\nlb(\ref{JLSSmainm}) immediately gives optimal immunity
for $\MIN^\m$ (see $\Pi_3$-Separation Theorem \ref{pi3separation}), and our
purported optimal result for the Turing degree of $\MIN^\m$,
Lemma~\ref{downonelevel}\nlb(\ref{MINmjoin0''is0'''}), follows directly from
the $\equiv_\m$-Completeness
Criterion~\ref{jockuschcompleteness}\nlb(\ref{jockuschcompletenessm}). Finally,
we have a satisfying proof of the fact that $\MIN^\m_\psi \equiv_\Tt \zero'''$
for some Kolmogorov numbering $\psi$ (Theorem~\ref{thickpheTuringcomplete}).
This same argument finds only a Turing degree for $\MIN^{\T^{(n)}}_\psi$.

\section{$\MIN$ vs. $\fMIN$}
We know that $\MIN \equiv_\T \zero'' \equiv_\T \fMIN$.  What can be said about
stronger reductions?  For example is it true, in general, that $\MIN_\phe
\equiv_\btt \fMIN_\phe$?  We know that there exists a Kolmogorov numbering
$\psi$ such that $\MIN_\psi \equiv_\Tt \fMIN_\psi$
(Theorem~\ref{thickpheTuringcomplete}\nlb(\ref{thickpheTuringcompleteMIN})),
and for any numbering $\phe$, there is a G\"{o}del numbering $\psi$ such that
$\MIN_\phe \not\equiv_\btt \fMIN_\psi$ (by Theorem~\ref{MINbttMIN}).  But do
there exists any G\"{o}del numberings $\phe$ and $\psi$ such that $\MIN_\phe
\equiv_\btt \fMIN_\psi$? Given a G\"{o}del numbering $\phe$, does there always
exist a G\"{o}del numbering $\psi$ such that $\MIN_\phe \equiv_\Tt \fMIN_\psi$?

\section{$\Pi_n$-immunity}
In Chapter~\ref{immunitychapter}, we prove optimal results with respect to
$\Sigma_n$-immunity. What about $\Pi_n$-immunity?  In particular, is $\MIN$
$\Pi_1$-immune? Is $\MIN^*$ $\Pi_2$-immune?  Is $\MIN^{\T^{(n)}}$
$\Pi_{n+3}$-immune?  We do have an optimal $\Pi_n$-immunity result for
$\MIN^\m$ ($\Pi_3$-Separation Theorem~\ref{pi3separation}), however this
argument does not generalize to other spectral sets.

Is there a weak form of the Arslanov Completeness Criterion which is equivalent
to immunity for $\MIN$-sets?  Does $\MIN$ contain a $\Delta_2$ spectral set, as
$\fMIN$ does?  Does there exist a direct reduction from $\fMIN$ to $\fR$ which
does not go through the halting set?

\section{A question of Friedman}
$\MIN^*$ and $\MIN^\T$ are not the only sets with short descriptions whose
Turing degree remains elusive.  Consider the set
\[A := \{n :  (\exists j < n)\: [\max(W_j) = n] \}.\]
Friedman asks \citep{Sch06}, what is the complexity of $A$?  A straightforward
argument shows that $A \leq_\T \zero'$, however it is surprising that we do not
know whether or not $A \in \Sigma_1$, $A \in \Pi_1$, or $A' \equiv_\T \zero'$.

On a related note, we visit a set reminiscent of the Kolmogorov random strings.
For any numbering $\phe$, let
\[\sR_\phe := \{x : (\forall j < x)\: [\phe_j (0) \neq x ] \}.\]
What are the possible degrees for $\sR_\phe$ when $\phe$ is a Kolmogorov
numbering?  Can $\sR_\phe$ be finite?

\section{Intermediate degrees}
\begin{prop}\label{joinprop}
Let $A_0, A_1, \dotsc$ be a sequence of sets, and let $I$ be a computable set.
Then
\[\ijoin_{i \in I} (A_i)' \leq_\T \left( \ijoin_{i \in I} A_i \right)'\]
\end{prop}
\begin{proof}
For all $k \in I$, \[A_k \leq_\T \ijoin_{i \in I} A_i. \] It then follows from
the Jump Theorem \citep{Soa87} that
\[(A_k)' \leq_\T \left( \ijoin_{i \in I} A_i \right)',\]
so
\[\ijoin_{i \in I} (A_i)' \leq_\T \left( \ijoin_{i \in I} A_i \right)' \join I \equiv_\T \left( \ijoin_{i \in I} A_i
\right)'. \qedhere\]
\end{proof}

In light of Theorem~\ref{intseq}, where we found a computable sequence
$\{x_k\}$ satisfying for all $n$ and $i$
\[\left( W_{x_i} \right)^{(n)} \not\leq_\T \ijoin_{j \neq i} \left( W_{x_j} \right)^{(n)},\]
Proposition~\ref{joinprop} leaves us with the burning question of whether or
not there is a computable sequence $\{z_k\}$ satisfying for all $n$ and $i$,
\[\left( W_{z_i} \right)^{(n)} \not\leq_\T \left( \ijoin_{j \neq i} W_{z_j} \right)^{(n)}.\]
It would be sufficient to show that for any computable $I \subseteq \omega$,
the indices $\{a_i\}$ from Lemma~\ref{sacksjumpseq} satisfy the additional
condition,
\begin{equation*}\label{sacksseqlance}
\left( \ijoin_{i \in I} \paren{W_{a_i(s)}^Y \join Y} \right)' \equiv_\T
\ijoin_{i \in I} \left( W_{s_i}^{Y'} \join Y' \right).
\end{equation*}
because then we can pull the jump operator outside of the join, as would be
needed in the proof of Theorem~\ref{intseq}.

\section{Other complexity measures and variants}
What can be said about approximability, autoreducibility, and size-minimal
indices of spectral sets?  These questions (and some answers to them) appear in
Schaefer's paper \citep{Sch98}.

\begin{defn}
A set $A$ is \emph{autoreducible} if, for all $x$, one can decide whether $x
\in A$ by querying only elements in $A - \{x\}$.
\end{defn}

Schaefer showed that $\fMIN$ is autoreducible \citep{Sch98}.  Are there other
spectral sets which are autoreducible?

\begin{defn}
\begin{defenum}
\item A set $A$ is \emph{$(1,k)$-computable} if there exists a computable
function $f$ such that for every set $X \subseteq \omega$ with $\size{X} = k$,
there is some $x \in X$ such that $f(x) = \chi_A(x)$.

\item A set $A$ is \emph{approximable} if $A$ is $(1,k)$-computable for some
$k$.
\end{defenum}
\end{defn}
$\fMIN$ is not $(1,2)$-computable, and there exists a G\"{o}del numbering
$\phe$ such that $\fMIN_\phe$ is not approximable \citep{Sch98}.  Does there
exists a G\"{o}del numbering $\psi$ such that $\MIN_\psi$ \textit{is}
approximable?

\begin{defn}
For a G\"{o}del numbering $\phe$ and a (total) size function $s$, define
\[\fMIN_{\phe,s} = \{e : (\forall i)\: [s(i) < s(e) \implies \phe_i \neq
\phe_e]\},\] to be the set of \emph{size-minimal indices of $\phe$}.
\end{defn}
In contrast to the $\MIN$-sets in this paper, there is a computable size
function $s$ (independent of the G\"{o}del numbering $\phe$) such that
$\fMIN_{\phe,s}$ is hyperimmune \citep{Sch98}.  It is an open problem to
determine whether $\fMIN_{\phe,s} \equiv_\T \zero''$ whenever $\phe$ is a
G\"{o}del numbering and $s$ is computable.

\section{Almost thickness} \label{almost thickness}
Ken Harris made the following observation, generalizing a familiar
representation for $\Sigma_3$ sets:

\begin{thm}[$\Sigma_n$-Representation, \citet{Har06}]
\label{harrisrepresentation}
 If $A \in \Sigma_n$, then there is a computable function $g$ such that
\begin{align*}
x \in A &\iff (\forall^\infty y_1)\, (\forall y_2)\, (\forall^\infty y_3)
\dotso \left[W_{g(x, \overline{y})} = \omega \right], \\
x \in \compliment{A} &\iff (\forall y_1)\, (\forall^\infty y_2)\, (\forall y_3)
\dotso \left[W_{g(x,\overline{y})} \text{ is finite} \right],
\end{align*}
where $\dotso$ denotes the remaining of the $n-3$ quantifiers.

Similarly, if $A \in \Pi_n$, then there is a computable function $g$ such that
\begin{align*}
x \in A &\iff (\forall y_1)\, (\forall^\infty y_2)\, (\forall y_3)
\dotso \left[W_{g(x, \overline{y})} = \omega \right], \\
x \in \compliment{A} &\iff (\forall^\infty y_1)\, (\forall y_2)\,
(\forall^\infty y_3) \dotso \left[W_{g(x,\overline{y})} \text{ is finite}
\right].
\end{align*}
\end{thm}

Here $\forall^\infty x$ means ``for all but finitely many x.''  This makes
$\aev$ a natural complement to the $\thick$ operator:
\begin{defn}
Let $\equiv_\alpha$ be an equivalence relation on sets.  Then
\[A \equiv_{\aev \equiv_\alpha} B \enskip \iff \enskip (\forall^\infty n) \left[ A^{[n]}
\equiv_\alpha B^{[n]} \right].\]
\end{defn}

Using the Representation Theorem~\ref{harrisrepresentation} and the methods
from Section~\ref{thickarithmetics}, we arrive at the following observations:.

\begin{prop}
\begin{thmenum}
\item $\MIN^{\aev =} \in \Pi_3 - \Sigma_3$,

\item $\MIN^{\aev *} \in \Pi_5 - \Sigma_5$,

\item $\MIN^{\aev \m} \in \Pi_5 - \Sigma_5$,

\item $\MIN^{\aev{\T^{(n)}}} \in \Pi_{n+6} - \Sigma_{n+6}$.
\end{thmenum}
\end{prop}

We remark that $\aev *$ is an example of an equivalence relation for which
$$\MIN^{\T^{(n)}} \not\supseteq \MIN^{\aev *} \not\supseteq \MIN^{\T^{(n+1)}}$$
for all $n$, however $\MIN^\aev$ sets are not hyperimmune.  What can be said
about the immunity of $\MIN^\aev$-sets?

\begin{prop}
$\MIN^{\aev =}$ is not $\Sigma_3$-immune. \label{ae=notimmune}
\end{prop}

\begin{proof}
Let
\[\omegaINF := \{ e: (\forall N)\, (\exists y \geq N)\, (\exists x)\, [ x \in
W_e^{[y]} ] \}. \] Similar to the proof of
Theorem~\ref{pi3separation}\nlb(\ref{pi3separation*}), let
\begin{align*}
P_k &:= \{ \pair{x,y} : y \text{ is a $k^\text{th}$ prime power} \}, \\
\intertext{and set}
A_k  &:= \left\{e : (\exists N)\, (\forall y \geq N) \left[W_e^{[y]} \subseteq P_k^{[y]} \right]\right\} \intersect \omegaINF, \\
A &:= \{e : (\exists k)\, (\forall j < e) \left[e \in A_k \enskip \& \enskip j
\not\in A_k \right]\}.
\end{align*}
Now $A_k \in \Delta_3$, because $\omegaINF \in \Pi_2$ and \[\left\{ e :
W_e^{[y]} \subseteq P_k^{[y]} \right\}\] is in $\Pi_1$.  Hence $A \in
\Sigma_3$.  $A$ is infinite because $A$ contains a member from each $A_k$, and
the $A_k$'s are pairwise disjoint.  Finally, $A \subseteq \MIN^{\aev =}$
because $e \in A$ implies $e \in A_k$ for some $k$, and any set which is equal
to $W_e$ must eventually be contained in $A_k$ (for sufficiently high rows). So
$\MIN^{\aev=}$ contains an infinite, $\Sigma_3$ subset.
\end{proof}

\begin{conj} For $n \geq 0$, \label{aenotimmune}
\begin{thmenum}
\item $\MIN^{\aev *}$ is not $\Sigma_4$-immune. \label{ae*notimmune}

\item $\MIN^{\aev \m}$ is not $\Sigma_5$-immune. \label{aemnotimmune}

\item $\MIN^{\aev{\T^{(n)}}}$ is not $\Sigma_{n+5}$-immune. \label{aeTnotimmune}
\end{thmenum}
\end{conj}

While it is clear that $\thick$ and $\aev$ can be combined iteratively to
obtain more equivalence relations (and hence more open questions), the author
currently considers this direction somewhat esoteric.

\section{Ershov hierarchy}

\begin{defn}
A set $X$ is \emph{d.c.e$.$ (difference of c.e$.$ sets)} if there exists c.e$.$
sets $A$ and $B$ such that $X = A - B$.  More generally, $X$ is
\emph{$n$-c.e$.$} if there exists a sequence of c.e$.$ sets $A_1, A_2, \dotsc,
A_n$ such that \[X = A_1 - A_2 \union A_3 - \dotsc \union A_n,\] where order of
operations in left to right.
\end{defn}

The d.c.e$.$ sets can be enumerated as pairs of indices for c.e$.$ sets, so it
is natural to consider minimal indices for d.c.e$.$ sets \citep{Har06}.  Let
$V_0, V_1, \dotsc$ be an enumeration of the d.c.e$.$ sets, and let
\[\text{2-$\MIN$} := \{e : (\forall j < e)\: [V_j \neq V_e] \}.\]
What is the Turing degree of 2-$\MIN$?  A similar question can be asked about
the $n$-c.e$.$ sets, whose indices are also enumerable.

\section{Polynomial-time}
Computational complexity intersects nontrivially with Kolmogorov complexity
\citep{For04}, so it is natural to ask what applications computational
complexity has in the generalized world of minimal indices (and vice-versa). We
examine, for example, a familiar notion from resource-bounded complexity within
the context of minimal indices.
\begin{defn} Let $A$ and $B$ be sets.  We write $A \leq_\p B$ if there exists
a computable algorithm $f$, running in time polynomial in input length, which
satisfies
\[x \in A \bigiff f(x) \in B\]
for all $x$.  If $A \leq_\p B$ and $B \leq_\p A$, we write $A \equiv_\p B$. Now
\[\MIN^\p := \{e : (\forall j < e)\: [W_j \not\equiv_\p W_e]\}.\]
\end{defn}
It is immediate that $\MIN^\p$ is $\Sigma_2$-immune, as $\MIN^* \supseteq
\MIN^\p$ (Theorem~\ref{RMINMIN*immune}\nlb(\ref{MIN*immune})), however the
Turing degree of $\MIN^\p$ is not known.